\documentclass[pdftex]{article}
\usepackage[top=1truein,bottom=1truein,left=1truein,right=1truein]{geometry}

\usepackage{graphicx}
\graphicspath{{./fig}}

\usepackage{authblk}
\usepackage{abstract}
\usepackage{amsmath}
\usepackage{mathtools}
\usepackage{amssymb}
\usepackage{amsfonts}

\usepackage{hyperref}
\usepackage{url}
\usepackage{orcidlink}

\usepackage{here}
\usepackage{tikz}
\usepackage{bm}
\usepackage{setspace}
\usepackage{enumerate}
\usepackage{arydshln}
\usepackage{caption}
\usepackage{color}
\usepackage{appendix}
\usepackage{lineno}
\usepackage{ulem}

\hypersetup{%
  colorlinks=true,
  linkcolor=blue,
  citecolor=blue,
  urlcolor=cyan
}

\DeclareMathOperator*{\argmin}{argmin}
\begin{document}

\title{\textbf{A pressure-projection formulation in a least-squares meshfree method for the incompressible Navier--Stokes equations using a staggered-variable arrangement}}

\author[ab]{Takeharu Matsuda\,\orcidlink{0009-0007-6263-2766}}
\author[b]{Satoshi Ii\,\orcidlink{0000-0002-5428-5385}\footnote{Corresponding author.\\ \textit{Email address:} ii.s.148c@m.isct.ac.jp (Satoshi Ii\,\orcidlink{0000-0002-5428-5385})}}

\affil[a]{
    {\small
      \textit{
        Department of Mechanical Engineering, Graduate School of Engineering, Kobe University, 1-1 Rokkodai-cho, Nada-ku, Kobe 657-8501, Japan
      }
    }
}
\affil[b]{
  {\small
    \textit{
      Department of Mechanical Engineering, School of Engineering, Institute of Science Tokyo, 2-12-1 Ookayama, Meguro-ku, Tokyo 152-8550, Japan
    }
  }
}

\date{}
\maketitle

\renewenvironment{abstract}
 {\quotation\small\noindent\rule{\linewidth}{.5pt}\par\smallskip
  {\centering\bfseries\abstractname\par}\medskip}
 {\par\noindent\rule{\linewidth}{.5pt}\endquotation}

\begin{abstract}
  Incompressible flow solvers based on strong-form meshfree methods represent arbitrary geometries without the need for a global mesh system. However, their local evaluations make it difficult to satisfy incompressibility at the discrete level. Moreover, the collocated arrangement of velocity and pressure variables tends to induce a zero-energy mode, leading to decoupling between two variables. In projection-based approaches, a spatial discretization scheme based on a conventional node-based moving least-squares method for the pressure causes inconsistency between discrete operators on both sides of the Poisson equation. Thus, a solenoidal velocity field cannot be ensured numerically. In this study, a numerical method for the incompressible Navier--Stokes equations is developed by introducing a local primal--dual grid into the mesh-constrained discrete point method, enabling consistent discrete operators. The constructed \textit{virtual} dual cell is defined solely from local node connectivity, and thus the method retains its meshfree nature. To consistently couple velocity and pressure variables under primal--dual arrangement, a numerical technique called time evolution converting, which evaluates interface quantities by considering the temporal evolution of flow, is used to update the velocity at cell interfaces. For numerical validation, a linear acoustic equation is solved to confirm the effectiveness of the staggered-variable arrangement based on the local primal--dual grid. Then, the incompressible Navier--Stokes equations are solved, and the proposed method is demonstrated to reduce the discrete velocity divergence to a level determined by the residual of the pressure Poisson equation, achieve the expected spatial convergence order, and accurately reproduce flow features over a wide range of Reynolds numbers.

  \vskip.5\baselineskip

\end{abstract}


\section*{Nomenclature}
\begin{center}
  \scriptsize
  \renewcommand{\arraystretch}{1.0}
  \begin{tabular}{@{}p{0.18\linewidth}p{0.57\linewidth}p{0.15\linewidth}@{}}
    \hline
    Symbol & Description & SI unit \\
    \hline
    \multicolumn{3}{@{}l}{\textit{Roman symbols}} \\
    $A$ & Amplitude of cylinder oscillation & $\mathrm{m}$ \\
    $\mathbf{A}$ & Radial projection of the convective flux & $\mathrm{m^{3}\,s^{-2}}$ \\
    $c$ & Acoustic speed & $\mathrm{m\,s^{-1}}$ \\
    $C$ & Dimensionless coefficient or number & -- \\
    $D$ & Cylinder diameter & $\mathrm{m}$ \\
    $\mathbf{D}$ & Radial projection of the velocity gradient & $\mathrm{m\,s^{-1}}$ \\
    $\mathbf{e}$ & Edge vector & $\mathrm{m}$ \\
    $f$ & Frequency of cylinder oscillation & $\mathrm{Hz}$ \\
    $H$ & Channel width & $\mathrm{m}$ \\
    $KC$ & Keulegan--Carpenter number & -- \\
    $l_0$ & DP spacing (background-mesh width) & $\mathrm{m}$ \\
    $L$ & Characteristic or domain length & $\mathrm{m}$ \\
    $L_1$, $L_2$, $L_\infty$ & Numerical error norms & $\mathrm{m\,s^{-1}}$ or $\mathrm{Pa}$ \\
    $N$ & Number of DPs in each coordinate direction & -- \\
    $\mathbf{n}$ & Outward unit normal vector & -- \\
    $p$ & Pressure (gauge pressure in the linear acoustic problem) & $\mathrm{Pa}$ \\
    $q$ & Boundary-normal velocity component & $\mathrm{m\,s^{-1}}$ \\
    $R$ & Radius & $\mathrm{m}$ \\
    $Re$ & Reynolds number & -- \\
    $St$ & Strouhal number & -- \\
    $t$ & Time & $\mathrm{s}$ \\
    $T$ & Period of cylinder oscillation & $\mathrm{s}$ \\
    $\mathbf{u}=(u,v)^{\top}$ & Fluid velocity & $\mathrm{m\,s^{-1}}$ \\
    $U$, $\mathbf{U}$ & Characteristic speed or velocity radial component & $\mathrm{m\,s^{-1}}$ or $\mathrm{m^{2}\,s^{-1}}$ \\
    $\mathbf{x}$ & Position vector & $\mathrm{m}$ \\
    \hline
    \multicolumn{3}{@{}l}{\textit{Greek symbols}} \\
    $\alpha$ & Damping or randomization parameter & $\mathrm{s^{-1}}$ or -- \\
    $\beta$ & TEC parameter & -- \\
    $\Delta t$ & Time interval & $\mathrm{s}$ \\
    $\delta\mathbf{x}$ & Random displacement vector of a DP & $\mathrm{m}$ \\
    $\delta P$ & Radial projection of the pressure gradient & $\mathrm{Pa}$ \\
    $\eta$ & Fluid viscosity & $\mathrm{Pa\,s}$ \\
    $\rho$ & Fluid density & $\mathrm{kg\,m^{-3}}$ \\
    $\boldsymbol{\tau}$ & Viscous stress tensor & $\mathrm{Pa}$ \\
    $\theta$ & Edge-interpolation parameter & -- \\
    \hline
    \multicolumn{3}{@{}l}{\textit{Subscripts}} \\
    $0$ & Reference or characteristic quantity & -- \\
    $\mathrm{c}$ & Center & -- \\
    $\mathrm{D}$ & Drag & -- \\
    $\mathrm{L}$ & Lift & -- \\
    $\mathrm{diff}$ & Diffusion & -- \\
    $\mathrm{cyl}$ & Cylinder & -- \\
    $\mathrm{damp}$ & Damping & -- \\
    $\mathrm{e}$ & Exact solution & -- \\
    $i,\,j$ & Node indices & -- \\
    $ij$ & Edge or staggered point between nodes $i$ and $j$ & -- \\
    $i\rightarrow$ & Radial projection originating from node $i$ & -- \\
    $\mathrm{rdm}$ & Random & -- \\
    \hline
    \multicolumn{3}{@{}l}{\textit{Superscripts}} \\
    $*$ & Intermediate or reference value & -- \\
    $h$ & Polynomial approximation & -- \\
    $n$ & Time step & -- \\
    \hline
  \end{tabular}
\end{center}

\section{Introduction}
Meshfree/meshless methods are a class of numerical methods that represent the computational domain using discrete points (DPs), and have been developed to simulate practical problems without the time-consuming and labor-intensive process of mesh generation.
Broadly, there are two types of meshfree methods. One is based on weak formulations, such as the element-free Galerkin method \cite{belytschko_ElementfreeGalerkinMethods_1994,wang_SuperiorityMixedElement_2012}, the meshless local Petrov--Galerkin method \cite{atluri_NewMeshlessLocal_1998,abbaszadeh_DirectMeshlessLocal_2020}, and reproducing kernel particle method \cite{liu_ReproducingKernelParticle_1995}. The other type is based on strong formulations, such as the smoothed particle hydrodynamics (SPH) \cite{gingold_SmoothedParticleHydrodynamics_1977}, the moving particle semi-implicit/simulation method \cite{koshizuka_MovingParticleSemiImplicitMethod_1996}, the finite-point method \cite{onate_StabilizedFinitePoint_1996}, and the generalized finite-difference method \cite{perrone_GeneralFiniteDifference_1975,liszka_FiniteDifferenceMethod_1980,benito_SolvingParabolicHyperbolic_2007,prieto_ApplicationGeneralizedFinite_2011,tang_LeastsquaresGeneralizedFinite_2025}.
The latter, strong-form meshfree methods, are considered truly meshfree because the governing equations are evaluated locally, thereby inherently avoiding use of global mesh structures.

Recently, incompressible flow solvers based on strong-form meshfree methods have been proposed \cite{dilts_MovingleastsquaresparticleHydrodynamicsConsistency_1999,suchde_MeshfreeGFDMSolvers_2018,tamai_LeastSquaresMoving_2014,matsunaga_ImprovedTreatmentWall_2020,vasyliv_SimulatingIncompressibleFlow_2020,matsuda_ParticlebasedMethodUsing_2022,matsuda_MeshconstrainedDiscretePoint_2025}. These use moving least-squares (MLS)\cite{lancaster_SurfacesGeneratedMoving_1981} for spatial discretization and Taylor series consistency.
Moreover, improved methods with more accurate higher-order temporal convergence have been developed \cite{matsunaga_StabilizedLSMPSMethod_2022,matsunaga_HighorderTimemarchingSchemes_2025}, and meshfree methods have shown increased usefulness in practice.
However, these strong-form meshfree methods for incompressible flows tend to suffer from unphysical checkerboard instability on flow fields, which arises from a zero-energy mode associated with the collocated arrangement of velocity and pressure variables. This makes it difficult to achieve the expected order of spatial convergence for both velocity and pressure.
The problem of zero-energy pressure modes has been known for a long time in mesh-based methods such as the finite difference, finite element, and finite volume methods.
This problem is due to the collocated arrangement of velocity and pressure variables.
In such cases, the pressure gradient and velocity are evaluated at inconsistent positions, which allows a nonphysical checkerboard pressure distribution \cite{hopman_QuantifyingCheckerboardProblem_2025}.
In the framework of meshfree methods, Swegle \textit{et al.} first reported the zero-energy mode problem in the context of SPH \cite{swegle_AnalysisSmoothedParticle_1994}.
The zero-energy mode is also reported to enhance tensile instability under negative pressure or strong vortices \cite{liu_DualParticleApproachIncompressible_2024}, and it is one of the difficult issues in meshfree methods.
Approaches employing polyharmonic spline radial basis functions (RBFs) augmented with appended polynomial terms have also been developed for the incompressible Navier--Stokes equations \cite{SHAHANE2021110623,Unnikrishnan2022,CHU2024112822,Unnikrishnan2024}. In these methods, a function is constructed by augmenting the radial basis function with polynomial terms of arbitrary order on a point cloud, which provides the advantage of achieving robust spatial convergence even when high-order functions are employed. However, the variables are arranged in a collocated manner, and it remains unclear how this affects numerical instabilities such as those described above, including odd-even decoupling.

Several improved meshfree schemes have been proposed, such as evaluating stresses at staggered positions relative to velocity points, and they have been effective in structural analysis \cite{dyka_ApproachTensionInstability_1995,randles_NormalizedSPHStress_2000,chalk_StressParticleSmoothedParticle_2020}. Similar staggered concepts have been introduced into meshfree methods for incompressible flows to suppress zero-energy modes and improve velocity--pressure coupling. He \textit{et al.} proposed staggered SPH formulations in which the velocity is represented at staggered particles or as a nonlocal two-point quantity, separately from the nodal pressure \cite{he_StaggeredMeshlessSolidfluid_2012,he_VariationalStaggeredParticle_2020}. Park \textit{et al.} introduced virtual staggered interpolation points and demonstrated the existence and stability of solutions for their meshfree Stokes formulation \cite{park_ExistenceStabilityVirtual_2016}. Liu \textit{et al.} proposed a dual-particle method in which virtual particles carrying pressure are regenerated at every time step \cite{liu_DualParticleApproachIncompressible_2024}. As a trade-off for the improved stability, their approach incurs additional computational cost owing to the generation and use of the virtual pressure particles. Chu and Schmidt \cite{CHU2023111756} proposed an RBF-based formulation employing staggered velocity and pressure nodes, which suppresses odd--even decoupling and achieves high computational accuracy. These formulations improve the stability of velocity--pressure coupling through different staggered or nonlocal constructions.

In addition, Trask \textit{et al.} \cite{trask_HighorderStaggeredMeshless_2017} have proposed a high-order staggered meshless method for elliptic boundary value problems that locally defines a virtual cell for each discrete point. This is motivated by the primal--dual relationship between cell-centered and cell-interface quantities in the finite volume method, and it enables consistent evaluation of discrete gradient and divergence operators in a truly meshfree manner. In particular, an MLS divergence reconstruction is applied to radial quantities evaluated along the edges, providing the basis for the radial-component reconstruction adopted. The formulation achieves high-order accuracy by increasing the polynomial order of the MLS reconstruction, with numerical results reported up to fourth-order accuracy. This method has also been extended to Stokes flow based on a monolithic formulation \cite{trask_CompatibleHighorderMeshless_2018}. In the Stokes formulation, incompressibility is incorporated through an MLS velocity reconstruction using polynomial basis functions that satisfy the divergence-free condition exactly.

The meshfree methods for incompressible flows discussed above can be broadly classified into monolithic and pressure-projection formulations according to how incompressibility is imposed. In a monolithic formulation, velocity and pressure are solved simultaneously as a coupled system, whereas a pressure-projection formulation \cite{chorin_NumericalSolutionNavierStokes_1968} separates the pressure solution from the velocity update. The latter avoids solving a single coupled system for velocity and pressure and therefore provides an attractive alternative to monolithic formulations. However, in conventional strong-form meshfree methods based on the pressure-projection formulation, the velocity divergence in the source term of the pressure Poisson equation, the pressure Laplacian, and the pressure gradient in the velocity-projection step may be evaluated using mutually inconsistent discrete operators. Consequently, solving the pressure Poisson equation does not necessarily ensure a divergence-free velocity field at the discrete level. To our knowledge, a strong-form meshfree pressure-projection formulation that ensures discrete consistency among these operators has not yet been established for the incompressible Navier--Stokes equations.

In this study, we propose a pressure-projection formulation within a least-squares meshfree method for the incompressible Navier--Stokes equations, using a staggered-variable arrangement that reduces the discrete velocity divergence to a level determined by the residual of the pressure Poisson equation.
The proposed method is based on the local primal--dual grid introduced by Trask \textit{et al.} \cite{trask_HighorderStaggeredMeshless_2017}; however, it employs this grid for the variable arrangement, in which the velocity and pressure are defined at nodes/DPs, while the radial components of the velocity are defined on the interfaces of the virtual dual cells.
The method allows a divergence-consistent MLS reconstruction with radial components for spatial discretization, which enables use of compatible discrete gradient and divergence operators.
Because of this property, the method achieves a local solenoidal velocity field on the order of the residual norm of the linear system for the pressure Poisson equation.
The pressure Poisson equation is solved in a staggered manner using the radial component of the velocity evaluated on virtual dual-cell interfaces.
In this framework, a time-evolution converting (TEC) formula \cite{xiao_SimpleCIPFinite_2004,xiao_NumericalSimulationsFreeinterface_2005,xiao_UnifiedFormulationCompressible_2006} is applied to update the radial component of the velocity at the virtual dual-cell interface with temporal information on the flow, and this time-evolved quantity is then used to obtain the pressure field.
Although the TEC formula differs from momentum-based interpolation (e.g., \cite{rhie_NumericalStudyTurbulent_1983,zhang_GeneralizedFormulationsRhie_2014}), it provides a strong coupling between velocity and pressure by temporally linking the two variables.
The proposed method extends the mesh-constrained discrete point (MCD) method, which is a type of strong-form meshfree method.
The MCD method has been developed as an efficient solver for incompressible flows with arbitrary boundaries \cite{matsuda_ParticlebasedMethodUsing_2022}, and it has already been extended to moving boundary problems \cite{matsuda_MeshconstrainedDiscretePoint_2025}.
In the MCD method, DPs are associated one-to-one with each mesh element of a structured mesh system (background mesh), which is laid over the analysis domain.
Each DP can be arbitrarily located inside its associated background mesh element, but it cannot exist outside of the cell boundary. This is so-called mesh constraint.
By arranging DPs along arbitrary boundary shapes under this mesh constraint, the method can realize both a faithful geometrical representation with DPs and efficient memory access using index information from the background mesh.
This discrete representation provides a uniform DP distribution at the level of the background mesh, and a constant stencil size can be expected.
Therefore, a second-order spatial discretization based on MLS can be stably evaluated with a compact $3 \times 3$ stencil.
The background mesh is used only for the mesh constraint and not for flow calculation or Gaussian integration, which makes the method different from existing mesh--particle hybrid methods (e.g., \cite{harlow_ParticleincellMethodNumerical_1962,brackbill_FLIPMethodAdaptively_1986,zhang_LeastsquaresMeshfreeMethod_2005,chandra_StabilizedMixedMaterial_2024,he_ImprovedMPMFormulation_2025,xu_CoupledSPHFVM_2021,nishiguchi_EulerianFiniteVolume_2024}).

The use of edge-based quantities, local connectivity, and virtual geometric entities in the proposed method is conceptually related to several conservative and finite-volume-like meshfree formulations developed to enforce or improve numerical flux conservation \cite{katz_meshless_2009,kwan-yu_chiu_conservative_2012,suchde_flux_2017,TRASK2020109187}. These formulations construct the discrete operators so as to enforce or improve a local balance of numerical fluxes. However, they are not specifically formulated to ensure a locally divergence-free velocity field through velocity--pressure coupling in a pressure-projection procedure. In contrast, the present formulation is designed to maintain discrete consistency among the operators used for velocity--pressure coupling, thereby directly relating the local discrete velocity divergence after projection to the residual of the pressure Poisson equation.

This paper contains the following sections.
In Section \ref{sec_flow_solver}, the fundamental formulation of the proposed method, comprising the spatial discretization scheme and velocity--pressure coupling method, is described.
Section \ref{sec_num_test} compares results obtained with the proposed methods with those of an existing meshfree method using the collocated variable arrangement, and investigates the effectiveness of the staggered-variable arrangement in the meshfree framework and applicability to incompressible Navier--Stokes flows.
In Section \ref{sec_conclusion}, the conclusion and remarks are given.

\section{A meshfree incompressible flow solver with consistent velocity--pressure coupling} \label{sec_flow_solver}
\subsection{Overview}
In this study, an incompressible flow solver is developed on the basis of the MCD method \cite{matsuda_ParticlebasedMethodUsing_2022,matsuda_MeshconstrainedDiscretePoint_2025}.
The spatial discretization scheme is based on MLS reconstruction with radial components \cite{trask_HighorderStaggeredMeshless_2017,trask_CompatibleHighorderMeshless_2018}, and spatial derivatives are evaluated using $3 \times 3$ compact stencils.
The MLS scheme satisfies Taylor series consistency, and settings such as the polynomial basis, scaling parameter, and influence radius are set as in \cite{matsuda_MeshconstrainedDiscretePoint_2025}.

Section \ref{subsec_spatial_discretization_stg} introduces the spatial discretization scheme constructed using the local primal--dual arrangement.
Section \ref{subsec_bd_treat_stg} describes the treatment of boundary conditions in the MLS reconstruction with radial components.
Section \ref{subsec_flow_solver_linear_acoustics} describes the basic concept of the proposed staggered formulation for coupling the velocity and pressure, where a linear acoustic equation is used as an essential example.
Section \ref{subsec_flow_solver_incomp_NS} presents the numerical solver for the unsteady incompressible Navier--Stokes equations based on the pressure-projection method.
Note that the MLS reconstruction based on the radial components of quantities introduced in \cite{trask_HighorderStaggeredMeshless_2017} are summarized in \ref{appendix_MLS_stg}.

\subsection{Spatial discretization with a local primal--dual arrangement} \label{subsec_spatial_discretization_stg}
\begin{figure} 
  \begin{center}
    \includegraphics[width=.8\linewidth]{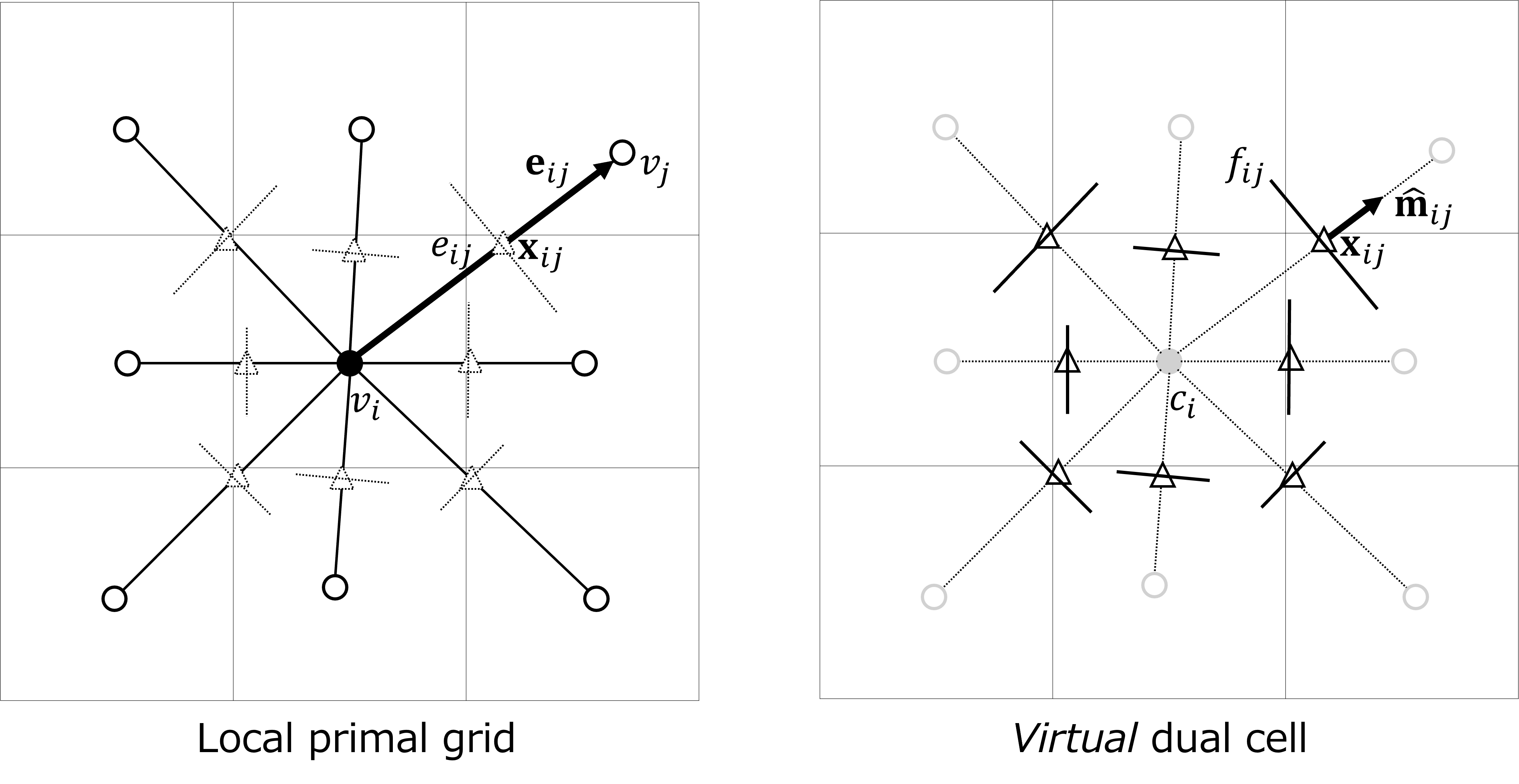}
    \caption{%
      Illustrations of the local primal grid (left) and \textit{virtual} dual cell (right).
    }
    \label{fig_arrng_local_primal_dual_grid}
  \end{center}
\end{figure}

In this study, the primal--dual arrangement \cite{trask_HighorderStaggeredMeshless_2017} is developed for MLS reconstruction with radial components (Fig. \ref{fig_arrng_local_primal_dual_grid}).
The domain $\Omega$ is discretized using a set of non-overlapping DP nodes (or vertices) $\mathcal{V} = \{ v_i \}^{N_{\rm DP}}_{i = 1}$, where each position is given as $\{ \mathbf{x}_{i} \}^{N_{\rm DP}}_{i = 1}$.
A set of local nodes $\mathcal{V}_{i}$ is introduced in the MCD framework as
\begin{equation}
  \mathcal{V}_{i}
  = \left\{
    v_{j} \in \mathcal{V} \,\middle|\, \mathbf{x}_{j} \in \mathcal{D}_{i}
  \right\},
  \label{eq_set_local_node}
\end{equation}
where $\mathcal{D}_{i}$ denotes a compact support domain for the node $v_{i}$ associated with the $3 \times 3$ background mesh, and a set of local edges $\mathcal{E}_{i}$ is introduced as
\begin{equation}
  \mathcal{E}_{i}
  = \left\{
    e_{ij} = \{ v_{i},\, v_{j} \} \,\middle|\, v_{j} \in \mathcal{V}_{i} ~\text{and}~  v_{i} \neq v_{j}
  \right\},
  \label{eq_set_local_edge}
\end{equation}
where $e_{ij}$ denotes the edge, which links $v_{i}$ to $v_{j}$.
The local primal--dual grid associated with $v_{i}$ is then given by $\mathcal{G}_{i} = (\mathcal{V}_{i},\, \mathcal{E}_{i})$ as a local graph consisting of $\mathcal{V}_{i}$ and $\mathcal{E}_{i}$ in Eqs. \eqref{eq_set_local_node} and \eqref{eq_set_local_edge}, respectively.
Regarding the local primal--dual grid, a virtual dual cell $c_{i}$ is introduced for each node $v_{i}$.
The boundary of $c_{i}$ consists of virtual faces $f_{ij}$, i.e., $\partial c_i = \bigcup_{e_{ij} \in \mathcal{E}_i} f_{ij}$, where each face intersects the edge $e_{ij}$ and has a unit normal vector $\hat{\mathbf{m}}_{ij}$ for the local primal--dual grid.
The intersection point $\mathbf{x}_{ij}$ between the edge $e_{ij}$ and face $f_{ij}$ can be arbitrarily defined as
\begin{equation}
  \mathbf{x}_{ij}
  = (1 - \theta_{ij}) \mathbf{x}_{i}
  + \theta_{ij} \mathbf{x}_{j},
  \label{eq_def_x_ij}
\end{equation}
where $\theta_{ij} \in (0,\, 1]$ denotes an arbitrary parameter for determining the intersection point $\mathbf{x}_{ij}$.
In this study, we set the parameter to

\begin{equation}
  \theta_{ij} = 
  \begin{cases}
    1, & \text{for}~ \mathbf{x}_{j} \in \Gamma_{d}, \\
    \frac{1}{2}, & \text{otherwise},
  \end{cases}
  \label{eq_def_theta_intersection}
\end{equation}

This indicates that the intersection point $\mathbf{x}_{ij}$ is set to coincide with the neighboring point $\mathbf{x}_{j}$ on the boundary $\Gamma_{d}$, where a normal component of a quantity is prescribed, while it is chosen as the midpoint of the edge $e_{ij}$ in the other cases. Accordingly, the unit normal vector of the virtual face $f_{ij}$ is set to $\hat{\mathbf{m}}_{ij}=\mathbf{n}_{j}$ for $\mathbf{x}_{j}\in \Gamma_{d}$, while it is set to $\hat{\mathbf{m}}_{ij}=\hat{\mathbf{e}}_{ij}$ in the other cases, where $\mathbf{n}_{j}$ is the unit normal vector of $\Gamma_{d}$ and $\hat{\mathbf{e}}_{ij} = \mathbf{e}_{ij}/\|\mathbf{e}_{ij}\|$ is the unit edge vector with $\mathbf{e}_{ij} = \mathbf{x}_{j} - \mathbf{x}_{i}$.
A set of intersection points $\mathcal{X}_{i}$ for the cell $c_{i}$ is defined as
\begin{equation}
  \mathcal{X}_{i}
  = \left\{
    \mathbf{x}_{ij} \,\middle|\, e_{ij} \in \mathcal{E}_{i}
  \right\}.
  \label{eq_Xi}
\end{equation}
It should be noted that, as in \cite{trask_HighorderStaggeredMeshless_2017}, the present study introduces a \textit{virtual} dual cell, which is based on the local connectivity of nodes, rather than a \textit{physical} one.
Therefore, the explicit geometric dimensions of the dual cell $c_{i}$ and its faces $f_{ij}$ are not required.
Although a common numerical flux per unit area is defined at each interface, the divergence is evaluated using the local MLS operator without cell volumes or face measures. Consequently, the discrete divergence is not expressed as a balance of surface-integrated fluxes, and the contributions from an internal interface do not necessarily cancel when the discrete equations associated with adjacent nodes are summed. Therefore, finite-volume-type local flux conservation is not guaranteed in the present formulation. Conservative meshfree formulations that introduce additional constructions to satisfy a discrete flux balance or divergence theorem have been proposed in the literature \cite{kwan-yu_chiu_conservative_2012, suchde_flux_2017, TRASK2020109187}.

By applying the primal--dual arrangement, we can express the radial component of the velocity $\mathbf{u} (\mathbf{x})$ at each intersection point $\mathbf{x}_{ij}$ as
\begin{equation}
  U_{i \rightarrow ij}
  = U_{i \rightarrow} (\mathbf{x}_{ij})
  = 2(\mathbf{x}_{ij} - \mathbf{x}_{i}) \cdot \mathbf{u}_{ij},
  \label{eq_def_rad_cmpt_Ui_ij}
\end{equation}
where $\mathbf{u}_{ij} = \mathbf{u} (\mathbf{x}_{ij})$.
When an approximate polynomial function $U^{h}_{i} (\mathbf{x})$ of $U_{i \rightarrow} (\mathbf{x}) = 2(\mathbf{x} - \mathbf{x}_{i}) \cdot \mathbf{u} (\mathbf{x})$ is introduced for a node $v_{i}$ (or cell $c_{i}$) at $\mathbf{x}_{i}$, $U^{h}_{i}$ can be approximated by the MLS reconstruction using $U_{i \rightarrow ij}$:
\begin{equation}
  \argmin_{{\bf c}_{i}}
  ~ \frac{1}{2} \sum_{\mathbf{x}_{ij} \in \mathcal{X}_{i}}
  w_{i}(\mathbf{x}_{ij})
  \left(
      U^{h}_{i} (\mathbf{x}_{ij}) - U_{i \rightarrow ij}
  \right)^2,
  \label{eq_def_obj_func_MLS_rad_cmpt}
\end{equation}
where the components of $\mathbf{c}_{i}$ are the coefficients of polynomial function $U^{h}_{i} (\mathbf{x})$.
Using a second-order polynomial basis, the approximation and its coefficient vector are given by $U^{h}_{i} (\mathbf{x}) = c_{i}^{10} (x - x_{i}) + c_{i}^{01} (y - y_{i}) + c_{i}^{20} (x - x_{i})^2 + c_{i}^{02} (y - y_{i})^2 + c_{i}^{11} (x - x_{i})(y - y_{i})$ and $\mathbf{c}_{i} = (c_{i}^{10}, c_{i}^{01}, c_{i}^{20}, c_{i}^{02}, c_{i}^{11})^{\mathsf{T}}$, respectively. Here, the constant term is omitted to enforce $U_i^h(\mathbf{x}_i)=0$, consistently with $U_{i\rightarrow}(\mathbf{x}_i)=0$.
Consequently, spatial derivatives at $\mathbf{x}_{i}$ are obtained as
\begin{equation}
  \mathbf{u}^{h} (\mathbf{x}_{i})
  = \frac{1}{2}\nabla U^{h}_{i} (\mathbf{x}_{i}),
  \label{eq_uh_at_xi}
\end{equation}
\begin{equation}
  \nabla \cdot \mathbf{u}^{h} (\mathbf{x}_{i})
  = \frac{1}{4} \nabla^{2} U^{h}_{i} (\mathbf{x}_{i}).
  \label{eq_div_uh_at_xi}
\end{equation}

These formulations can be extended to the tensor variables such as by approximating $\nabla\cdot{\bm \sigma}$, where ${\bm \sigma} \in \mathbb{R}^{d \times d}$ is a second-order tensor (e.g., Cauchy stress tensor). For a second-order tensor, polynomials are reconstructed for a vector quantity in the radial projection ${\bm \Sigma}_{i \rightarrow}(\mathbf{x})=2(\mathbf{x}-\mathbf{x}_{i})\cdot{\bm \sigma}(\mathbf{x}) \in \mathbb{R}^{d}$, and $\nabla \cdot \boldsymbol{\sigma}$ is approximated as
\begin{equation}
  \nabla \cdot \boldsymbol{\sigma} (\mathbf{x}_{i})
  \approx \nabla \cdot \boldsymbol{\sigma}^{h} (\mathbf{x}_{i})
  = \frac{1}{4} \nabla^{2} \boldsymbol{\Sigma}^{h}_{i} (\mathbf{x}_{i}).
\end{equation}

As described above, in the MLS reconstruction based on the primal--dual grid used in this study, the local neighborhood is defined as a set of nodes and edges. Therefore, the present formulation is expected to remain applicable to neighborhoods larger than those considered here, provided that the resulting MLS problem has full rank. In \cite{trask_CompatibleHighorderMeshless_2018}, fourth-order reconstructions of both velocity and pressure were applied to Stokes flow by using neighborhoods containing a sufficient number of points for the MLS reconstruction. The effectiveness of such high-order MLS reconstructions within the proposed formulation for incompressible flows on a staggered primal--dual grid remains to be investigated in future work.

\subsection{Treatment of boundary conditions on a primal--dual arrangement} \label{subsec_bd_treat_stg}
\begin{figure} 
  \begin{center}
    \includegraphics[width=.45\linewidth]{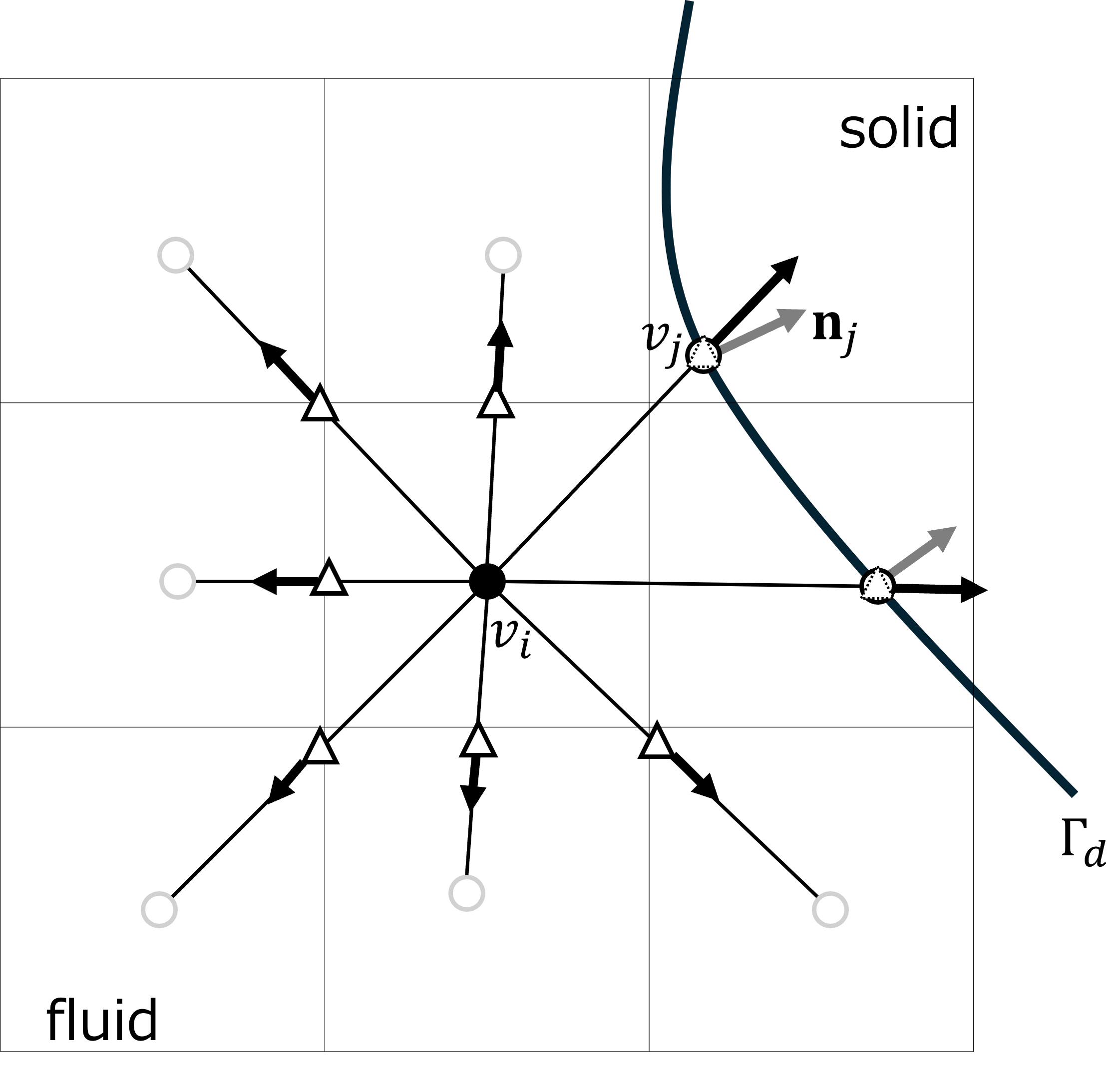}
    \caption{%
      Illustration of applying the normal component of any vector (or tensor) quantity as a boundary condition on $\Gamma_{d}$.
    }
    \label{fig_Imposing_Dirichlet_bd}
  \end{center}
\end{figure}

We extend the formulation using the primal-dual arrangement in the case that a normal component of any quantity is imposed as a boundary condition on the boundary $\Gamma_{d}$.
The set of intersection points $\mathcal{X}_i$ is decomposed into the boundary-related subset and its complement as
\begin{equation}
  \mathcal{X}_{i}^{d}
  =
  \left\{
    \mathbf{x}_{ij} \in \mathcal{X}_{i}
    \,\middle|\,
    \mathbf{x}_{j} \in \Gamma_{d}
  \right\},
  \quad \textrm{and} \quad
  \tilde{\mathcal{X}}^{d}_{i}
  = \mathcal{X}_{i}\setminus\mathcal{X}_{i}^{d}.
  \label{eq_Xi_Gammad}
\end{equation}
Using Eqs. \eqref{eq_dU} and \eqref{eq_ddU}, the gradient of $U_{i \rightarrow} (\mathbf{x})$ is expressed as
\begin{equation}
  \begin{split}
  \nabla U_{i \rightarrow}
  &= 2\mathbf{u}
    + 2\nabla\mathbf{u} \cdot \mathbf{e}_{i}, \\
  & = 2\mathbf{u}
    + ((\nabla\mathbf{u}+\nabla\mathbf{u}^{\top}) + (\nabla\mathbf{u}-\nabla\mathbf{u}^{\top})) \cdot \mathbf{e}_{i}, \\
  & = 2\mathbf{u}
    + \frac{1}{2}\nabla\nabla U_{i \rightarrow} \cdot \mathbf{e}_{i}
    - \mathbf{e}_{i} \cdot \nabla\nabla\mathbf{u}\cdot\mathbf{e}_{i}
    + (\nabla\mathbf{u}-\nabla\mathbf{u}^{\top}) \cdot \mathbf{e}_{i},
  \end{split}
  \label{eq_dU_expand}
\end{equation}
where $\mathbf{e}_{i}(\mathbf{x})=\mathbf{x}-\mathbf{x}_{i}$ is the radial vector.
Assuming $\mathbf{x}_{ij} \in \mathcal{X}^{d}_{i}$ (Fig. \ref{fig_Imposing_Dirichlet_bd}), the normal projection using the normal vector $\mathbf{n}_{j}$ on $\Gamma_{d}$ for Eq. \eqref{eq_dU_expand} is given by
\begin{equation}
  \mathbf{n}_{j} \cdot \nabla U_{i \rightarrow}
  = 2\mathbf{n}_{j}\cdot\mathbf{u}
    + \frac{1}{2}(\mathbf{n}_{j}\mathbf{e}_{i}:\nabla\nabla) U_{i \rightarrow}
    - (\mathbf{n}_{j}\mathbf{e}_{i} : \nabla\nabla)\mathbf{u}\cdot\mathbf{e}_{i}
    + \mathbf{n}_{j} \cdot(\nabla\mathbf{u}-\nabla\mathbf{u}^{\top}) \cdot \mathbf{e}_{i}.
\label{eq_n_dot_dU}
\end{equation}
At $\mathbf{x}_{ij}$, the third term on the right-hand side of Eq. \eqref{eq_n_dot_dU} is then of order $O(\|\mathbf{e}_{ij}\|^2)$, provided that $\nabla\nabla\mathbf{u}$ is bounded.
The last term vanishes when the edge direction $\mathbf{e}_{ij}$ is aligned with the boundary normal $\mathbf{n}_{j}$, owing to the skew-symmetry of $\nabla\mathbf{u}-\nabla\mathbf{u}^{\top}$. More generally, this term is induced by the tangential component of $\mathbf{e}_{ij}$ and can be interpreted as a non-orthogonality correction. Although this correction may scale with $\|\mathbf{e}_{ij}\|$ in the same way as the second term, we omit it in the present approximation to keep the boundary treatment local and simple.
Therefore, Eq. \eqref{eq_n_dot_dU} is approximated as
\begin{equation}
  \mathbf{n}_{j} \cdot \nabla U_{i \rightarrow}(\mathbf{x}_{ij})
  \approx 
  2(\mathbf{n}\cdot\mathbf{u})_{j}
  + \frac{1}{2} (\mathbf{n}_{j}\mathbf{e}_{ij}:\nabla\nabla) U_{i \rightarrow}(\mathbf{x}_{ij}).
\label{eq_n_dot_dU_approx}
\end{equation}
Using Eq. \eqref{eq_n_dot_dU_approx} to impose a boundary condition $q_{j} = (\mathbf{n} \cdot \mathbf{u})_{j}$ on $\Gamma_{d}$, we then construct the approximate polynomial $U^{h}_{i}(\mathbf{x})$ extended from Eq. \eqref{eq_def_obj_func_MLS_rad_cmpt} as
\begin{equation}
\begin{split}
  \argmin_{{\bf c}_{i}}
  ~ &\frac{1}{2} \sum_{\substack{\mathbf{x}_{ij}\in\tilde{\mathcal{X}}^{d}_{i}}}
  w_{i}(\mathbf{x}_{ij})
  \left(
    U^{h}_{i} (\mathbf{x}_{ij}) - U_{i \rightarrow ij}
  \right)^2 \\
  &+ \frac{1}{2} \sum_{\mathbf{x}_{ij}\in\mathcal{X}^{d}_{i}}
  w_{i}(\mathbf{x}_{ij}) r^{2}_{s}
  \left(
    \mathbf{n}_{j}\cdot\nabla U^{h}_{i}(\mathbf{x}_{ij})
    - \frac{1}{2} (\mathbf{n}_{j}\mathbf{e}_{ij}:\nabla\nabla) U^{h}_{i}(\mathbf{x}_{ij})
    - 2q_{j}
  \right)^2,
\end{split}
\label{eq_def_obj_func_with_bd_treatment}
\end{equation}
where $r_{s}$ denotes the scaling parameter.
The second term of Eq. \eqref{eq_def_obj_func_with_bd_treatment} reflects the boundary condition, which is inspired by the formulation that applies the Neumann condition for $U^{h}_{i}$ proposed in \cite{matsunaga_ImprovedTreatmentWall_2020,matsuda_ParticlebasedMethodUsing_2022,matsuda_MeshconstrainedDiscretePoint_2025}.

\subsection{Formulation of the linear acoustic equation} \label{subsec_flow_solver_linear_acoustics}
We formulate the coupling of the velocity and pressure in the staggered-variable arrangement for a linear acoustic equation.
The governing equations are expressed as
\begin{align}
  \frac{\partial p}{\partial t}
  &= - \rho c^{2} \nabla \cdot \mathbf{u}, \label{eq_dpdt_la} \\
  \rho \frac{\partial \mathbf{u}}{\partial t}
  &= - \nabla p, \label{eq_dudt_la}
\end{align}
where $\partial /\partial t$ and $\nabla = \partial /\partial \mathbf{x}$ denote the temporal and spatial derivative operators, respectively. $c$ denotes the acoustic speed, $\rho$ denotes the fluid density, and $\mathbf{u}$ and $p$ denote the velocity and pressure.

\begin{figure} 
  \begin{center}
    \includegraphics[width=.35\linewidth]{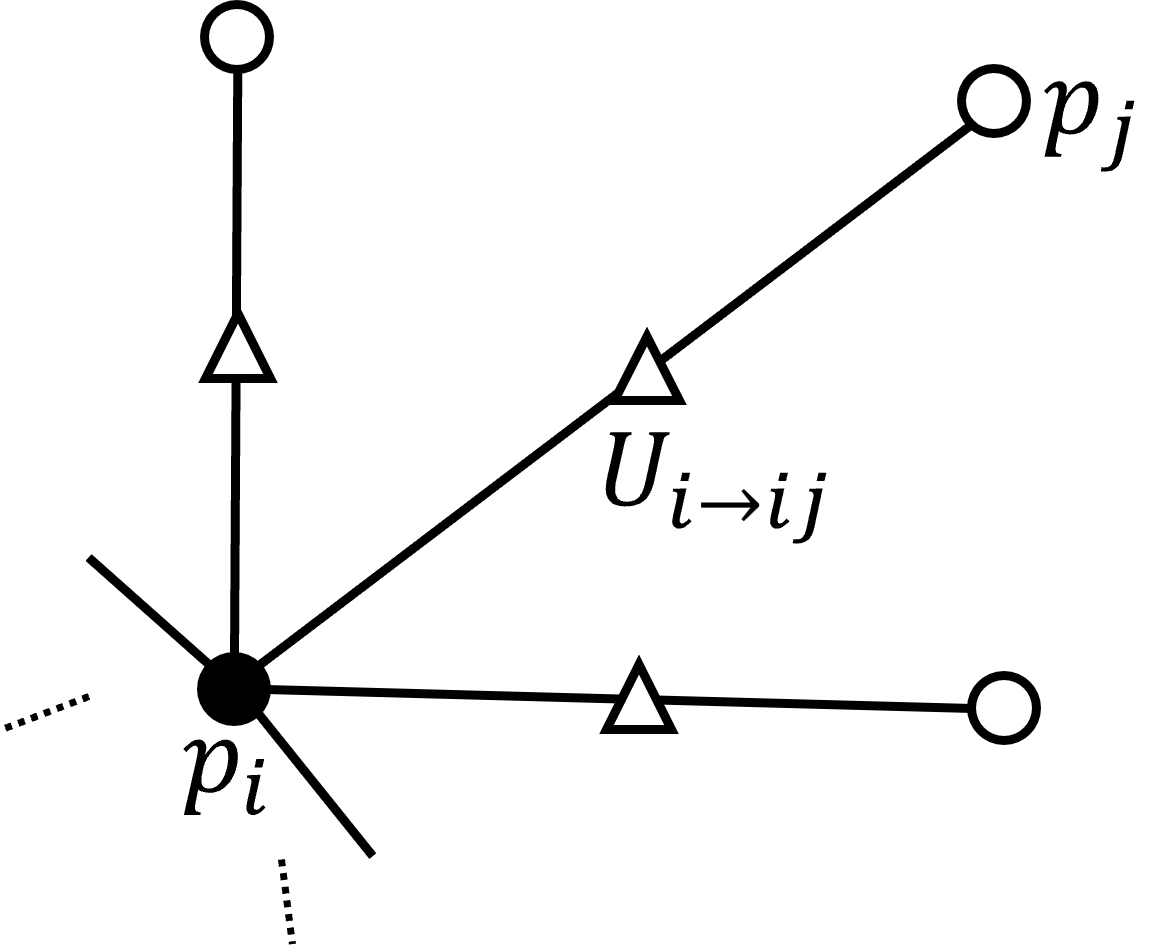}
    \caption{%
      Staggered-variable arrangement for the 2D linear acoustic equation on local primal--dual grids.
    }
    \label{fig_variable_linear_acoustics}
  \end{center}
\end{figure}

These equations are discretized using the staggered arrangement of velocity and pressure variables based on the local primal--dual grids (Fig. \ref{fig_variable_linear_acoustics}).
While the pressure is defined on each DP, the radial component $U_{i \rightarrow ij} = 2(\mathbf{x}_{ij} - \mathbf{x}_{i}) \cdot \mathbf{u} (\mathbf{x}_{ij})$ of the velocity is given on each midpoint $\mathbf{x}_{ij}$ of the edge vector $\mathbf{e}_{ij}$ (i.e., a staggered point).
According to the variable arrangement of velocity and pressure, the semi-discrete forms of Eqs. \eqref{eq_dpdt_la} and \eqref{eq_dudt_la} can be written as
\begin{equation}
  \frac{p^{n+1}_{i} - p^{n}_{i}}{\Delta t}
  = - \rho c^{2} \nabla \cdot \left. \mathbf{u}^{n} \right|_{i}
  = - \frac{\rho c^{2}}{4} \nabla^{2} \left. U^{h}_{i} \right|^{n}_{i},
  \label{eq_semi_dscrt_dpdt_la}
\end{equation}
\begin{equation}
  \rho \frac{U^{n+1}_{i \rightarrow ij} - U^{n}_{i \rightarrow ij}}{\Delta t}
  = - 2(\mathbf{x}_{ij} - \mathbf{x}_{i}) \cdot \nabla \left. p^{n+1} \right|_{ij}.
  \label{eq_semi_dscrt_dudt_la}
\end{equation}

In this study, we apply the second-order approximation for the MLS reconstruction of $U^{h}_{i}({\bf x})$ and the radial projection of the pressure gradient on the right side of Eq. \eqref{eq_semi_dscrt_dudt_la}.
The Taylor series expansions of $p_i = p (\mathbf{x}_{i})$ and $p_j = p (\mathbf{x}_{j})$ corresponding to the intersection point $\mathbf{x}_{ij}$ can be written as
\begin{align}
  p_{i}
  &= p_{ij}
  -\theta_{ij}\mathbf{e}_{ij} \cdot \nabla \left. p \right|_{ij}
  + \frac{1}{2} \theta^{2}_{ij}\mathbf{e}_{ij}\mathbf{e}_{ij} : \left. \nabla\nabla p \right|_{ij}
  + \mathcal{O} (\| \mathbf{e}_{ij} \|^{3}),
  \label{eq_pi_taylor_apprx} \\
  p_{j}
  &= p_{ij}
  + (1-\theta_{ij})\mathbf{e}_{ij} \cdot \nabla \left. p \right|_{ij}
  + \frac{1}{2} (1-\theta_{ij})^{2}\mathbf{e}_{ij}\mathbf{e}_{ij} : \left. \nabla\nabla p \right|_{ij}
  + \mathcal{O} (\| \mathbf{e}_{ij} \|^{3}).
  \label{eq_pj_taylor_apprx}
\end{align}
Here, using Eq. \eqref{eq_def_x_ij}, we have $\mathbf{x}_{i}-\mathbf{x}_{ij} = -\theta_{ij}\mathbf{e}_{ij}$ and $\mathbf{x}_{j}-\mathbf{x}_{ij} = (1-\theta_{ij})\mathbf{e}_{ij}$, where $\mathbf{e}_{ij}=\mathbf{x}_{j}-\mathbf{x}_{i}$.
By subtracting Eq. \eqref{eq_pi_taylor_apprx} from Eq. \eqref{eq_pj_taylor_apprx}, we obtain
\begin{equation}
  \begin{split}
    p_{j} - p_{i} 
    &= \mathbf{e}_{ij} \cdot \left. \nabla p \right|_{ij}
    + \frac{1}{2} (1-2\theta_{ij}) \mathbf{e}_{ij}\mathbf{e}_{ij} : \left. \nabla\nabla p \right|_{ij}
    + \mathcal{O}(\| \mathbf{e}_{ij} \|^{3}).
  \end{split}
  \label{eq_grad_p_general}
\end{equation}
When $\mathbf{x}_{ij}$ is chosen as the midpoint of the edge $e_{ij}$, namely $\theta_{ij}=1/2$, the coefficient $1-2\theta_{ij}$ in Eq. \eqref{eq_grad_p_general} vanishes. Hence, in this case, $p_j-p_i$ approximates the radial projection $\mathbf{e}_{ij}\cdot\left.\nabla p\right|_{ij}$ with a truncation error of $\mathcal{O}(\|\mathbf{e}_{ij}\|^3)$ as in \cite{trask_HighorderStaggeredMeshless_2017}.
Thus, $2(\mathbf{x}_{ij} - \mathbf{x}_{i}) \cdot \nabla \left. p \right|_{ij} = \mathbf{e}_{ij}\cdot \nabla \left. p \right|_{ij} = p_{j}-p_{i}$ (for $\theta_{ij}=1/2$), and Eq. \eqref{eq_semi_dscrt_dudt_la} can be approximated by
\begin{equation}
  \rho \frac{U^{n+1}_{i \rightarrow ij} - U^{n}_{i \rightarrow ij}}{\Delta t}
  = - (p^{n+1}_{j} - p^{n+1}_{i}).
  \label{eq_semi_dscrt_dudt_la2}
\end{equation}

\subsection{Formulation of the incompressible Navier--Stokes equations} \label{subsec_flow_solver_incomp_NS}
The proposed formulation is applied to the incompressible Navier--Stokes equations given in the following conservative form:
\begin{equation}
  \nabla \cdot \mathbf{u} = 0, 
  \label{eq_continuity}
\end{equation}
\begin{equation}
  \frac{\partial (\rho \mathbf{u})}{\partial t} + \nabla \cdot (\rho \mathbf{u} \mathbf{u})
  = -\nabla p + \nabla \cdot \boldsymbol{\tau}, 
  \label{eq_NS}
\end{equation}
where $\rho$ denotes the fluid density, $\mathbf{u}$ and $p$ denote the velocity and pressure, and $\boldsymbol{\tau} = \eta \left( \nabla \mathbf{u} + (\nabla \mathbf{u})^{\top} \right)$ denotes the viscous stress tensor.
According to Eq. \eqref{eq_continuity}, the viscous term can be written as $\nabla \cdot \boldsymbol{\tau} = \eta \nabla \cdot \nabla \mathbf{u}$. Assuming a constant density ($\rho = const.$), Eq. \eqref{eq_NS} can be rewritten as
\begin{equation}
  \rho \left( \frac{\partial \mathbf{u}}{\partial t} + \nabla \cdot (\mathbf{u} \mathbf{u}) \right)
  = - \nabla p + \eta \nabla \cdot \nabla \mathbf{u}.
  \label{eq_NS2}
\end{equation}

Applying the pressure-projection method with the first-order explicit method, we can express the semi-discrete forms. The prediction step for the intermediate velocity $\mathbf{u}^{*}$ is given by
\begin{equation}
  \rho \left( \frac{\mathbf{u}^{*}_{i} - \mathbf{u}^{n}_{i}}{\Delta t}
  + \nabla \cdot \left.\left( \mathbf{u}^{n} \mathbf{u}^{n} \right)\right|_{i} \right)
  = \eta \nabla \cdot \nabla \left. \mathbf{u}^{n} \right|_{i},
  \label{eq_prediction_1}
\end{equation}
and then, the TEC formula is applied to update the radial component of the velocity at the virtual dual-cell interface as
\begin{equation}
  U^{*}_{i \rightarrow ij}
  = \mathrm{TEC} \left( U^{n}_{i \rightarrow ij},\, \mathbf{u}^{n}_{i},\, \mathbf{u}^{n}_{j},\, \mathbf{u}^{*}_{i},\, \mathbf{u}^{*}_{j} \right).
  \label{eq_prediction_stg_1}
\end{equation}
The pressure Poisson equation is then solved in a staggered manner using the radial component of the velocity evaluated on virtual dual-cell interfaces:
\begin{equation}
  \nabla \cdot \nabla \left. p^{n+1} \right|_{i}
  = \frac{\rho}{\Delta t} \nabla \cdot \left. \mathbf{u}^{*} \right|_{i}.
  \label{eq_pressure_1}
\end{equation}
Finally, the projection step for updating the velocity and its radial component is given by
\begin{equation}
  \rho \frac{U^{n+1}_{i \rightarrow ij} - U^{*}_{i \rightarrow ij}}{\Delta t}
  = - 2(\mathbf{x}_{ij} - \mathbf{x}_{i}) \cdot \nabla \left. p^{n+1} \right|_{ij},
  \label{eq_projection_stg_1}
\end{equation}
\begin{equation}
  \rho \frac{\mathbf{u}^{n+1}_{i} - \mathbf{u}^{*}_{i}}{\Delta t}
  = - \nabla \left. p^{n+1} \right|_{i}.
  \label{eq_projection_1}
\end{equation}

\begin{figure} 
  \begin{center}
    \includegraphics[width=.35\linewidth]{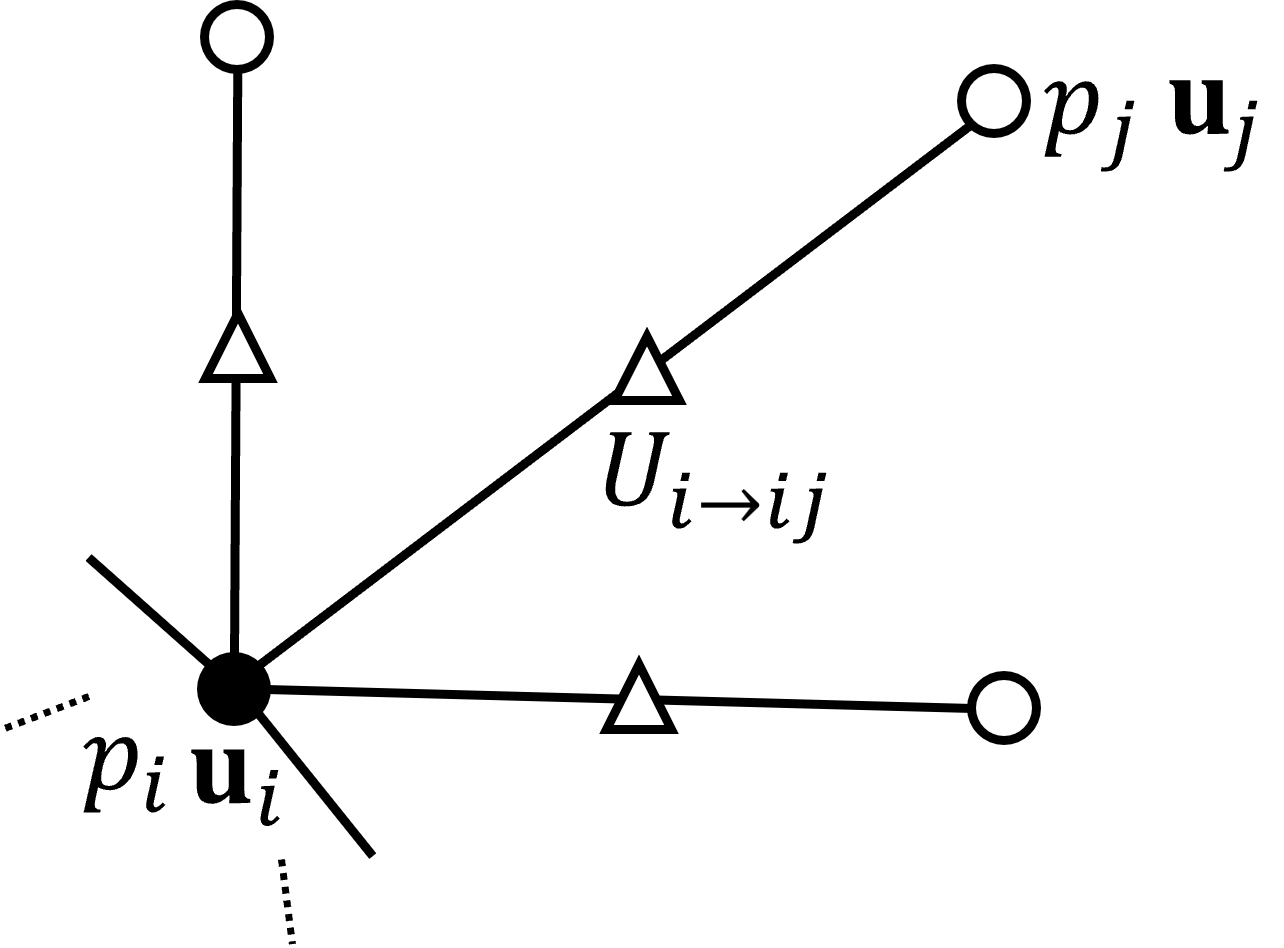}
    \caption{%
      Staggered-variable arrangement for the 2D Navier--Stokes equations on local primal--dual grids.
    }
    \label{fig_variable_NS}
  \end{center}
\end{figure}

For the spatial discretization in \eqref{eq_prediction_1}--\eqref{eq_projection_1}, the MLS reconstruction with radial components is applied based on the local primal--dual arrangement.
Fig. \ref{fig_variable_NS} illustrates staggered arrangements of $p_{i}$, $\mathbf{u}_{i}$, and $U_{i \rightarrow ij}$.
When the MLS reconstruction is applied with local primal--dual grids, the following radial components are defined:
\begin{align}
  U_{i \rightarrow ij}
  &= 2(\mathbf{x}_{ij} - \mathbf{x}_{i}) \cdot \left. \mathbf{u} \right|_{ij}, \label{eq_def_U_rad_cmpt} \\
  {\delta}P_{i \rightarrow ij}
  &= 2(\mathbf{x}_{ij} - \mathbf{x}_{i}) \cdot \nabla \left. p \right|_{ij}, \label{eq_def_nab_p_rad_cmpt} \\
  \mathbf{A}_{i \rightarrow ij}
  &= 2(\mathbf{x}_{ij} - \mathbf{x}_{i}) \cdot \left. (\mathbf{u} \mathbf{u}) \right|_{ij}, \label{eq_def_A_rad_cmpt} \\
  \mathbf{D}_{i \rightarrow ij}
  &= 2(\mathbf{x}_{ij} - \mathbf{x}_{i}) \cdot \nabla \left. \mathbf{u} \right|_{ij}. \label{eq_def_D_rad_cmpt}
\end{align}

When these definitions are used, the prediction step for the intermediate velocity $\mathbf{u}^{*}$ [Eq. \eqref{eq_prediction_1}] can be expressed as
\begin{equation}
  \rho \left( \frac{\mathbf{u}^{*}_{i} - \mathbf{u}^{n}_{i}}{\Delta t}
  + \frac{1}{4} \nabla^{2} \left. \mathbf{A}^{h}_{i} \right|^{n}_{i} \right)
  = \frac{1}{4} \eta \nabla^{2} \left. \mathbf{D}^{h}_{i} \right|^{n}_{i}.
  \label{eq_prediction_2}
\end{equation}
Here, $\mathbf{A}_{i \rightarrow ij}$ is approximated with an upwind technique:
\begin{equation}
  \mathbf{A}_{i \rightarrow ij}
  = \left\{
    \begin{array}{ll} \displaystyle
      U_{i \rightarrow ij} \left. \mathbf{u}^{h}_{i} \right|_{ij}, & \mathrm{for} \,\, U_{i \rightarrow ij} \geq 0, \\[6pt]
      U_{i \rightarrow ij} \left. \mathbf{u}^{h}_{j} \right|_{ij}, & \mathrm{otherwise},
    \end{array}
  \right.
\end{equation}
where $\left. \mathbf{u}^{h}_{i} \right|_{ij} = \mathbf{u}^{h}_{i} (\mathbf{x}_{ij})$ or $\left. \mathbf{u}^{h}_{j} \right|_{ij} = \mathbf{u}^{h}_{j} (\mathbf{x}_{ij})$ denotes the polynomial approximations of velocities at $\mathbf{x}_{ij}$ from $\mathbf{x}_{i}$ or $\mathbf{x}_{j}$.
These interpolations are carried out using the standard nodal-based MLS scheme (not local primal--dual grids) at nodes $v_{i}$ and $v_{j}$, as described in \cite{matsuda_ParticlebasedMethodUsing_2022,matsuda_MeshconstrainedDiscretePoint_2025}.
The pressure-projection parts [Eqs. \eqref{eq_pressure_1}--\eqref{eq_projection_1}] can be expressed as
\begin{equation}
  \nabla^{2} \left. {\delta}P^{h}_{i} \right|^{n+1}_{i}
  = \frac{\rho}{\Delta t} \nabla^{2} \left. U^{h}_{i} \right|^{*}_{i},
  \label{eq_pressure_2}
\end{equation}
\begin{equation}
  \rho \frac{U^{n+1}_{i \rightarrow ij} - U^{*}_{i \rightarrow ij}}{\Delta t}
  = - {\delta}P^{n+1}_{i \rightarrow ij},
  \label{eq_projection_stg_2}
\end{equation}
\begin{equation}
  \rho \frac{\mathbf{u}^{n+1}_{i} - \mathbf{u}^{*}_{i}}{\Delta t}
  = - \frac{1}{2}\nabla \left. {\delta}P^{h}_{i} \right|^{n+1}_{i}.
  \label{eq_projection_2}
\end{equation}
In a procedure analogous to that for the linear acoustic equation, the following approximation is employed:
\begin{align}
  {\delta}P_{i \rightarrow ij}
  &= p_{j} - p_{i}, \label{eq_deltaP_ij} \\
  \mathbf{D}_{i \rightarrow ij}
  &= \mathbf{u}_{j} - \mathbf{u}_{i}, \label{eq_D_ij}
\end{align}

The TEC formula \cite{xiao_SimpleCIPFinite_2004,xiao_NumericalSimulationsFreeinterface_2005,xiao_UnifiedFormulationCompressible_2006} is based on the following concept:
\begin{equation}
  \frac{\partial U_{i \rightarrow ij}}{\partial t}
  = \frac{1}{2} \left(
    \frac{\partial \left. \mathbf{U}^{h}_{i} \right|_{ij}}{\partial t}
    + \frac{\partial \left. \mathbf{U}^{h}_{j} \right|_{ij}}{\partial t}
  \right),
  \label{eq_idea_TEC}
\end{equation}
where $\left. \mathbf{U}^{h}_{i} \right|_{ij} = 2(\mathbf{x}_{ij}-\mathbf{x}_{i}) \cdot \left. \mathbf{u}^{h}_{i} \right|_{ij}$ and $\left. \mathbf{U}^{h}_{j} \right|_{ij} = 2(\mathbf{x}_{ij}-\mathbf{x}_{j}) \cdot \left. \mathbf{u}^{h}_{j} \right|_{ij}$ are the radial components of the interpolated velocities $\mathbf{u}^{h}_{i}(\mathbf{x}_{ij})$ and $\mathbf{u}^{h}_{j}(\mathbf{x}_{ij})$ from the nodes $v_{i}$ and $v_{j}$, respectively.
Applying the first-order method for discretizing Eq. \eqref{eq_idea_TEC}, we can write Eq. \eqref{eq_prediction_stg_1} as
\begin{equation}
  \begin{split}
  &\mathrm{TEC} \left(
    U^{n}_{i \rightarrow ij},\, \mathbf{u}^{n}_{i},\, \mathbf{u}^{n}_{j},\, \mathbf{u}^{*}_{i},\, \mathbf{u}^{*}_{j}
  \right) \\
  & \ = \beta \left(
    U^{n}_{i \rightarrow ij}
    + \frac{
        \left. \mathbf{U}^{h}_{i} \right|^{*}_{ij}
        - \left. \mathbf{U}^{h}_{i} \right|^{n}_{ij}
        + \left. \mathbf{U}^{h}_{j} \right|^{*}_{ij}
        - \left. \mathbf{U}^{h}_{j} \right|^{n}_{ij}
      }{2}
  \right) \\
  &\quad + (1 - \beta) \left(
    \frac{
      \left. \mathbf{U}^{h}_{i} \right|^{*}_{ij}
      + \left. \mathbf{U}^{h}_{j} \right|^{*}_{ij}
    }{2}
  \right).
  \label{eq_def_TEC}
  \end{split}
\end{equation}
Eq. \eqref{eq_def_TEC} considers the temporal variation in $U_{i \rightarrow ij}$, thereby linking the velocity and pressure in time. Here, the parameter $\beta \in [0,\, 1]$ is introduced to avoid the inconsistent coupling between $U_{i \rightarrow ij}$ and the nodal velocities $\mathbf{u}_{i}$ and $\mathbf{u}_{j}$. This parameter is set to $\beta=0.99$ in this study. From the viewpoint of suppressing the odd-even instability associated with the interpolation procedure, it is theoretically desirable to choose a value as close to $\beta=1$ as possible. However, when $\beta=1$, $U_{i \rightarrow ij}$ is updated solely through its temporal evolution without a contribution from the instantaneous nodal velocities, which may weaken the discrete coupling between $U_{i \rightarrow ij}$ and the nodal velocity field. This may lead to numerical instability, particularly at high Reynolds numbers. For this reason, we adopted $\beta=0.99$ so that an averaging operation on the velocity itself slightly contributes to the temporal evolution.

The boundary treatments for the radial quantities $U$, $\delta P$, $\mathbf{A}$, and $\mathbf{D}$ are described in detail in \ref{sec_bc_radial_components}.

\section{Numerical tests} \label{sec_num_test}
In this section, several numerical validations are presented. The simulations were performed using in-house codes written in Fortran 90 and executed with shared-memory parallelization using OpenMP. In each numerical validation, the number of iterations required to solve the pressure Poisson equation was on the order of several tens to several hundreds.
\subsection{Linear acoustics} \label{subsec_linear_acoustics}
Our first investigation was validating the influence of the staggered arrangements of the velocity and pressure variables in the meshfree framework.
The 2D linear acoustics problem with velocity and pressure variables was solved in the square domain $[0,\, L] \times [0,\, L]$.
Parameters were set to $L = 1$, $\rho = 1050$, and $c = 1000$.
The domain was discretized using $N \times N$ uniformly distributed DPs, and the calculations were performed at $N = 17$, $33$, $65$, $129$, and $257$.
The DP spacing $l_{0}$ (the background mesh width) can be written as $l_{0} = L / (N - 1)$.
The time interval was set to $\Delta t = 5 \times 10^{-8}$, so the Courant number corresponding to the acoustic speed was equal to $c \Delta t / l_{0} = 1.28 \times 10^{-2}$ at $N = 257$.
The calculation was executed until $t = 3 \times 10^{-4}$.
The initial velocity field was set to zero, and the initial pressure field was given as
\begin{equation}
  p (\mathbf{x},\, 0) = \left\{
    \begin{array}{ll} \displaystyle
      \frac{1}{2} \left( 1 + \cos{\left( \frac{\pi \| \mathbf{x} - \mathbf{x}_{\rm 0} \|}{R} \right)} \right), & if\, \| \mathbf{x} - \mathbf{x}_{\rm 0} \| < R, \\ [6mm]
      0, & otherwise.
    \end{array}
  \right.
  \label{eq_init_p_la}
\end{equation}
The parameters were set to $R = L/10$ and $\mathbf{x}_{\rm 0} = (L/2,\, L/2)^{\top}$.

\begin{figure} 
  \begin{center}
    \includegraphics[width=\linewidth]{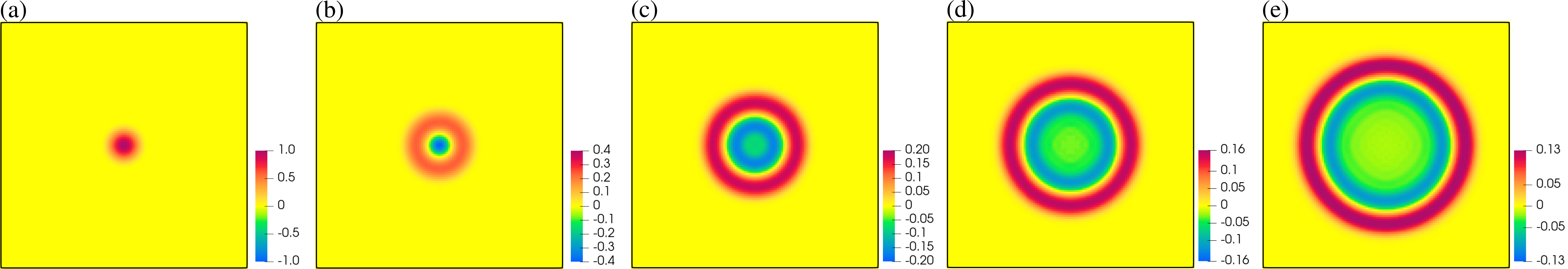}
    \caption{%
      Snapshots of the pressure field for $N = 129$ at five times: (a) $t = 0$, (b) $t = 7.5 \times 10^{-5}$, (c) $t = 1.5 \times 10^{-4}$, (d) $t = 2.25 \times 10^{-4}$, and (e) $t = 3 \times 10^{-4}$, calculated using the staggered-variable arrangement.
    }
    \label{fig_snapshot_p_linear_acoustics}
  \end{center}
\end{figure}

\begin{figure} 
  \begin{center}
    \includegraphics[width=\linewidth]{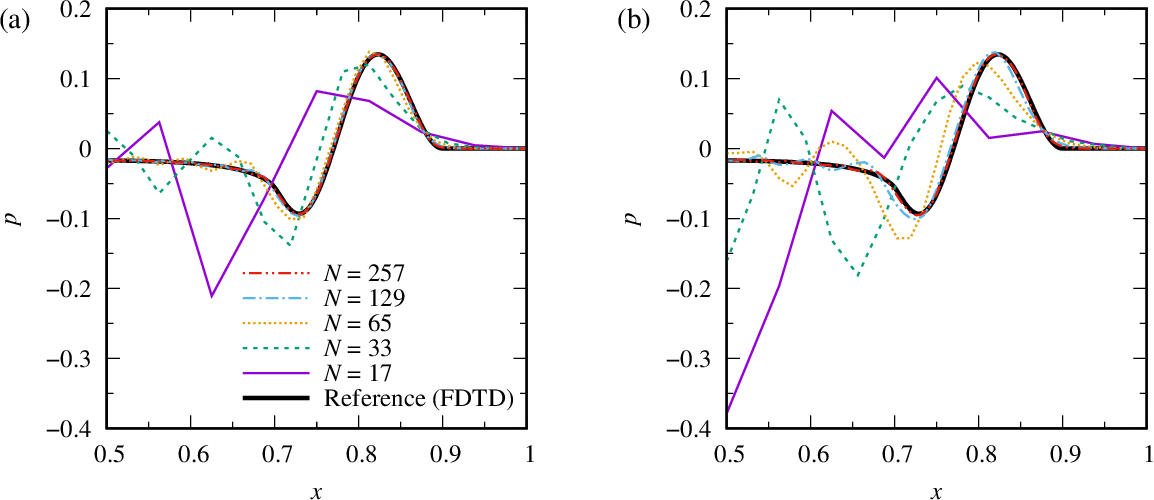}
    \caption{%
      Comparisons of the pressure $p$ along the $x$-axis for $x \in [L/2,\,L]$ using the MCD methods with the (a) staggered arrangement and (b) collocated arrangement.
      The black solid curves indicate the reference solution calculated with the FDTD method.
      The purple solid, green dashed, yellow dotted, blue dash-dot, and red dash-dot-dot curves correspond to resolutions $N = 17$, 33, 65, 129, and 257, respectively.
    }
    \label{fig_p_ref_linear_acoustics}
  \end{center}
\end{figure}

Fig. \ref{fig_snapshot_p_linear_acoustics} shows the snapshots of pressure distributions calculated using the proposed staggered method with $N = 129$ at five time instants.
These confirm that the pressure propagates axisymmetrically in a spatiotemporally smooth manner.
Fig. \ref{fig_p_ref_linear_acoustics} shows comparisons of the pressure distributions along the $x$-axis centerline in the range $[L/2,\, L]$ with the high-resolution reference solution obtained using the finite-difference time-domain (FDTD) method with $N = 2048$.
Fig. \ref{fig_p_ref_linear_acoustics} indicates that numerical solutions obtained using the proposed staggered method are closer to the high-resolution FDTD reference solution than those obtained via the conventional collocated method.

\begin{figure} 
  \begin{center}
    \includegraphics[width=\linewidth]{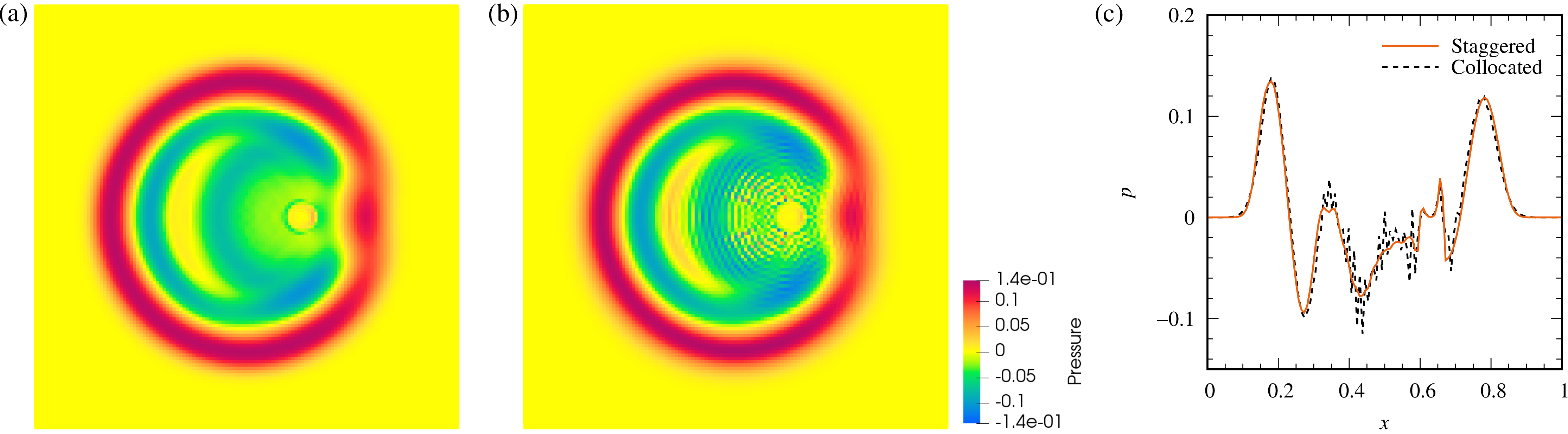}
    \caption{%
      Snapshots of pressure fields ($N = 129$) at $t = 3 \times 10^{-4}$ in the presence of an obstacle, calculated using (a) the staggered-variable arrangement and (b) collocated arrangement. (c) Comparison of the pressure distributions along the $x$-axis, where the orange solid and black dashed lines correspond to the staggered and collocated arrangements, respectively.
    }
    \label{fig_p_damping_linear_acoustics}
  \end{center}
\end{figure}

Fig. \ref{fig_p_damping_linear_acoustics} shows (a) a snapshot of the pressure field in a system with a cylindrical obstacle ($N=129$, $\Delta t = 10^{-7}$) at $t = 3 \times 10^{-4}$ obtained using the proposed staggered method, (b) one obtained via the conventional collocated method, and (c) the pressure distributions along the $x$-axis centerline.
Here, a cylindrical obstacle was modeled using the momentum equation [Eq. \eqref{eq_dudt_la}] with an additional damping term proportional to the velocity:
\begin{equation}
  \rho \frac{\partial \mathbf{u}}{\partial t} = -\nabla p - \alpha_{\rm damp} \rho \mathbf{u},
  \label{eq_la_damp}
\end{equation}
where $\alpha_{\rm damp}$ denotes the damping coefficient, which was defined using a constant $\alpha_{0}$ as
\begin{equation}
\begin{split}
  &\alpha_{\rm damp} (\mathbf{x}) \\
  &= \left\{
    \begin{array}{ll} \displaystyle
      \frac{\alpha_{0}}{2} \left( 1 + \cos{\left( \frac{\pi \| \mathbf{x} - \mathbf{x}_{\rm c, damp} \|}{R_{\rm damp}} \right)} \right), & if\, \| \mathbf{x} - \mathbf{x}_{\rm c, damp} \| < R_{\rm damp}, \\ [6mm]
      0, & otherwise.
    \end{array}
  \right.
\end{split}
  \label{eq_def_damp}
\end{equation}
In this study, the parameters were set to $\alpha_{0} = 10^{7}$, $R_{\rm damp} = L/20$, and $\mathbf{x}_{\rm c, damp} = (5L/8,\, L/2)^{\top}$.
According to Fig. \ref{fig_p_damping_linear_acoustics}, the proposed staggered method represents the pressure field without unphysical numerical oscillations.
In contrast, the conventional collocated method yields unphysical checkerboard patterns in the vicinity of the cylindrical obstacle (b), and the pressure distribution is no longer represented appropriately because the large numerical oscillations are sufficient to obscure the physical characteristics of the phenomenon (c).
These results confirm the effectiveness of the staggered arrangement of velocity and pressure variables in the meshfree framework.


\subsection{Flow around a stationary cylinder} \label{subsec_stat_cyl}
\begin{figure} 
  \begin{center}
    \includegraphics[width=.65\linewidth]{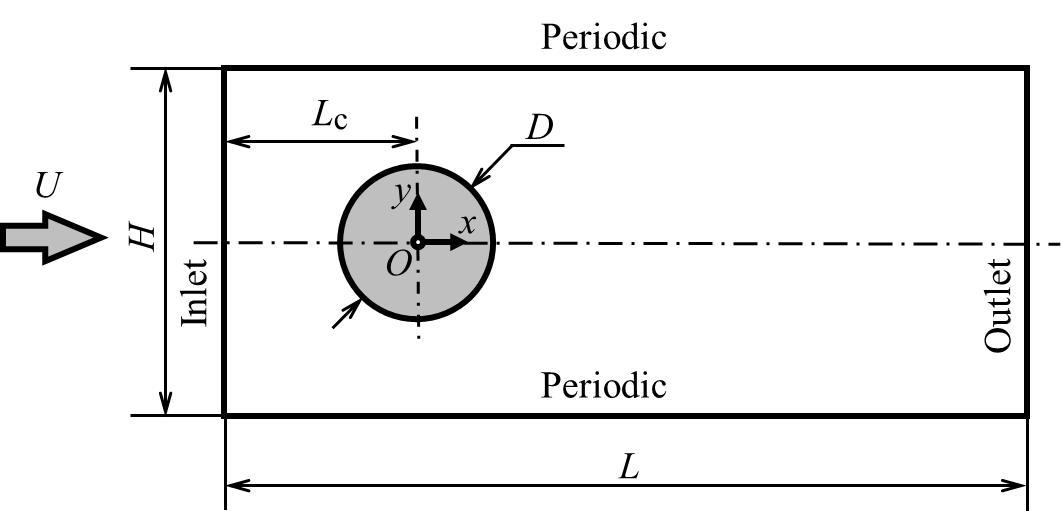}
    \caption{%
      Illustration of the fluid flow problem around a cylinder.
    }
    \label{fig_schematic_illust_Karman}
  \end{center}
\end{figure}

We calculated incompressible flows around a stationary cylinder of diameter $D$ (Fig. \ref{fig_schematic_illust_Karman}) and investigated the effectiveness of the proposed velocity--pressure coupling strategy for incompressible Navier--Stokes equations.
In this study, the cylinder diameter was set to $D = 1$, and the width and length of the channel were set to $H = 4D$ and $L = 8D$, respectively.
The center of the cylinder was at the origin of the coordinate system, and the offset distance of the cylinder center from the inlet was set to $L_{\rm c} = 2D$.
The no-slip condition was imposed on the cylinder wall ($\mathbf{u} = \mathbf{0}$).
Furthermore, the left-side channel boundary was subjected to the inlet condition ($\mathbf{u} = (U,\, 0)^{\top}$), and the right side was subjected to the outlet condition ($-p{\bf n} +\eta\partial \mathbf{u} / \partial n = \mathbf{0}$).
To eliminate the influence of boundaries other than the cylinder wall, the side walls of the channel were subjected to periodic boundary conditions.
The parameters were set to $U = 1$, $\rho = 1$, and $\eta = 5 \times 10^{-2}$, yielding Reynolds number $Re = \rho U D / \eta = 200$.
The background mesh width was set to $l_{0} = 1/32$, and the rectangular channel was represented using $257 \times 130$ DPs.
The reference time interval was set to $\Delta t^{*} = 5 \times 10^{-3}$, so the Courant number was equal to $U \Delta t^{*} / l_{0} = 0.16$.
In this study, the time intervals $\Delta t = \Delta t^{*}$ and $0.01 \Delta t^{*}$ were employed to investigate the influence of time resolution on the flow field.
The linear system of the pressure Poisson equation was solved using the Bi-CGSTAB method \cite{vandervorst_BiCGSTABFastSmoothly_1992}.
The iterations were stopped when the relative $L_{2}$ residual norm of the linear system fell below the convergence criterion $\varepsilon_{\rm crit}$, which was set to $\varepsilon_{\rm crit} = 10^{-5}$ in this study.

\begin{figure} 
  \begin{center}
    \includegraphics[width=\linewidth]{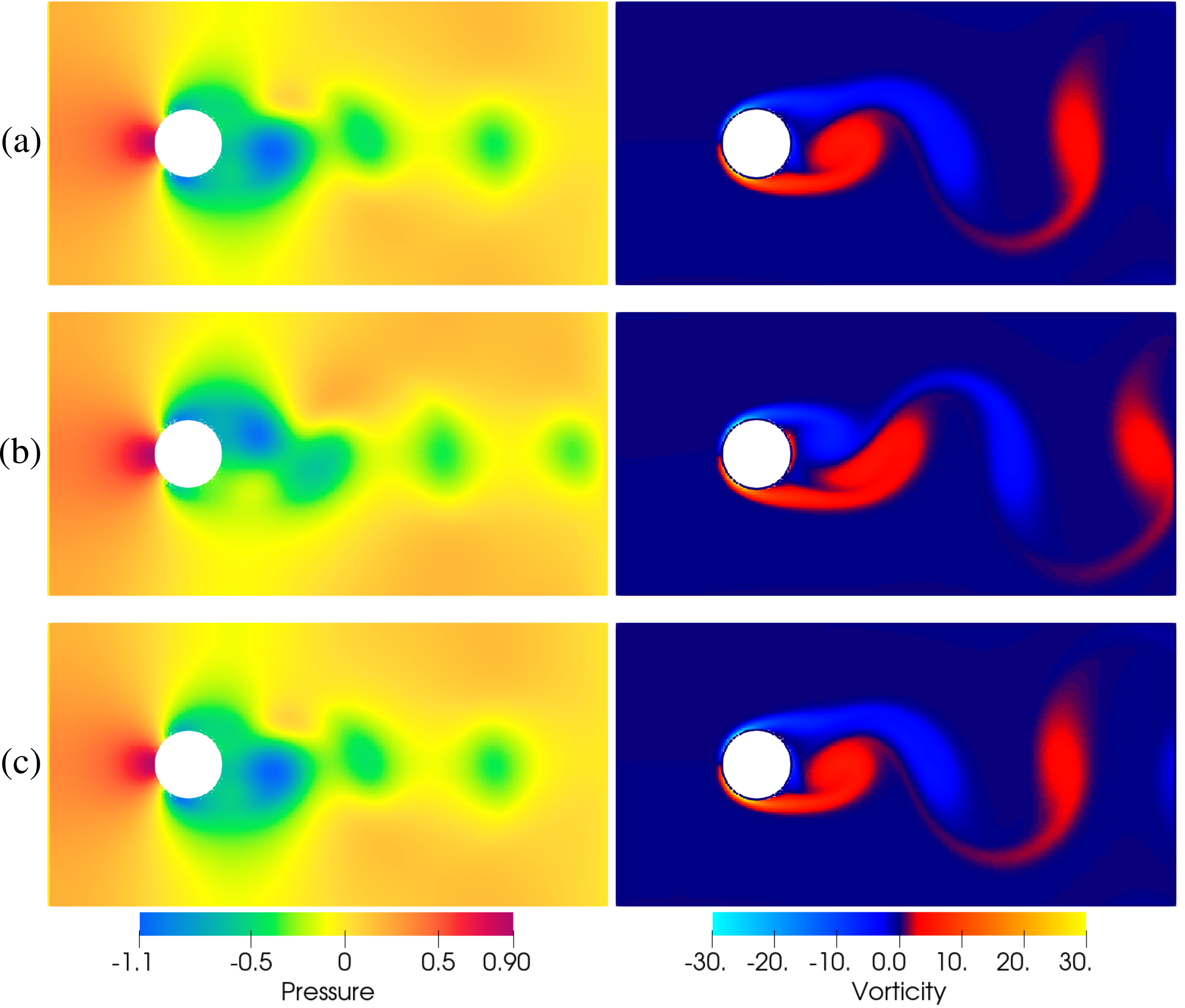}
    \caption{%
      Snapshots at $t = 100$ of pressure and vorticity fields around a cylinder at $Re = 200$ for $\Delta t = \Delta t^{*}$. (a) Proposed staggered arrangement with $\beta = 0.99$, (b) staggered arrangement with $\beta = 0$, and (c) conventional collocated arrangement.
    }
    \label{fig_p_vort_dt_Karman}
  \end{center}
\end{figure}

\begin{figure} 
  \begin{center}
    \includegraphics[width=\linewidth]{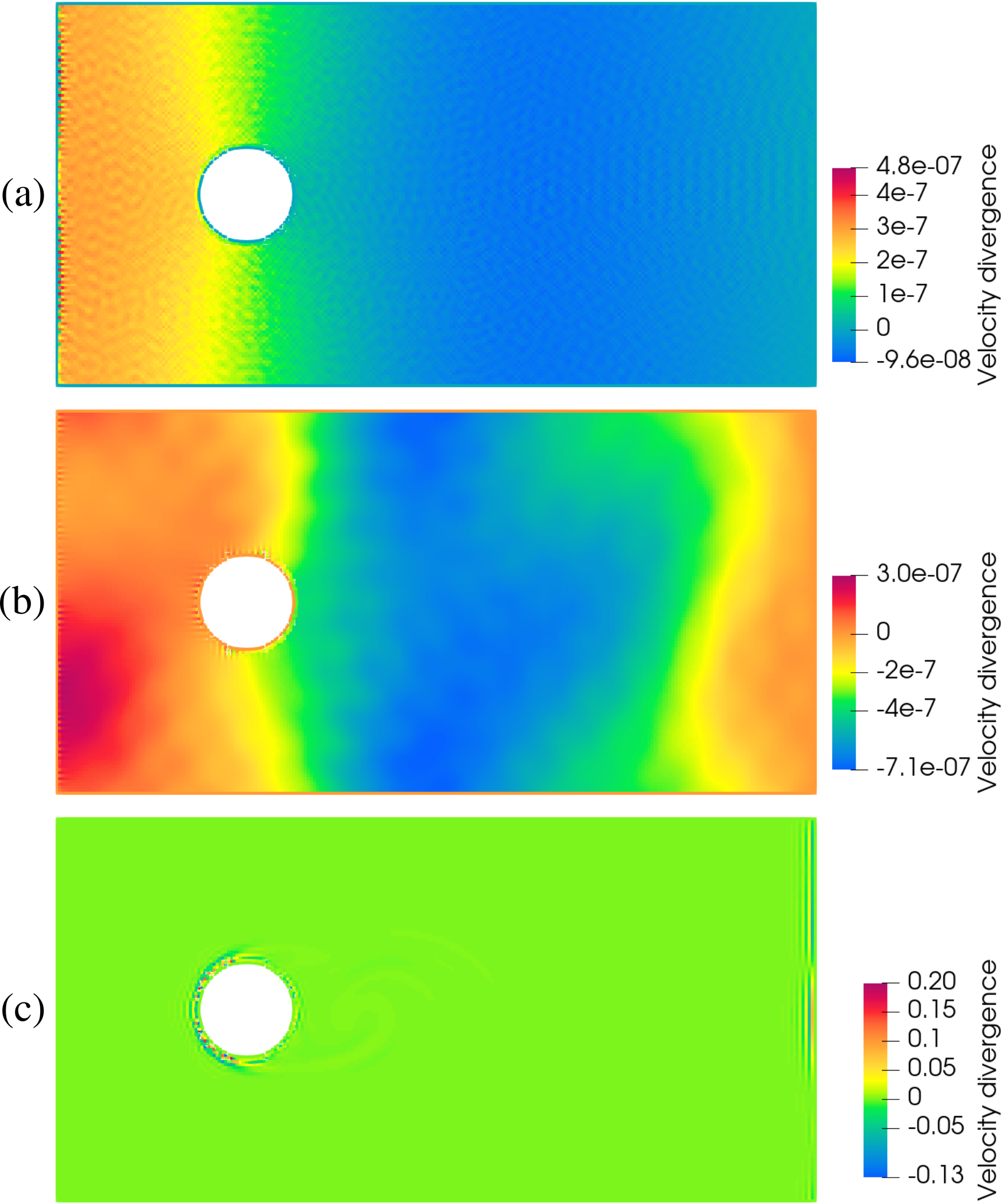}
    \caption{%
      Snapshots at $t = 100$ of velocity divergence around a cylinder at $Re = 200$ for $\Delta t = \Delta t^{*}$.
      (a) Proposed staggered arrangement with $\beta = 0.99$, (b) staggered arrangement with $\beta = 0$, and (c) conventional collocated arrangement.
    }
    \label{fig_vel_div_dt_Karman}
  \end{center}
\end{figure}

Fig. \ref{fig_p_vort_dt_Karman} shows the pressure and vorticity fields obtained using (a) the proposed staggered method with the TEC parameter $\beta = 0.99$, (b) the staggered method with $\beta = 0$, and (c) the conventional collocated method.
Fig. \ref{fig_vel_div_dt_Karman} shows velocity divergence fields.
According to Fig. \ref{fig_p_vort_dt_Karman}, the calculated flow fields are spatiotemporally smooth, and vortices are periodically shed in the cylinder wake at the coarse time resolution of $\Delta t = \Delta t^{*}$.
Fig. \ref{fig_vel_div_dt_Karman} indicates that, in the proposed staggered method (a, b), the velocity divergence has an order of magnitude consistent with that of the convergence criterion of the pressure Poisson equation, regardless of whether the TEC is enabled.
These results are attributed to the spatial discretization scheme based on the MLS reconstruction with radial components, which guarantees algebraic consistency between the discrete gradient and divergence operators (more details are shown in \ref{appendix_div_u_residual}).
In contrast, the result of the conventional collocated method [Fig. \ref{fig_vel_div_dt_Karman} (c)], which uses ordinary MLS for spatial discretization, yields a velocity divergence field with a much larger order of magnitude than that expected from the convergence criterion of the pressure Poisson equation.
These results show that the proposed staggered method achieves incompressible fluid flow at the discrete level in the meshfree framework.

\begin{figure} 
  \begin{center}
    \includegraphics[width=\linewidth]{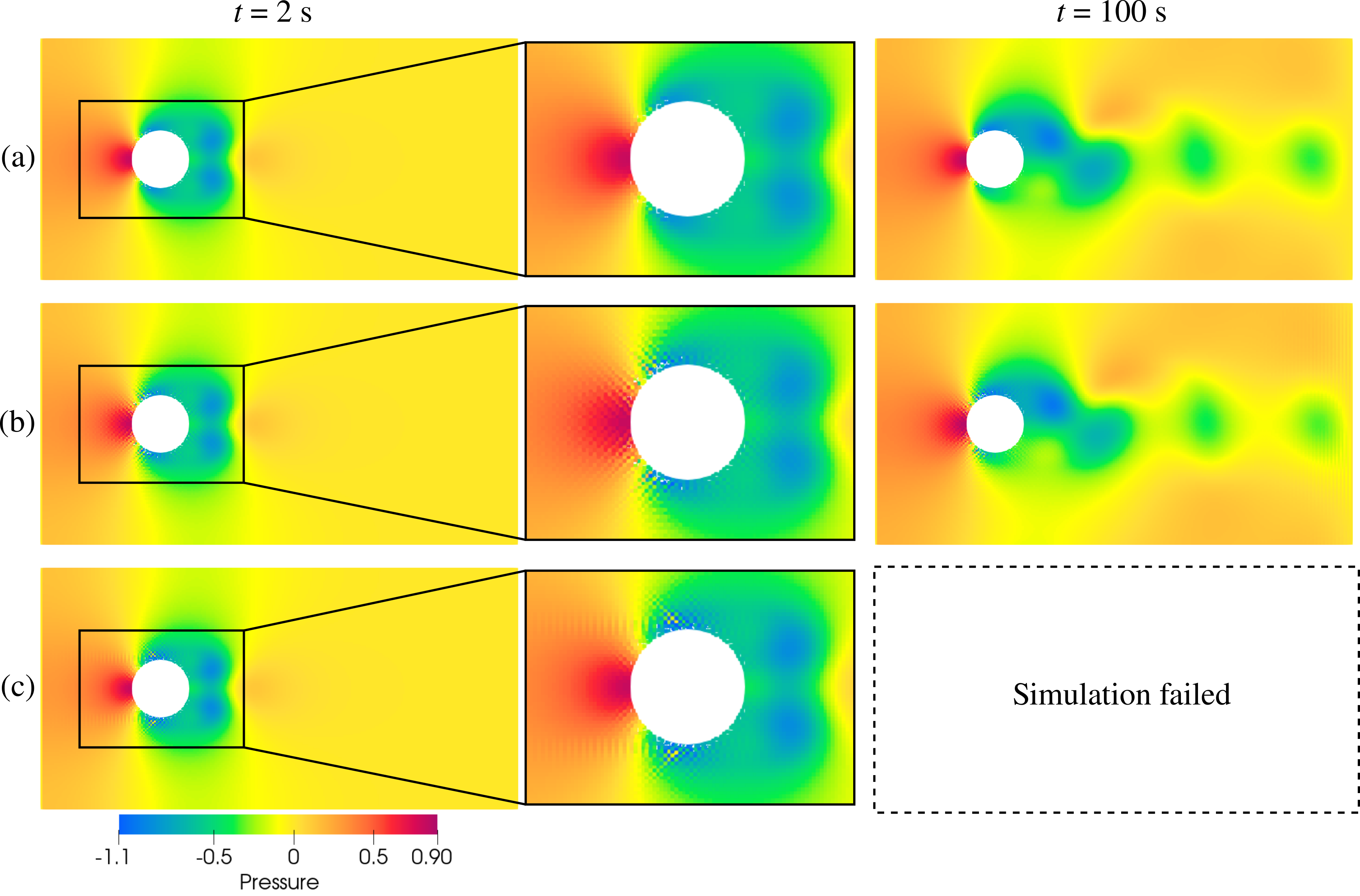}
    \caption{%
      Pressure distribution around a circular cylinder at $Re=200$ for $\Delta t = 0.01 \Delta t^{*}$.
      (a) Staggered arrangement with TEC ($\beta = 0.99$), (b) staggered arrangement without TEC ($\beta = 0$), and (c) conventional collocated arrangement.
    }
    \label{fig_p_0.01dt_Karman}
  \end{center}
\end{figure}

\begin{figure} 
  \begin{center}
    \includegraphics[width=\linewidth]{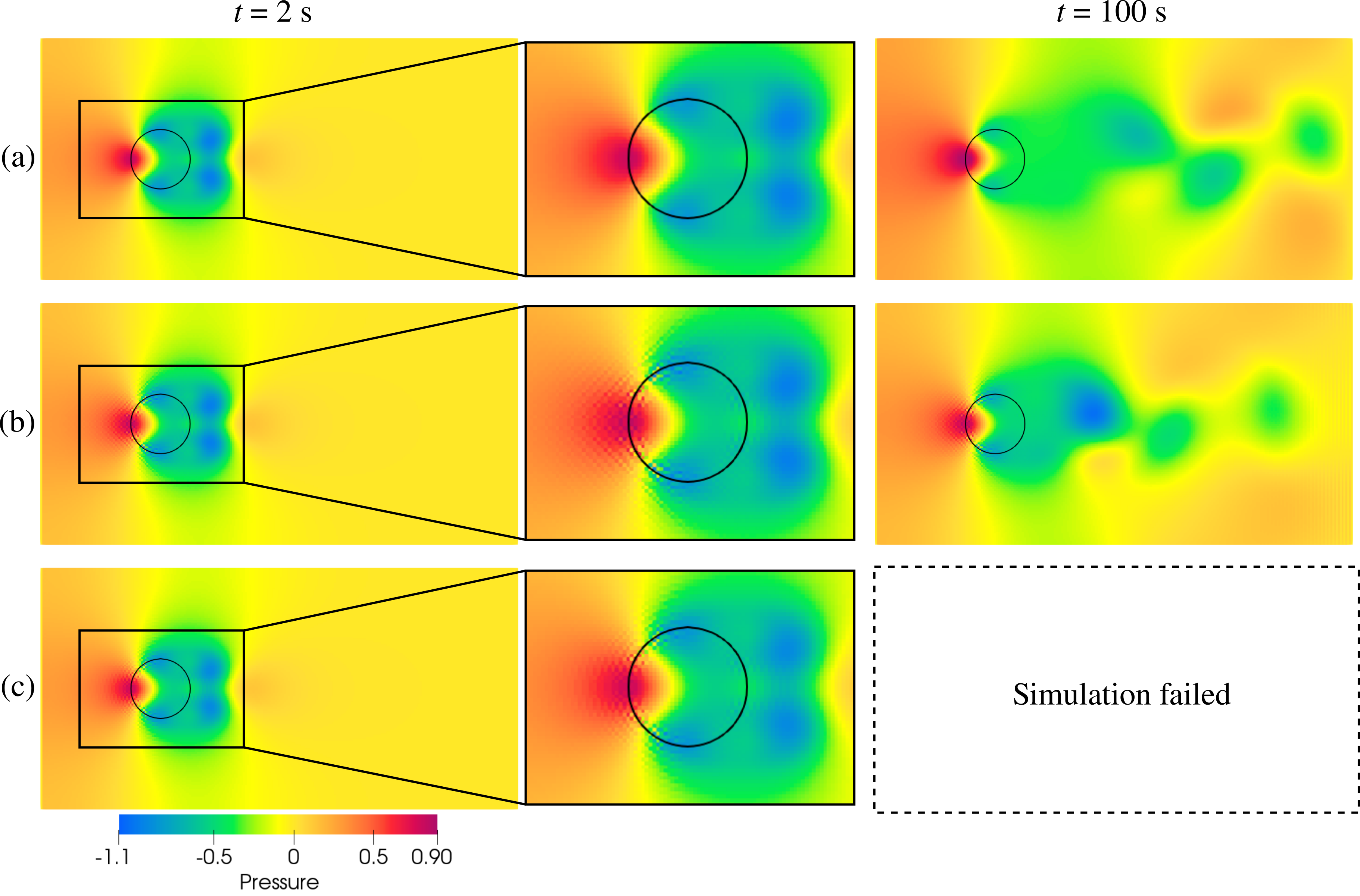}
    \caption{%
      Pressure distribution around a circular obstacle ($Re=200$, $\Delta t = 0.01 \Delta t^{*}$) represented by a damping term proportional to velocity.
      The black solid curve indicates the nominal position of the obstacle boundary.
      (a) Staggered arrangement with TEC ($\beta = 0.99$), (b) staggered arrangement without TEC ($\beta = 0$), and (c) conventional collocated arrangement.
    }
    \label{fig_p_0.01dt_damping_Karman}
  \end{center}
\end{figure}

Fig. \ref{fig_p_0.01dt_Karman} shows the snapshots of pressure fields around the cylinder at $t = 2$ and $100$ for $\Delta t = 0.01 \Delta t^{*}$ obtained using (a) the proposed staggered method with $\beta = 0.99$, (b) the proposed staggered method with $\beta = 0$, and (c) the conventional collocated method.
Fig. \ref{fig_p_0.01dt_damping_Karman} shows the pressure fields, where a cylindrical obstacle is represented by a uniform DP arrangement using the momentum equation [Eq. \eqref{eq_NS2}] with a damping term proportional to the velocity:
\begin{equation}
  \frac{\partial (\rho \mathbf{u})}{\partial t}
  + \nabla \cdot (\rho \mathbf{u} \mathbf{u})
  = -\nabla p + \nabla \cdot \boldsymbol{\tau} - \alpha_{\rm damp} \rho \mathbf{u},
  \label{eq_NS_damp}
\end{equation}
where $\alpha_{\rm damp}$ is also calculated from Eq. \eqref{eq_def_damp}.
The parameters were set to $\alpha_{0} = 10^{4}$, $R_{\rm damp} = D / 2$, and $\mathbf{x}_{\rm c, damp} = (0,\, 0)^{\top}$. This obstacle treatment enables us to investigate the influence of the wall boundary on the fluid flow while excluding the nonuniform DP arrangement along the cylinder.
According to Figs. \ref{fig_p_0.01dt_Karman} and \ref{fig_p_0.01dt_damping_Karman}, the proposed staggered method with $\beta = 0.99$ obtained a spatiotemporally smooth pressure field.
In contrast, the pressure showed checkerboard instability in the vicinity of the cylinder at $t = 2$ in the proposed staggered method with $\beta = 0$ (b) and the conventional collocated method (c), regardless of whether the DP arrangement was uniform.
Meanwhile, the flow calculation did not fail even if the conventional collocated method was used for an empty channel.
Thus, the presence of the obstacle itself causes the checkerboard instability in the pressure.
These results confirm that the velocity and pressure are temporally linked by the TEC formula and are consistently coupled.
Although the proposed method with $\beta = 0$ stably reproduced periodic vortex shedding despite checkerboard pressure oscillations near the cylinder, the conventional collocated method failed shortly after $t = 2$.
We attribute the successful calculation in the staggered method to the spatial discretization scheme, which satisfies the algebraic consistency of the discrete gradient and divergence operator, resulting in numerical incompressibility at the discrete level.

These results confirm that introducing a staggered arrangement of velocity and pressure variables together with a consistent velocity--pressure coupling with TEC in the meshfree framework enables stable calculation of the temporal evolution of incompressible flows.


\begin{figure} 
  \begin{center}
    \includegraphics[width=\linewidth]{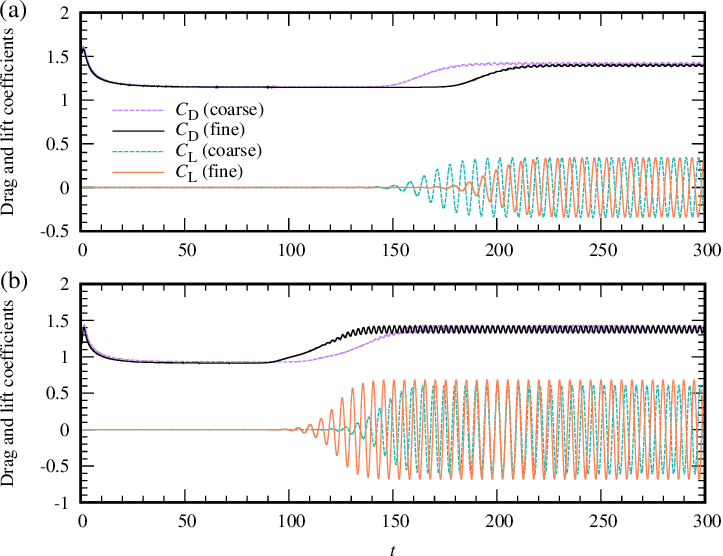}
    \caption{%
      Time change of $C_{\rm D}$ and $C_{\rm L}$ for the flow around a stationary cylinder at (a) $Re = 100$ and (b) $Re = 200$ calculated by the proposed staggered formulation with $\beta = 0.99$. Purple and blue-green dashed lines indicate the results with $l_{0} = 0.05$ and $\Delta t = 0.01$ (coarse), and black and orange solid lines indicate the results with $l_0=0.025$ and $\Delta t=0.005$ (fine).
    }
    \label{fig_cd_cl_Karman_20D}
  \end{center}
\end{figure}

In addition, the proposed staggered formulation with $\beta=0.99$ was applied to flow around a stationary cylinder in a wider computational domain to validate the drag and lift coefficients. The domain dimensions were $D=1$, $H=20D$, $L=40D$, and $L_{\rm c}=10D$. Simulations were performed at $Re=100$ and 200 using two spatial and temporal resolutions: $l_0=0.05$ and $\Delta t=0.01$ for the coarse resolution, and $l_0=0.025$ and $\Delta t=0.005$ for the fine resolution. Figure~\ref{fig_cd_cl_Karman_20D} shows the time histories of the drag coefficient $C_{\rm D}$ and lift coefficient $C_{\rm L}$ at (a) $Re=100$ and (b) $Re=200$. Table~\ref{tab:cd_cl_st_Karman} summarizes the mean drag coefficient, the fluctuation amplitudes of the drag and lift coefficients, and the Strouhal number $St$. For comparison, Liu \textit{et al.} \cite{liu_PreconditionedMultigridMethods_1998} used a finite-difference method with a body-fitted mesh, whereas Ding \textit{et al.} \cite{ding_NumericalSimulationFlows_2007} used a least-squares meshfree finite-difference method. Under the fine-resolution condition, the relative differences in $C_{\rm D}$, $C_{\rm L}$, and $St$ from the results of Liu \textit{et al.} were 3.5\%, $-0.29$\%, and 3.0\%, respectively, at $Re=100$, and 5.1\%, $-1.4$\%, and 4.7\%, respectively, at $Re=200$. Thus, the proposed method provided fluid-force and vortex-shedding characteristics in good agreement with the benchmark results obtained using a body-fitted mesh. Under the coarse-resolution condition, the relative difference in $C_{\rm L}$ from the result of Liu \textit{et al.} was 0.59\% at $Re=100$ but increased to $-13$\% at $Re=200$. As discussed in Section~\ref{subsec_bd_treat_stg}, the MLS formulation approximates the normal-velocity boundary condition by neglecting the skew-symmetric term $\nabla\mathbf{u}-(\nabla\mathbf{u})^{\top}$, which introduces an error of $\mathcal{O}(\|\mathbf{e}_{ij}\|)$. At low spatial resolution, this boundary-discretization error may have a non-negligible influence on the near-wall flow as the Reynolds number and associated vorticity increase. The improved agreement obtained with the fine resolution indicates that this influence is reduced by spatial refinement.

\begin{table}[t]
    \centering
    \caption{Mean values and fluctuation amplitudes of $C_{\rm D}$ and $C_{\rm L}$, and Strouhal number $St$ for the flow around a stationary cylinder at $Re=100$ and 200. ``coarse'' is the result with $l_0=0.05$ and $\Delta t=0.01$, and ``fine'' is the result with $l_0=0.025$ and $\Delta t=0.005$. Results marked with ``*'' are taken from experimental studies.}
    \label{tab:cd_cl_st_Karman}
    \begin{tabular}{lllll}
        \hline
        ~~ & ~~ & $C_{\rm D}$ & $C_{\rm L}$ & $St$ \\
        \hline
        $Re=100$~~ & Prandtl \textit{et al.}* \cite{prandtl_ErgebnisseAerodynamischenVersuchsanstalt_2009}~~ & 1.43                & --            & --    \\
        ~~         & Williamson* \cite{williamson_ObliqueParallelModes_1989}~~                               & --                  & --            & 0.164 \\
        ~~         & Liu \textit{et al.} \cite{liu_PreconditionedMultigridMethods_1998}~~                    & $1.350 \pm 0.012$~~ & $\pm 0.339$~~ & 0.165 \\
        ~~         & Ding \textit{et al.} \cite{ding_NumericalSimulationFlows_2007}~~                        & $1.356 \pm 0.010$~~ & $\pm 0.287$~~ & 0.166 \\
        ~~         & Present (coarse)~~                                                                      & $1.418 \pm 0.012$~~ & $\pm 0.341$~~ & 0.168 \\
        ~~         & Present (fine)~~                                                                        & $1.397 \pm 0.010$~~ & $\pm 0.338$~~ & 0.170 \\
        ~~ & ~~ & ~~ & ~~ & ~~ \\
        $Re=200$~~ & Prandtl \textit{et al.}* \cite{prandtl_ErgebnisseAerodynamischenVersuchsanstalt_2009}~~ & 1.29                & --            & --    \\
        ~~         & Williamson* \cite{williamson_ObliqueParallelModes_1989}~~                               & --                  & --            & 0.197 \\
        ~~         & Liu \textit{et al.} \cite{liu_PreconditionedMultigridMethods_1998}~~                    & $1.31 \pm 0.049$~~  & $\pm 0.69$~~  & 0.192 \\
        ~~         & Ding \textit{et al.} \cite{ding_NumericalSimulationFlows_2007}~~                        & $1.348 \pm 0.050$~~ & $\pm 0.659$~~ & 0.196 \\
        ~~         & Present (coarse)~~                                                                      & $1.391 \pm 0.040$~~ & $\pm 0.603$~~ & 0.195 \\
        ~~         & Present (fine)~~                                                                        & $1.377 \pm 0.049$~~ & $\pm 0.680$~~ & 0.201 \\
        \hline
    \end{tabular}
\end{table}

\subsection{Periodic Taylor--Green vortex} \label{subsec_TGV}
\begin{figure} 
  \begin{center}
    \includegraphics[width=\linewidth]{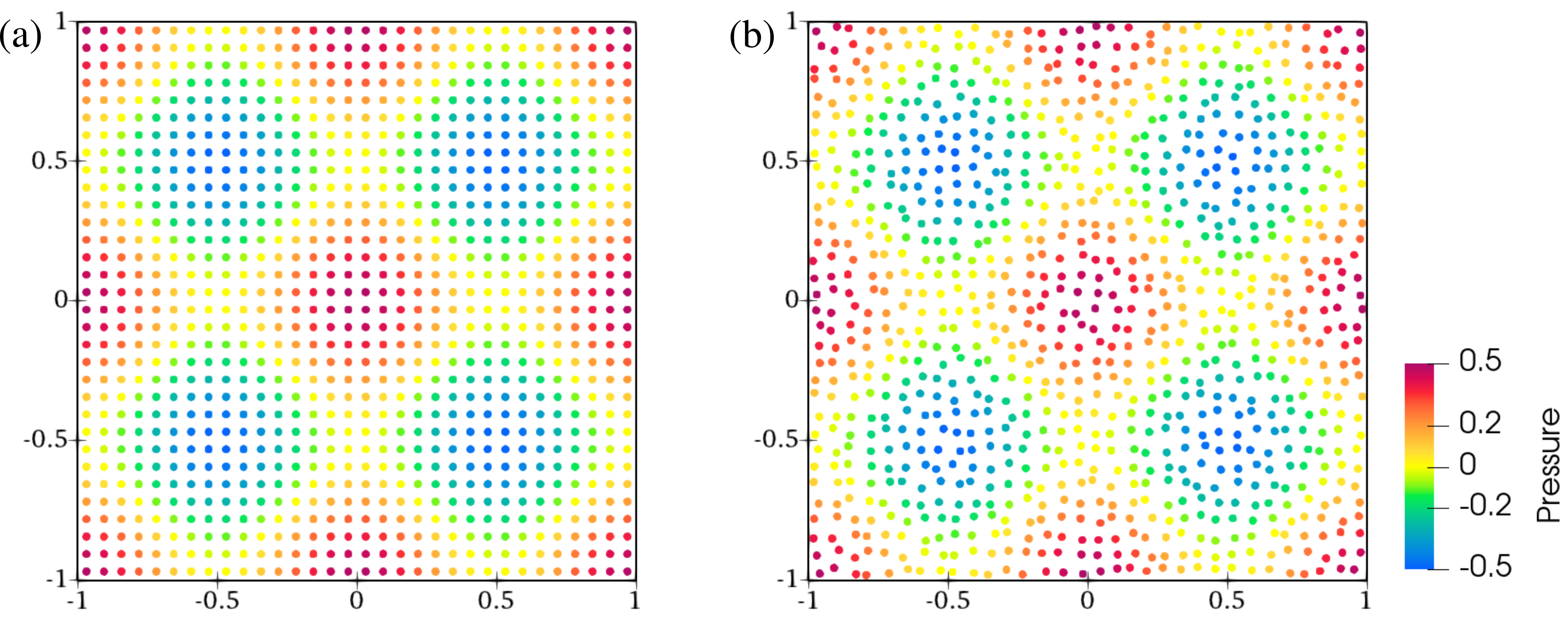}
    \caption{%
      Computational domain of the Taylor--Green vortex ($N = 32$) for (a) uniform DP arrangement and (b) randomized DP arrangement. Periodic boundary conditions are imposed on all boundaries. The color map indicates the initial pressure field.
    }
    \label{fig_DP_arrng_TGV_N32}
  \end{center}
\end{figure}

In this section, we investigate the spatial convergence accuracy by solving the Taylor--Green vortex problem in the square domain $[-L,L]\times[-L,L]$. Periodic boundary conditions are imposed on all boundaries. The exact solutions of the $u$ and $v$ components of velocities and pressure at $t$ are 
\begin{align}
  u_e(\mathbf{x},\, t) &= U \exp{\left( -\frac{2 \pi^{2} \eta t}{\rho L^{2}} \right)} \sin{\left( \frac{\pi x}{L} \right)} \cos{\left( \frac{\pi y}{L} \right)}, \label{eq_u_exact_TGV} \\[6pt]
  v_e(\mathbf{x},\, t) &= -U \exp{\left( -\frac{2 \pi^{2} \eta t}{\rho L^{2}} \right)} \cos{\left( \frac{\pi x}{L} \right)} \sin{\left( \frac{\pi y}{L} \right)}, \label{eq_v_exact_TGV} \\[6pt]
  p_e(\mathbf{x},\, t) &= \frac{\rho U^{2}}{4} \exp{\left( -\frac{4 \pi^{2} \eta t}{\rho L^{2}} \right)} \left[ \cos{\left( \frac{\pi x}{L} \right)} + \cos{\left( \frac{\pi y}{L} \right)}\right], \label{eq_p_exact_TGV}
\end{align}
where $U$ denotes the reference velocity magnitude.
The initial velocity field is given by the exact solution at $t = 0$.
In this study, we set the parameters to $U=1$, $L = 1$, $\rho = 1$, and $\eta = 0.01$.
The Reynolds number, defined as $Re = \rho U L / \eta$, was set to 100.
The DP spacing (or width of a background mesh element defined as $l_{0} = 2L/N$, where $N$ is the number of DPs in each direction) was doubly increased from $N =$ 8 to 1024.
The time step was set to $\Delta t = 1 \times 10^{-5}$, which was independent of the spatial resolution, so the diffusion number was $C_{\rm diff} = \eta \Delta t / (\rho l_{0}^{2}) = 0.0262$ at $N = 1024$.
The linear system of the pressure Poisson equation was solved using the Bi-CGSTAB method with the Lagrangian multiplier to constrain the pressure field to have zero mean over the domain.
The iteration was terminated when the absolute $L_{2}$ residual norm fell below the convergence criterion $\varepsilon_{\rm crit} = 10^{-7}$.
We calculated the temporal change in fluid flow for the uniform DP arrangement [Fig. \ref{fig_DP_arrng_TGV_N32}(a)] and the randomized DP arrangement [Fig. \ref{fig_DP_arrng_TGV_N32}(b)].
The random DP arrangement was generated by shifting each DP using random displacement vector $\delta \mathbf{x} = (\delta x,\, \delta y)^{\top}$ from the center of the background mesh element, where $\delta x$ and $\delta y$ are independently and identically distributed according to the uniform distribution $\mathcal{U}(-\alpha_{\rm rdm} l_{0}/2,\,\alpha_{\rm rdm} l_{0}/2)$.
The amplification coefficient $\alpha_{\rm rdm}$ was set to 0.5.

The numerical errors in $L_1$, $L_2$, and $L_\infty$ norms are defined as 
\begin{align}
  L_{1} &= \frac{1}{N}\sum_{\mathbf{x}_{i} \, \in \, \Omega}{\left| \mathbf{\phi}(\mathbf{x}_{i})-\mathbf{\phi}_{e}(\mathbf{x}_{i}) \right|}, 
  \label{eq_L1_TGV} \\[6pt]
  L_{2} &= \sqrt{\frac{1}{N}\sum_{\mathbf{x}_{i} \, \in \, \Omega}{\left| \mathbf{\phi}(\mathbf{x}_{i})-\mathbf{\phi}_{e}(\mathbf{x}_{i}) \right|^{2}}}, 
  \label{eq_L2_TGV} \\[6pt]
  L_{\infty} &= \max_{\mathbf{x}_{i} \, \in \, \Omega}{\left| \mathbf{\phi}(\mathbf{x}_{i})-\mathbf{\phi}_{e}(\mathbf{x}_{i}) \right|},
  \label{eq_Linf_TGV}
\end{align}
where $\phi \in \{u,\,v,\,p\}$ is the numerical solution calculated using the proposed method, and $\phi_{e}$ is the corresponding exact solution ($u_{e},\,v_{e},\,p_{e}$).

\begin{figure} 
  \begin{center}
    \includegraphics[width=\linewidth]{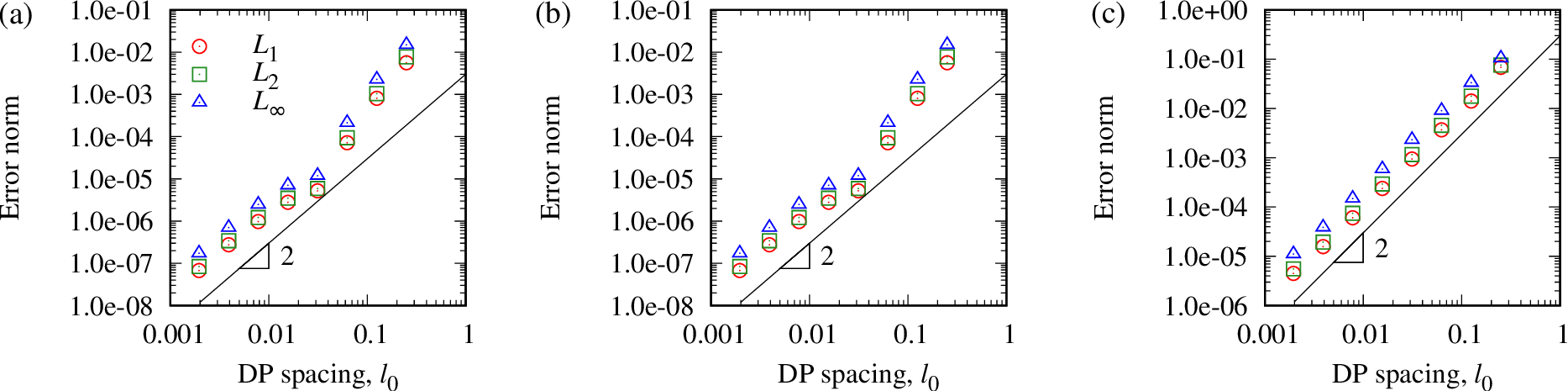}
    \caption{%
      Spatial convergence of $L_{1}$, $L_{2}$, and $L_{\infty}$ error norms at $t = 0.1$ for (a) the $x$ component of velocity, (b) the $y$ component of velocity, and (c) the pressure for a uniform DP arrangement.
    }
    \label{fig_err_TGV_uniform}
  \end{center}
\end{figure}

\begin{figure} 
  \begin{center}
    \includegraphics[width=\linewidth]{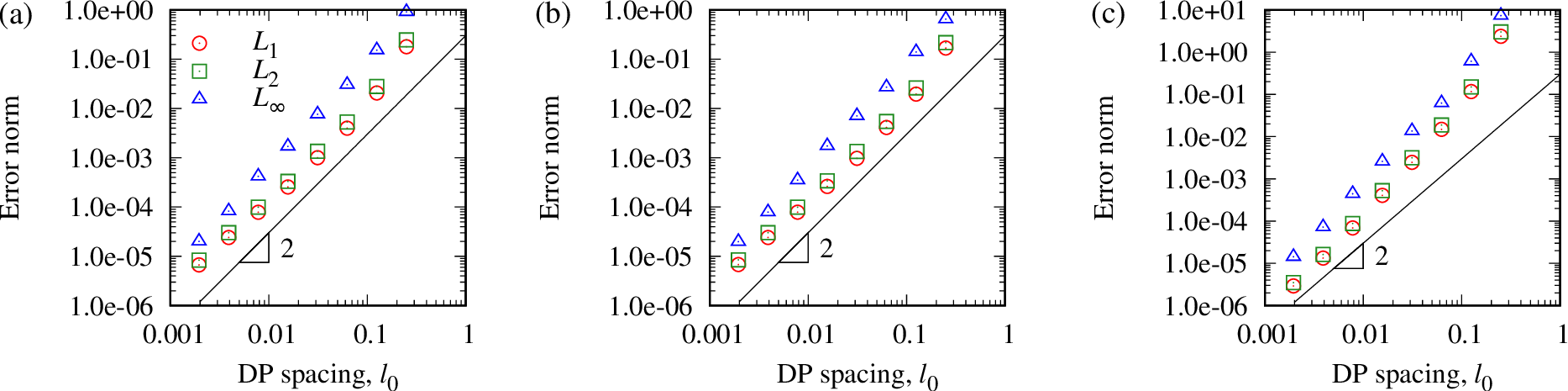}
    \caption{%
      Spatial convergence of $L_{1}$, $L_{2}$, and $L_{\infty}$ error norms at $t = 0.1$ for (a) the $x$ component of velocity, (b) the $y$ component of velocity, and (c) the pressure for a randomized DP arrangement.
    }
    \label{fig_err_TGV_random}
  \end{center}
\end{figure}

Figs. \ref{fig_err_TGV_uniform} and \ref{fig_err_TGV_random} show spatial convergences of $L_{1}$, $L_{2}$, and $L_{\infty}$ error norms at $t = 0.1$ for (a) $u$, (b) $v$, and (c) $p$ for the uniform and randomized DP arrangements, respectively.
Moreover, Table \ref{tab:convergence_orders} summarizes the convergence orders estimated from the numerical results at spatial resolutions, $N=512$ and 1024, as $L_{q-{\rm Order}}^{(1024)} = \log_{2}{\left( L_{q}^{(512)} / L_{q}^{(1024)} \right)}$.
These results indicate that the proposed second-order MLS reconstruction shows approximately second-order spatial convergence accuracy for velocity and pressure in both the uniform and random DP arrangements.


\begin{table}[t]
    \centering
    \caption{Observed convergence orders for $u$, $v$, and $p$, estimated from the numerical results at spatial resolutions, $N=512$ and 1024.}
    \label{tab:convergence_orders}
    \begin{tabular}{lllll}
        \hline
        ~~ & ~~ & $L_{1-{\rm Order}}$ & $L_{2-{\rm Order}}$ & $L_{\infty-{\rm Order}}$ \\
        \hline
        Uniform~~
        & $u$~~ & 2.02~~ & 2.02~~ & 2.01 \\
        & $v$~~ & 2.02~~ & 2.02~~ & 2.01 \\
        & $p$~~ & 1.80~~ & 1.80~~ & 1.80 \\
        ~~ & ~~ & ~~ & ~~ & ~~ \\
        Randomized~~
        & $u$~~ & 1.85~~ & 1.85~~ & 2.04 \\
        & $v$~~ & 1.83~~ & 1.83~~ & 1.99 \\
        & $p$~~ & 2.19~~ & 2.21~~ & 2.34 \\
        \hline
    \end{tabular}
\end{table}

\subsection{Lid-driven cavity flow} \label{subsec_cavity}
We calculated the two-dimensional lid-driven cavity flow to demonstrate the applicability of the proposed method over a wide range of Reynolds numbers.
This problem is a classical benchmark for incompressible Navier--Stokes solvers and has been validated by many researchers \cite{ghia_HighReSolutionsIncompressible_1982,hou_SimulationCavityFlow_1995,peng_Transition2DLiddriven_2003,li_FlowStructuresTwodimensional_2025}.
In this study, the flow was calculated in a $[-L/2,\,L/2] \times [-L/2,\,L/2]$ square domain, where $L$ is the characteristic length set to $L = 1$.
The boundary condition $\mathbf{u} = (U,\,0)^{\top}$ for velocity was imposed on the upper wall, which drives the flow. The quantity $U$ is the characteristic flow speed, which was set to $U = 1$.
The other three walls were subject to the no-slip boundary condition ($\mathbf{u} = \mathbf{0}$).
We set the fluid density to $\rho = 1$, and set the viscosity $\eta$ by setting the Reynolds number $Re = \rho UL / \eta$ as $Re = 100,\, 1000$, and $5000$.
The square domain was discretized using $N \times N$ DPs.
The flows at $Re = 100$ and $1000$ were calculated with $N = 33$, $65$, and $129$, while the flows at $Re = 5000$ were calculated with $N = 65$, $129$, and $257$.
The time interval was set such that the Courant number was equal to 0.128.
The linear system of the pressure Poisson equation was solved using Bi-CGSTAB with the Lagrangian multiplier constraint that the mean pressure over the domain was 0.
Bi-CGSTAB iterations were stopped when the relative $L_{2}$ residual norm was less than the convergence criterion $\varepsilon_{\rm crit} = 10^{-7}$.
The flow was simulated until the maximum temporal variation rate of each velocity component was less than $10^{-3}$.

\begin{figure} 
  \begin{center}
    \includegraphics[width=\linewidth]{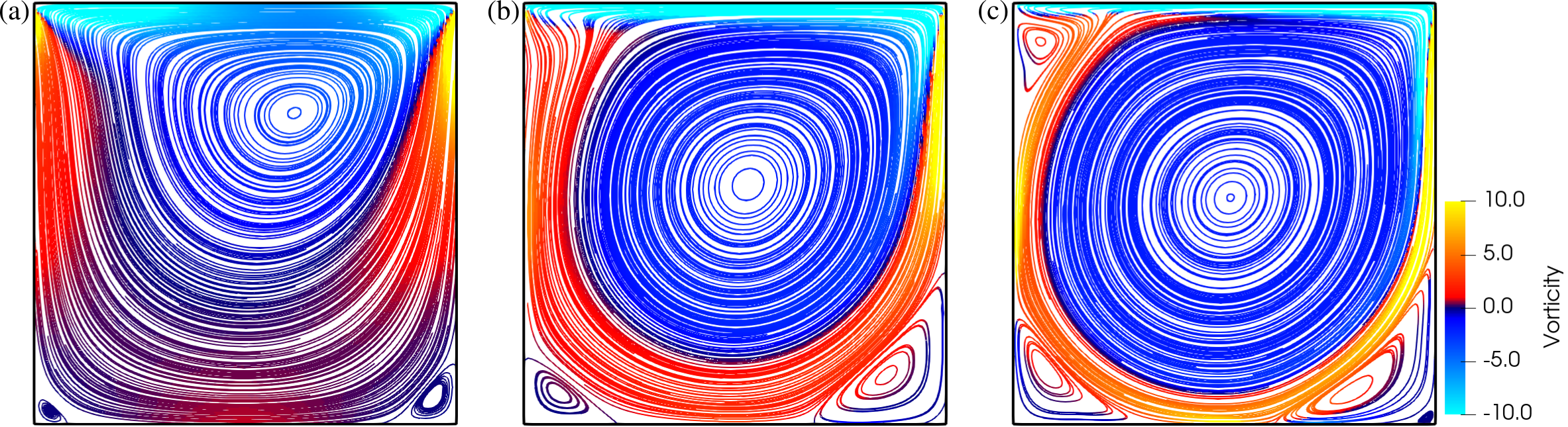}
    \caption{%
      Streamline visualizations of the lid-driven cavity flow calculated using the proposed staggered MCD method at (a) $Re=100$ with $129 \times 129$ (snapshot at $t = 9.322$), (b) $Re = 1000$ with $129 \times 129$ ($t = 34.229$), and (c) $Re = 5000$ with $257 \times 257$ ($t = 176.6615$).
      The vorticity field is visualized with contours limited to the range $[-10,\, 10]$ to highlight the vortex structures in the cavity interior, while higher magnitudes near the walls are omitted.
    }
    \label{fig_stream_line_cavity}
  \end{center}
\end{figure}

\begin{figure} 
  \begin{center}
    \includegraphics[width=\linewidth]{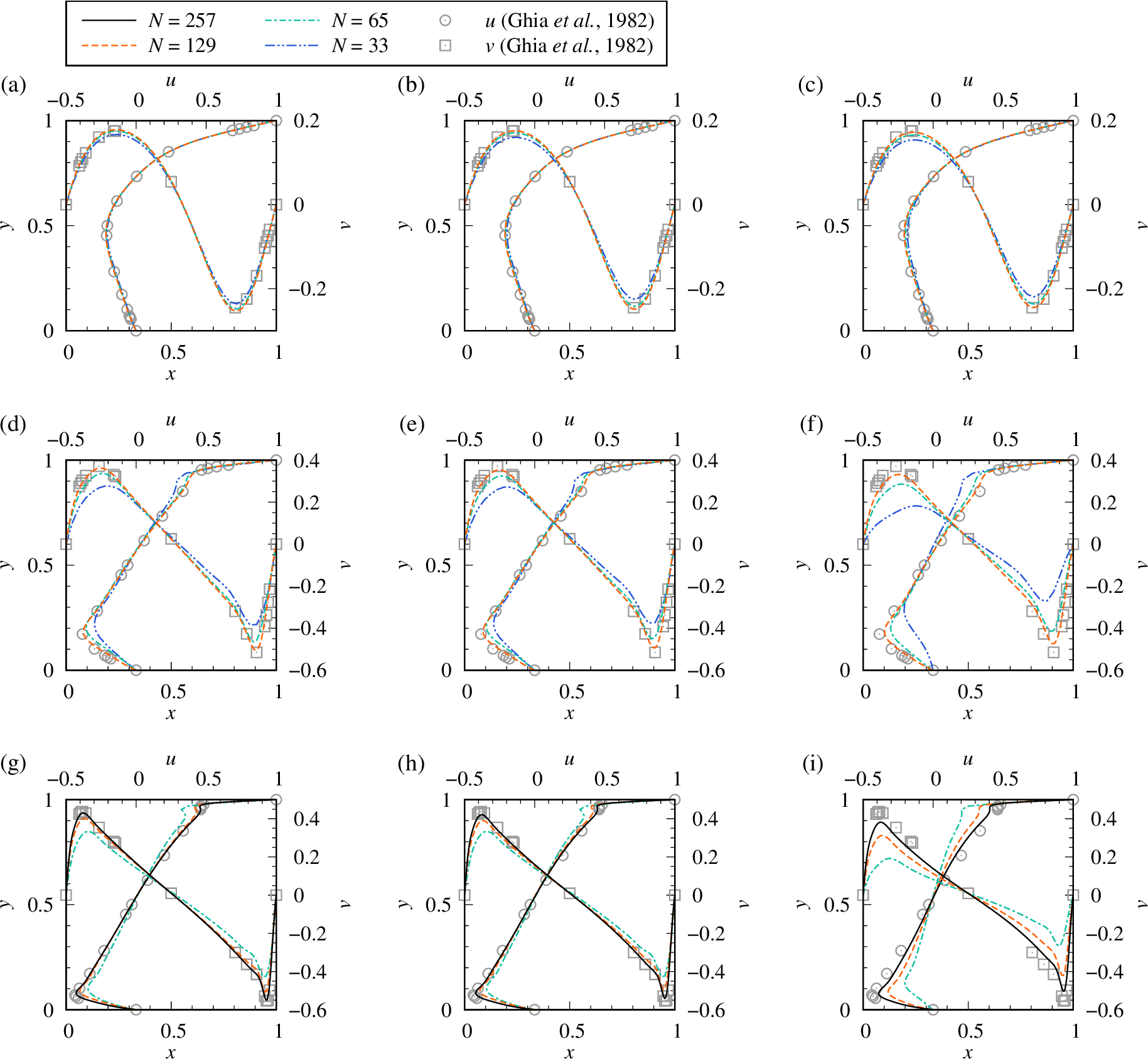}
    \caption{%
      Centerline velocity profiles of the lid-driven cavity flow compared with the reference data of Ghia \textit{et al.} \cite{ghia_HighReSolutionsIncompressible_1982} for (a--c) $Re = 100$, (d--f) $Re = 1000$, and (g--i) $Re = 5000$.
      Panels (a, d, g) show the proposed staggered method with $\beta_{\rm TEC} = 0.99$, (b, e, h) use $\beta_{\rm TEC} = 0$, and (c, f, i) show the previous collocated MCD method.
      Black solid, orange dashed, green dash-dot, and blue dash-dot-dot curves represent $N =$ 257, 129, 65, and 33, respectively.
      Gray circles and squares are the reference data for $u$ and $v$, respectively.
    }
    \label{fig_uv_cavity}
  \end{center}
\end{figure}

Fig. \ref{fig_stream_line_cavity} shows the streamlines of the cavity flows at (a) $Re = 100$, (b) $Re = 1000$, and (c) $Re = 5000$ at the final timestep.
In these snapshots, the primary vortex becomes larger and shifts toward the center of the square domain as the Reynolds number increases.
Secondary vortices developed in the lower left and right corners.
At $Re = 5000$, an additional secondary vortex appeared in the upper left corner, and a tertiary vortex formed in the lower right corner.
These tendencies are in good agreement with previously reported numerical results \cite{ghia_HighReSolutionsIncompressible_1982,hou_SimulationCavityFlow_1995,li_FlowStructuresTwodimensional_2025}.
Fig. \ref{fig_uv_cavity} shows the velocity profiles along the horizontal and vertical centerlines of the square domain.
The proposed staggered approach obtains a more accurate solution than that of the conventional collocated formulation as the Reynolds number increases.
At $Re = 5000$, the results of the staggered approach at $N = 129$ [Fig. \ref{fig_uv_cavity}(g, h)] are closer to the reference data of Ghia \textit{et al.} \cite{ghia_HighReSolutionsIncompressible_1982} than those of the $N = 257$ solution of the collocated method [Fig. \ref{fig_uv_cavity}(i)].
These convincingly show that the consistency of the discrete gradient and divergence operators (in satisfying incompressibility at the discrete level) is important for simulating the higher-Reynolds-number flow in the projection-based meshfree framework.
Moreover, comparing $\beta = 0.99$ and $\beta = 0$ in the staggered method indicates that the TEC formula provides more accurate flow characteristics for all Reynolds numbers.

To evaluate the computational performance of the staggered and collocated formulations, their computational times were compared. The simulations were performed on a workstation equipped with an AMD Ryzen 9 9950X 16-core processor. For the flow with $257 \times 257$ DPs and $Re=5000$, the average computational times per time step for the staggered ($\beta = 0.99$) and collocated formulations were $0.120$ s and $0.109$ s, respectively, and the increase due to the staggered formulation was approximately $10\%$. The average numbers of iterations for the pressure Poisson equation were $565$ and $529$, respectively. Therefore, even including the additional evaluation of the radial velocity components in the proposed method, the increase in computational time relative to the collocated formulation was limited.


\subsection{Flow around a periodically oscillating cylinder} \label{subsec_cyl_oscill}
As a final test case, we calculated the flow around a periodically oscillating cylinder of diameter $D$.
Although this problem has been used as a benchmark to validate a flow solver in moving boundary problems \cite{dutsch_LowReynoldsnumberFlowOscillating_1998,guilmineau_NUMERICALSIMULATIONVORTEX_2002,chi_DirectionalGhostcellImmersed_2020,ghomizad_SharpInterfaceDirectforcing_2021,matsuda_MeshconstrainedDiscretePoint_2025}, in this study, the flow was calculated in a coordinate system $\tilde{\mathbf{x}}$ fixed to the moving cylinder to eliminate the effect of updating the DP arrangement around the moving cylinder and to validate the flow solver itself through comparisons with both available experimental and numerical results.
According to the Galilean transformation, the moving coordinate system $\tilde{\mathbf{x}}$ was defined as $\tilde{\mathbf{x}} = \mathbf{x} - \mathbf{x}_{\rm cyl}$, where $\mathbf{x}_{\rm cyl} = (x_{\rm cyl},\,y_{\rm cyl})^{\top}$ is the cylinder center.
In the coordinate system fixed to the moving cylinder, the momentum equation can be written as
\begin{equation}
  \frac{\partial (\rho \mathbf{u})}{\partial t} + \nabla \cdot (\rho \mathbf{u} \mathbf{u}) = -\nabla p + \nabla \cdot \boldsymbol{\tau} - \rho \ddot{\mathbf{x}}_{\rm cyl}.
  \label{eq_NS_single_cyl_oscill_fix}
\end{equation}
In this study, the cylinder motion was given as $x_{\rm cyl} = -A \sin{(2 \pi f t)}$, $y_{\rm cyl} = 0$, where $A = U / 2 \pi f$ denotes the amplitude coefficient of the cylinder motion, and $U$ and $f$ are the characteristic velocity magnitude of the cylinder and the frequency.
The Reynolds number $Re$ and Keulegan--Carpenter number $KC$ of the flow were defined as $Re = \rho UD / \eta$ and $KC = U / fD$, and were set to $Re = 100$ and $KC = 5$.
Parameters were set to $D = 1$, $f = 1$, $U = 5$, $\rho = 1$, and $\eta = 5 \times 10^{-2}$.
The side length of the computational domain was set to $L = 20D$.
The width of the background mesh was set to $l_{0} = 2.5 \times 10^{-2}$, and the time interval to $\Delta t = 1/900$.
The no-slip boundary conditions were imposed on all the walls.
We calculated the linear system of the pressure Poisson equation using Bi-CGSTAB with the Lagrangian multiplier constraint that the mean pressure is zero.
In this problem, the iteration for solving pressure was terminated when the relative $L_{2}$ error norm of Bi-CGSTAB fell below the convergence criterion $\varepsilon_{\rm crit} = 10^{-5}$.
The TEC parameter $\beta$ was set to 0.99.

\begin{figure} 
  \begin{center}
    \includegraphics[width=\linewidth]{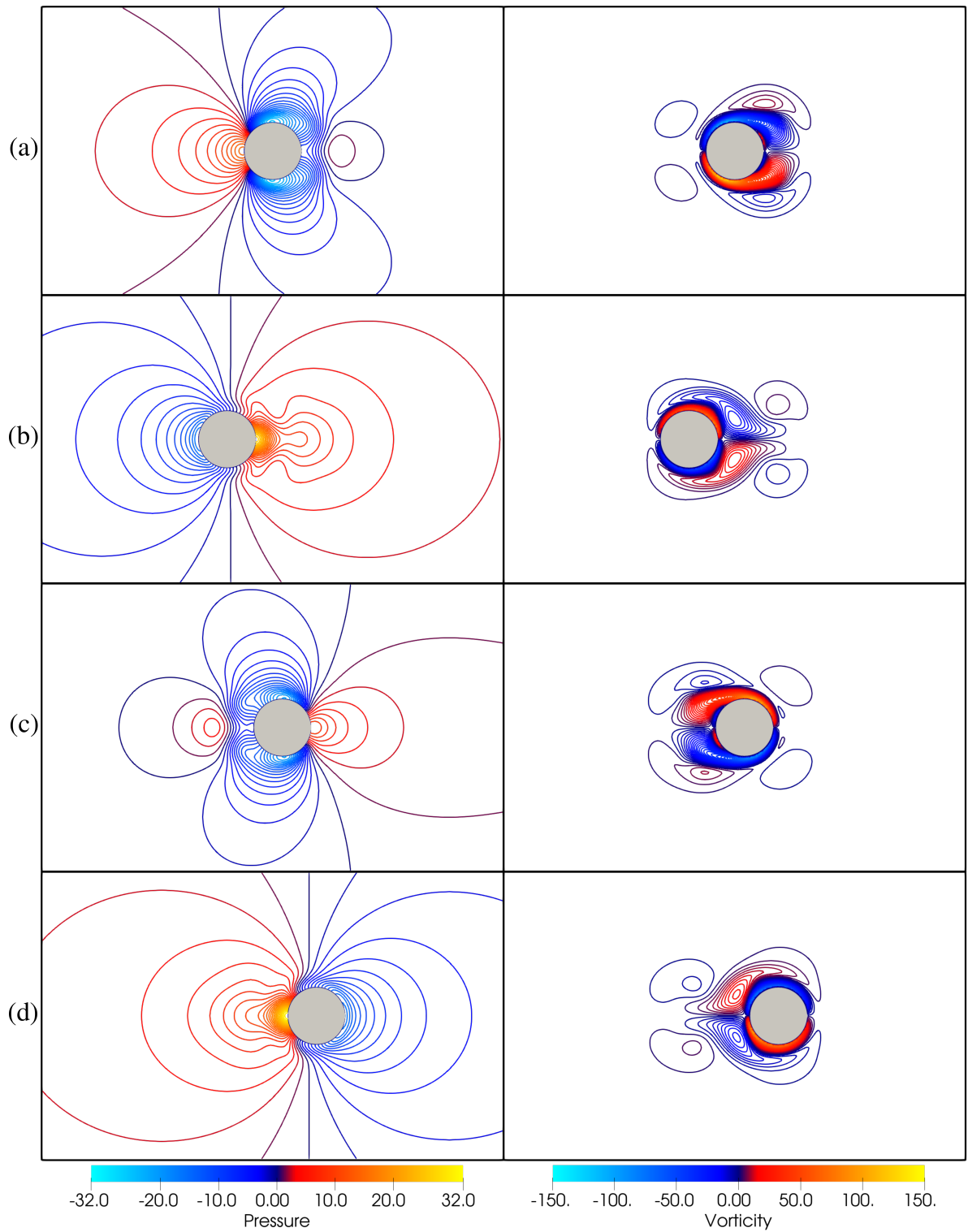}
    \caption{%
      Numerical results for the pressure (left) and vorticity (right) isolines around an oscillating cylinder. The phase angles of the cylinder position are (a) 0\textdegree, (b) 96\textdegree, (c) 192\textdegree, and (d) 288\textdegree.
      The pressure contour levels are \textminus 32 to 32 in increments of 1.28, and the vorticity contour levels are \textminus 150 to 150 in increments of 2.4.
    }
    \label{fig_contour_single_cyl_oscill_fix}
  \end{center}
\end{figure}

\begin{figure} 
  \begin{center}
    \includegraphics[width=\linewidth]{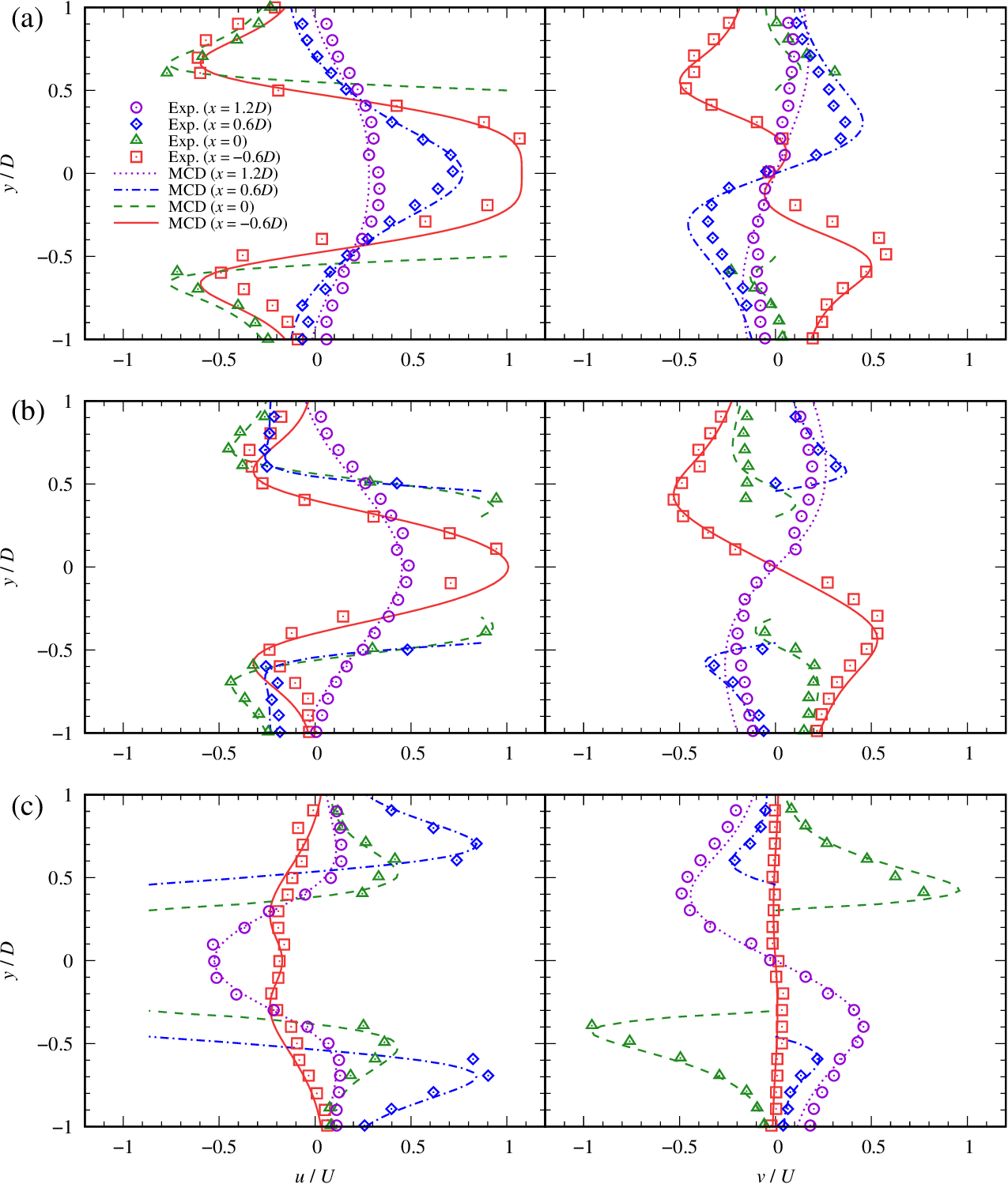}
    \caption{%
      Comparisons of $x$ and $y$ components of velocities for four cross-sections: $x = -0.6D$ (red curve/squares), $x = 0$ (green dashed curve/triangles), $x = 0.6D$ (blue dash-dotted curve/diamonds), and $x = 1.2D$ (purple dotted curve/circles).
      The phase angles of the cylinder position are (a) 180\textdegree, (b) 210\textdegree, and (c) 330\textdegree. Experimental results are taken from \cite{dutsch_LowReynoldsnumberFlowOscillating_1998}.
    }
    \label{fig_uv_profile_single_cyl_oscill_fix}
  \end{center}
\end{figure}

\begin{figure} 
  \begin{center}
    \includegraphics[width=.65\linewidth]{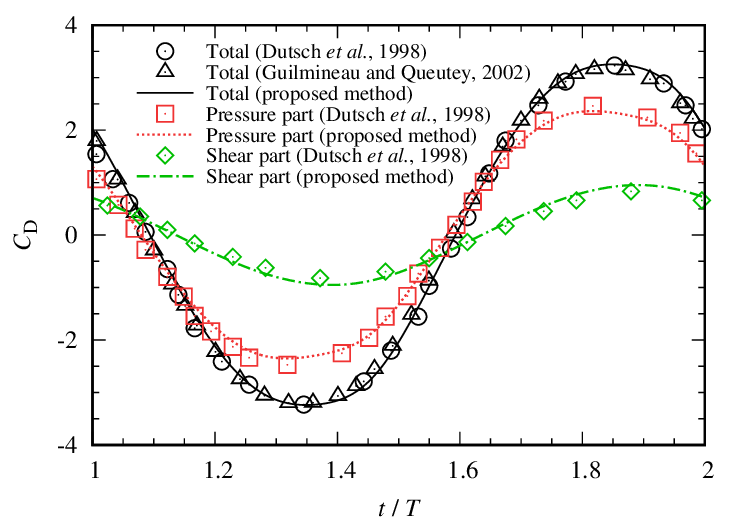}
    \caption{%
      Comparison of the time variation in the drag force coefficient and its components.
      $T = 1/f$ is the period of cylinder motion.
      The circles/triangles/squares/diamonds indicate previous results \cite{dutsch_LowReynoldsnumberFlowOscillating_1998,guilmineau_NUMERICALSIMULATIONVORTEX_2002}.
    }
    \label{fig_cd_single_cyl_oscill_fix}
  \end{center}
\end{figure}

Fig. \ref{fig_contour_single_cyl_oscill_fix} shows the contours of the pressure and vorticity around the cylinder at four phase angles calculated with the proposed method.
The results show the asymmetric pressure distribution between the upstream and downstream sides of the cylinder caused by inertial effects, as well as the periodic vortex shedding in the wake behind the cylinder.
Both the pressure and vorticity are in good agreement with the well-validated numerical solution of D\"utsch \textit{et al.} using a fixed body-fitted mesh \cite{dutsch_LowReynoldsnumberFlowOscillating_1998}.
Fig. \ref{fig_uv_profile_single_cyl_oscill_fix} shows the comparisons of the $x$ and $y$ components of the velocity profiles with published experimental results \cite{dutsch_LowReynoldsnumberFlowOscillating_1998} at four cross-sections and three phase angles.
These graphs show that the velocity calculated with the proposed method is in good agreement with both the experimental and numerical results.
Fig. \ref{fig_cd_single_cyl_oscill_fix} shows the temporal change in the drag coefficient $C_{\rm D}$ and its components (pressure and viscosity), and the results are in good agreement with the reference numerical solutions \cite{dutsch_LowReynoldsnumberFlowOscillating_1998,guilmineau_NUMERICALSIMULATIONVORTEX_2002}.


\section{Conclusion} \label{sec_conclusion}
This study proposed a least-squares strong-form meshless method for the incompressible Navier--Stokes equations.
The proposed method adopts MLS reconstruction with radial components defined at staggered points (i.e., midpoints between each DP and its neighbors) for the spatial discretization.
The proposed method consistently evaluates the discrete gradient and divergence operators, thereby guaranteeing flow incompressibility consistent with the error norm of the pressure Poisson equation when the linear system converges.
Achieving such consistency has been inherently challenging in strong-form collocated least-squares meshfree frameworks employing projection-based velocity--pressure coupling.
Moreover, in the proposed method, the radial velocity components with a nonlinear advection term are evaluated on staggered points using the TEC formula, which links the temporal variation between the velocity and pressure.
This achieves a consistent velocity--pressure coupling and suppresses the odd-even instability causing the checkerboard pressure instability commonly observed in collocated methods.
Several numerical tests indicate that the proposed method achieves approximately second-order spatial convergence accuracy for both velocity and pressure, as expected from the formulation, and reproduces incompressible flows more accurately than the conventional collocated method over a wide range of Reynolds numbers.
The formulation is implemented within the MCD framework using compact 3 $\times$ 3 stencils and demonstrates accurate calculations.
These results demonstrate that the proposed method offers a promising direction for further developing stable, accurate, and high-fidelity incompressible flow solvers in strong-form meshfree frameworks.

There are several limitations. First, the radial components of the velocity evaluated at staggered points must be interpolated from the quantities defined at DPs. Therefore, dispersive and dissipative errors are unavoidable due to the interpolation/reconstruction procedure. An investigation into stabilization techniques and higher-order reconstruction would be necessary to address these issues. Second, the proposed method employs a second-order upwind scheme based on the sign of the flux evaluated on each edge. This formulation has been confirmed to reproduce the lid-driven cavity flow well for Reynolds numbers up to $Re \sim O(10^3)$, showing good agreement with existing numerical results. However, at higher Reynolds numbers ($O(10^4)$ or larger), the numerical viscosity may become excessive, which could reduce the fidelity of the flow reproduction. A higher-order scheme for the convective term could improve the accuracy. Third, although the present formulation provides stable velocity--pressure coupling in practical computations, it is not guaranteed that the formulation can still suppress checkerboard-type pressure instability when the time step becomes extremely small. It should be investigated whether the TEC formulae or other more rigorous formulations for velocity--pressure coupling would remain effective under such extreme conditions.

\section*{Acknowledgments}
The authors thank Dr. Hiro Wakimura (Institute of Science Tokyo) for his constructive comments on the numerical scheme. This work was supported by JST SPRING Grant Number JPMJSP2156; JSPS KAKENHI Grant Numbers JP24KJ1851, JP22H00190, JP24K02557, JP26K00867, and JP26KJ0221; the MEXT Program for Promoting Researches on the Supercomputer Fugaku (Development of human digital twins for cerebral circulation using Fugaku, Grant Number JPMXP1020230118). This study used computational resources of the supercomputer Fugaku provided by the RIKEN Center for Computational Science (project IDs: hp230208, hp240220, hp250236, hp240080, hp240458, and hp250282). Some computations were also carried out using the supercomputer Flow at the Information Technology Center, Nagoya University. The work was partly supported by the Research Organization for Information Science \& Technology (RIST) under the HPCI User Support Program. Mark R. Kurban from Edanz (https://www.jp.edanz.com/ac) edited a draft of this paper.

%
\section*{CRediT authorship contribution statement}
\textbf{Takeharu Matsuda:} Writing--original draft, Visualization, Validation, Software, Methodology, Investigation, Data curation, Conceptualization. \textbf{Satoshi Ii:} Writing--review \& editing, Validation, Software, Methodology, Conceptualization.

\section*{Declaration of competing interest}
The authors declare that they have no known competing financial interests or personal relationships that could have appeared to influence the work reported in this paper.

\section*{Data Availability Statement}
Data will be made available on request.

%
\begin{appendix}
\renewcommand{\thesection}{Appendix \Alph{section}}
\renewcommand{\theequation}{\Alph{section}.\arabic{equation}}  
\setcounter{equation}{0}  
\renewcommand{\thefigure}{\Alph{section}.\arabic{figure}}  
\setcounter{figure}{0}  

\section{MLS reconstruction using radial components} \label{appendix_MLS_stg}
\setcounter{equation}{0}

As in \cite{trask_HighorderStaggeredMeshless_2017}, a vector variable $\mathbf{u}(\mathbf{x})$ in $\mathbf{x} \in \mathbb{R}^d$, where $d$ is the dimension number ($d=1, 2, 3$), is converted to a radial component from a point $\mathbf{x}_{c}$, $U_{c \rightarrow}(\mathbf{x}) \in \mathbb{R}^{1}$:
\begin{equation}
    U_{c \rightarrow}(\mathbf{x})
    = 2(\mathbf{x}-\mathbf{x}_{c}) \cdot \mathbf{u}(\mathbf{x}).
    \label{eq_radial_U}
\end{equation}
Here, the following relations hold:
\begin{align}
    \nabla U_{c \rightarrow}(\mathbf{x})
    &= 2\mathbf{u} (\mathbf{x})
    + 2 \nabla \mathbf{u} (\mathbf{x}) \cdot (\mathbf{x} - \mathbf{x}_{c}), \label{eq_dU} \\
    \rightarrow \quad
    \nabla U_{c \rightarrow} (\mathbf{x}_{c})
    &= 2\mathbf{u} (\mathbf{x}_{c}), \label{eq_dU_at_xc} \\
    \nabla \nabla U_{c \rightarrow} (\mathbf{x})
    &= 2\nabla \mathbf{u} (\mathbf{x})
    + 2\left( \nabla \mathbf{u} (\mathbf{x}) \right)^{\top}
    + 2\nabla \nabla \mathbf{u} (\mathbf{x}) \cdot (\mathbf{x} - \mathbf{x}_{c}), \label{eq_ddU} \\
    \rightarrow \quad
    \nabla\nabla U_{c \rightarrow}(\mathbf{x}_{c})
    &= 2\nabla \mathbf{u} (\mathbf{x}_{c}) + 2\left( \nabla \mathbf{u} (\mathbf{x}_{c}) \right)^{\top}, \label{eq_ddU_at_xc}
\end{align}
and
\begin{align}
  \nabla^2 U_{c \rightarrow} (\mathbf{x})
  &= 4 \nabla \cdot \mathbf{u} (\mathbf{x})
  + 2 \nabla^2 \mathbf{u} (\mathbf{x}) \cdot (\mathbf{x} - \mathbf{x}_{c}), \label{eq_div_dU} \\
  \rightarrow \quad \nabla^{2} U_{c \rightarrow} (\mathbf{x}_{c})
  &= 4 \nabla \cdot \mathbf{u} (\mathbf{x}_{c}). \label{eq_div_dU_at_xc}
\end{align}

A local polynomial function $U^{h}_{c}(\mathbf{x})= {\bf c}_{c}\cdot{\bf p}_{c}(\mathbf{x})$ is introduced as an approximate function of $U_{c \rightarrow}(\mathbf{x})$, where the components of ${\bf c}_{c}$ are the coefficients and ${\bf p}_{c}(\mathbf{x})$ is the arbitrary order basis function, which is reconstructed by solving the optimization problem
\begin{equation}
    \argmin_{{\bf c}_{c}}
    \ \frac{1}{2} \sum_{\mathbf{x} \in \Omega_{c}}
    w_{c}(\mathbf{x})
    \left( 
        U^{h}_{c} (\mathbf{x}) - U_{c \rightarrow}(\mathbf{x})
    \right)^2,
\end{equation}
where $w_{c}(\mathbf{x})$ is the weight function from $\mathbf{x}_{c}$.

Thus, once the approximate function $U^{h}_{c}(\mathbf{x})$ is reconstructed, the spatial derivatives at $\mathbf{x}_{c}$ are approximated as
\begin{align}
  \mathbf{u} (\mathbf{x}_{c})
  &\approx \mathbf{u}^{h} (\mathbf{x}_{c}) 
  = \frac{1}{2} \nabla U^{h}_{c} (\mathbf{x}_{c}), \\
  \nabla \cdot \mathbf{u} (\mathbf{x}_{c})
  &\approx \nabla \cdot \mathbf{u}^{h} (\mathbf{x}_{c})
  = \frac{1}{4} \nabla^2 U^{h}_{c} (\mathbf{x}_{c}),
\end{align}
where $\mathbf{u}^{h}$ is the approximate function of $\mathbf{u}$.

\section{Boundary treatment for the incompressible Navier--Stokes equations} \label{sec_bc_radial_components}
\setcounter{equation}{0}

The boundary treatment used in the MLS reconstruction depends on both the physical boundary condition for the incompressible Navier--Stokes equations and the radial quantity to be reconstructed. Two types of physical boundaries are considered: a prescribed-velocity boundary, such as a no-slip wall or an inlet boundary, denoted by $\Gamma_{u}$, and a prescribed-pressure boundary, denoted by $\Gamma_{p}$. On these boundaries, the physical conditions are written as
\begin{align}
  \mathbf{u} 
  &= \mathbf{u}_{d}, 
  && \text{on } \Gamma_{u},
  \label{eq_bc_velocity_prescribed} \\
  -p\mathbf{n} + \eta \mathbf{n}\cdot\nabla\mathbf{u} 
  &= -p_{d}\mathbf{n}, 
  && \text{on } \Gamma_{p}. 
  \label{eq_bc_pressure_traction}
\end{align}
Since the pressure-projection method is employed, the condition in Eq.~\eqref{eq_bc_pressure_traction} is imposed separately on the pressure and viscous (or velocity-gradient) terms in the present formulation.
\begin{equation}
  p=p_d,
  \qquad \mathbf{n}\cdot\nabla\mathbf{u}=\mathbf{0},
  \qquad \text{on } \Gamma_p.
  \label{eq_bc_pressure_prescribed}
\end{equation}
Here, $p_{d}=0$ is used in the present simulations. Hereafter, $v_i$ denotes the node at which a differential operator is evaluated, whereas $v_j$ denotes a neighboring node located on either $\Gamma_u$ or $\Gamma_p$. The designation of $\Gamma_d$ in the boundary-constrained MLS reconstruction depends on the radial quantity being reconstructed: it represents the boundary on which the corresponding normal component is prescribed as in Eq. \eqref{eq_def_obj_func_with_bd_treatment}. Therefore, $\Gamma_d$ does not necessarily coincide with the same physical boundary for all of $U$, $\delta P$, $\mathbf{A}$, and $\mathbf{D}$.

\subsection*{Prescribed-velocity boundary $\Gamma_u$}
For $\mathbf{x}_j\in\Gamma_u$, the prescribed normal velocity
\begin{equation}
  q^{U}_{j}
  = (\mathbf{n}\cdot\mathbf{u}_{d})_{j},
  \qquad \text{on } \Gamma_u.
  \label{eq_bc_q_U_velocity}
\end{equation}
is imposed in the boundary-constrained MLS reconstruction of $U_i^h$ through Eq.~\eqref{eq_def_obj_func_with_bd_treatment}. In particular, the prescribed velocity vector is not directly substituted into the radial-component definition in Eq.~\eqref{eq_def_U_rad_cmpt}. Instead, only its normal component in Eq.~\eqref{eq_bc_q_U_velocity} is incorporated as the boundary constraint. This distinction is essential because the reconstructed $U_i^h$ is used to evaluate the velocity divergence, to which the velocity component normal to the physical boundary must contribute.

For the reconstruction of $\delta P_i^h$, the homogeneous normal pressure gradient is imposed as
\begin{equation}
  q^{\delta P}_{j}
  = (\mathbf{n}\cdot\nabla p)_{j}
  = 0,
  \qquad \text{on } \Gamma_u.
  \label{eq_bc_q_deltaP_velocity}
\end{equation}
Thus, both $U_i^h$ and $\delta P_i^h$ use the boundary-constrained MLS reconstruction on $\Gamma_u$, with the prescribed normal velocity and the homogeneous normal pressure gradient, respectively.

The advective radial component $\mathbf{A}$ is also incorporated through its normal component at the physical boundary. Using the prescribed velocity, its boundary value is
\begin{equation}
  \mathbf{q}^{A}_{j}
  = (\mathbf{n}\cdot\mathbf{u}_{d})_{j}(\mathbf{u}_{d})_{j}
  = q^{U}_{j}(\mathbf{u}_{d})_{j},
  \qquad \text{on } \Gamma_u.
  \label{eq_bc_q_A_velocity}
\end{equation}
This value is imposed in the boundary-constrained MLS reconstruction of $\mathbf{A}_i^h$ instead of evaluating $\mathbf{A}_{i\rightarrow ij}$ at the edge midpoint.

In contrast, $\mathbf{D}_{i\rightarrow ij}$ is evaluated at the edge midpoint using the prescribed boundary velocity as the value at $v_j$:
\begin{equation}
  \mathbf{D}_{i\rightarrow ij}
  = (\mathbf{u}_{d})_{j}-\mathbf{u}_{i},
  \qquad \text{for } \mathbf{x}_j\in\Gamma_u.
  \label{eq_bc_D_velocity}
\end{equation}
Accordingly, no additional normal-component constraint is imposed in the MLS reconstruction of $\mathbf{D}_i^h$ on $\Gamma_u$.

When an interpolated velocity associated with $v_j$ is required at $\mathbf{x}_{ij}$ in the TEC or upwind evaluation, the prescribed Dirichlet velocity is assigned by zeroth-order interpolation as
\begin{equation}
  \left.\mathbf{u}_{j}^{h}\right|_{ij}
  = (\mathbf{u}_{d})_{j},
  \qquad \text{for } \mathbf{x}_j\in\Gamma_u.
  \label{eq_bc_u_interp_velocity}
\end{equation}

\subsection*{Prescribed-pressure boundary $\Gamma_p$}
For $\mathbf{x}_j\in\Gamma_p$, the radial pressure increment $\delta P_{i\rightarrow ij}$ at the edge midpoint is evaluated using the prescribed boundary pressure at $v_j$:
\begin{equation}
  \delta P_{i\rightarrow ij}
  = (p_d)_{j}-p_{i},
  \qquad \text{for } \mathbf{x}_j\in\Gamma_p.
  \label{eq_bc_deltaP_pressure}
\end{equation}

For the advective radial component $\mathbf{A}$, the value at the edge midpoint $\mathbf{x}_{ij}$ is evaluated using the radial velocity component $U_{i\rightarrow ij}$ and the MLS approximation associated with $v_i$:
\begin{equation}
  \mathbf{A}_{i\rightarrow ij}
  = U_{i\rightarrow ij}
    \left.\mathbf{u}_{i}^{h}\right|_{ij},
  \qquad \text{for } \mathbf{x}_j\in\Gamma_p.
  \label{eq_bc_A_pressure}
\end{equation}
This boundary treatment assumes that no backflow occurs through $\Gamma_p$.

Similarly, when the two interface values appearing in the TEC formula \eqref{eq_def_TEC} are required, the value associated with $v_j$ is set equal to that evaluated using the MLS approximation associated with $v_i$:
\begin{equation}
  \left.\mathbf{U}_{j}^{h}\right|_{ij}
  = \left.\mathbf{U}_{i}^{h}\right|_{ij},
  \qquad \text{for } \mathbf{x}_j\in\Gamma_p.
  \label{eq_bc_TEC_pressure}
\end{equation}
Thus, the two interface values in the TEC formula coincide at the prescribed-pressure boundary, without altering the local support used in the MLS construction.

Finally, the homogeneous normal velocity-gradient condition in Eq.~\eqref{eq_bc_pressure_prescribed} is imposed in the boundary-constrained MLS reconstruction of $\mathbf{D}_i^h$ as
\begin{equation}
  \mathbf{q}^{D}_{j}
  = (\mathbf{n}\cdot\nabla\mathbf{u})_{j}
  = \mathbf{0},
  \qquad \text{on } \Gamma_p.
  \label{eq_bc_q_D_pressure}
\end{equation}

Note that the radial velocity $U_{i\rightarrow ij}$ at the edge midpoint for $\mathbf{x}_j\in\Gamma_p$ is updated through the TEC formula and the projection step.

\subsection*{Compatibility of the pressure projection at the boundaries}
The boundary treatments for $U_{i\rightarrow ij}$ and $\delta P_{i\rightarrow ij}$ are paired so as to preserve the discrete compatibility between the divergence operator used on the right-hand side of the pressure Poisson equation, the pressure Laplacian, and the pressure-gradient operator used in the projection step. On $\Gamma_u$, the prescribed normal velocity and the corresponding homogeneous normal pressure-gradient condition are incorporated through the boundary-constrained MLS reconstruction. On $\Gamma_p$, the midpoint value $U_{i\rightarrow ij}$ is projected using $\delta P_{i\rightarrow ij}=(p_d)_{j}-p_i$. Consequently, the pressure projection reduces the residual of the discrete velocity divergence to the level corresponding to the residual of the pressure Poisson equation, apart from the scaling factor $\Delta t/\rho$, as described in \ref{appendix_div_u_residual}.

\subsection*{MLS reconstruction}
The boundary conditions described above are incorporated into the MLS reconstructions of $U_i^h$, $\delta P_i^h$, $\mathbf{A}_i^h$, and $\mathbf{D}_i^h$. To explicitly describe this variable-dependent treatment, we first define the sets of intersection points associated with the prescribed-velocity and prescribed-pressure boundaries and their respective complements as
\begin{align}
  \mathcal{X}_{i}^{u}
  &=
  \left\{
    \mathbf{x}_{ij} \in \mathcal{X}_{i}
    \,\middle|\,
    \mathbf{x}_{j} \in \Gamma_{u}
  \right\},
  \quad \textrm{and} \quad
  \tilde{\mathcal{X}}^{u}_{i} 
  = \mathcal{X}_{i}\setminus\mathcal{X}_{i}^{u},
  \\
  \mathcal{X}_{i}^{p}
  &=
  \left\{
    \mathbf{x}_{ij} \in \mathcal{X}_{i}
    \,\middle|\,
    \mathbf{x}_{j} \in \Gamma_{p}
  \right\},
  \quad \textrm{and} \quad
  \tilde{\mathcal{X}}^{p}_{i}
  = \mathcal{X}_{i}\setminus\mathcal{X}_{i}^{p}.
\end{align}
The boundary-constrained MLS reconstructions of the radial quantities are then specified separately as follows.

\paragraph{Reconstruction of $U_i^h$}
At a prescribed-velocity boundary, the prescribed normal velocity $q_{j}^{U}=(\mathbf{n}\cdot\mathbf{u}_{d})_{j}$ is imposed through
\begin{equation}
\begin{split}
  \argmin_{\mathbf{c}_{i}^{U}}
  ~&\frac{1}{2}
  \sum_{\mathbf{x}_{ij}\in\tilde{\mathcal{X}}^{u}_{i}}
  w_{i}(\mathbf{x}_{ij})
  \left(
    U_{i}^{h}(\mathbf{x}_{ij})-U_{i\rightarrow ij}
  \right)^{2}
  \\
  &+
  \frac{1}{2}
  \sum_{\mathbf{x}_{ij}\in\mathcal{X}_{i}^{u}}
  w_{i}(\mathbf{x}_{ij})r_{s}^{2}
  \left(
    \mathbf{n}_{j}\cdot\nabla U_{i}^{h}(\mathbf{x}_{ij})
    -\frac{1}{2}
    (\mathbf{n}_{j}\mathbf{e}_{ij}:\nabla\nabla)
    U_{i}^{h}(\mathbf{x}_{ij})
    -2q_{j}^{U}
  \right)^{2}.
\end{split}
\label{eq_mls_U_with_bc}
\end{equation}
\paragraph{Reconstruction of $\delta P_i^h$}
At a prescribed-velocity boundary, the homogeneous normal pressure-gradient condition $q_{j}^{\delta P}=(\mathbf{n}\cdot\nabla p)_{j}=0$ is imposed through
\begin{equation}
\begin{split}
  \argmin_{\mathbf{c}_{i}^{\delta P}}
  ~&\frac{1}{2}
  \sum_{\mathbf{x}_{ij}\in\tilde{\mathcal{X}}^{u}_{i}}
  w_{i}(\mathbf{x}_{ij})
  \left(
    \delta P_{i}^{h}(\mathbf{x}_{ij})
    -\delta P_{i\rightarrow ij}
  \right)^{2}
  \\
  &+
  \frac{1}{2}
  \sum_{\mathbf{x}_{ij}\in\mathcal{X}_{i}^{u}}
  w_{i}(\mathbf{x}_{ij})r_{s}^{2}
  \left(
    \mathbf{n}_{j}\cdot\nabla \delta P_{i}^{h}(\mathbf{x}_{ij})
    -\frac{1}{2}
    (\mathbf{n}_{j}\mathbf{e}_{ij}:\nabla\nabla)
    \delta P_{i}^{h}(\mathbf{x}_{ij})
    -2q_{j}^{\delta P}
  \right)^{2}.
\end{split}
\label{eq_mls_deltaP_with_bc}
\end{equation}
Note that $\delta P_{i\rightarrow ij}=(p_d)_{j}-p_{i}$ for $\mathbf{x}_{ij}\in\mathcal{X}_{i}^{p}$.

The MLS approximations of $\mathbf{A}_i^h$ and $\mathbf{D}_i^h$ are constructed componentwise, where $\alpha=1,\ldots,d$ denotes the spatial-component index.
\paragraph{Reconstruction of $\mathbf{A}_i^h$}
At a prescribed-velocity boundary, the $\alpha$-th component of the normal advective flux is given by $q_{j,\alpha}^{A}=(q^{U}u_{d,\alpha})_{j}$, as defined in Eq.~\eqref{eq_bc_q_A_velocity}.
The MLS approximation $A_{i,\alpha}^{h}$ is then constructed independently for each spatial component $\alpha=1,\ldots,d$ through
\begin{equation}
\begin{split}
  \argmin_{\mathbf{c}_{i}^{A,\alpha}}
  ~&\frac{1}{2}
  \sum_{\mathbf{x}_{ij}\in\tilde{\mathcal{X}}^{u}_{i}}
  w_{i}(\mathbf{x}_{ij})
  \left(
    A_{i,\alpha}^{h}(\mathbf{x}_{ij})
    -A_{i\rightarrow ij,\alpha}
  \right)^{2}
  \\
  &+
  \frac{1}{2}
  \sum_{\mathbf{x}_{ij}\in\mathcal{X}_{i}^{u}}
  w_{i}(\mathbf{x}_{ij})r_{s}^{2}
  \left(
    \mathbf{n}_{j}\cdot\nabla A_{i,\alpha}^{h}(\mathbf{x}_{ij})
    -\frac{1}{2}
    (\mathbf{n}_{j}\mathbf{e}_{ij}:\nabla\nabla)
    A_{i,\alpha}^{h}(\mathbf{x}_{ij})
    -2q_{j,\alpha}^{A}
  \right)^{2}.
\end{split}
\label{eq_mls_A_with_bc}
\end{equation}
Note that $A_{i\rightarrow ij,\alpha} = U_{i\rightarrow ij} \left.u_{i,\alpha}^{h}\right|_{ij}$ for $\mathbf{x}_{ij}\in\mathcal{X}_{i}^{p}$, according to the boundary treatment in Eq. \eqref{eq_bc_A_pressure}.

\paragraph{Reconstruction of $\mathbf{D}_i^h$}
At a prescribed-pressure boundary, the homogeneous normal velocity-gradient condition for the $\alpha$-th component is given by $q_{j,\alpha}^{D}=(\mathbf{n}\cdot\nabla u_{\alpha})_{j}=0$, as defined in Eq.~\eqref{eq_bc_q_D_pressure}.
The MLS approximation $D_{i,\alpha}^{h}$ is then constructed independently for each spatial component $\alpha=1,\ldots,d$ through
\begin{equation}
\begin{split}
  \argmin_{\mathbf{c}_{i}^{D,\alpha}}
  ~&\frac{1}{2}
  \sum_{\mathbf{x}_{ij}\in\tilde{\mathcal{X}}^{p}_{i}}
  w_{i}(\mathbf{x}_{ij})
  \left(
    D_{i,\alpha}^{h}(\mathbf{x}_{ij})
    -D_{i\rightarrow ij,\alpha}
  \right)^{2}
  \\
  &+
  \frac{1}{2}
  \sum_{\mathbf{x}_{ij}\in\mathcal{X}_{i}^{p}}
  w_{i}(\mathbf{x}_{ij})r_{s}^{2}
  \left(
    \mathbf{n}_{j}\cdot\nabla D_{i,\alpha}^{h}(\mathbf{x}_{ij})
    -\frac{1}{2}
    (\mathbf{n}_{j}\mathbf{e}_{ij}:\nabla\nabla)
    D_{i,\alpha}^{h}(\mathbf{x}_{ij})
    -2q_{j,\alpha}^{D}
  \right)^{2}.
\end{split}
\label{eq_mls_D_with_bc}
\end{equation}
Note that $D_{i\rightarrow ij,\alpha} = (u_{d,\alpha})_{j}-u_{i,\alpha}$ for $\mathbf{x}_{ij}\in\mathcal{X}_{i}^{u}$, according to the boundary treatment in Eq. \eqref{eq_bc_D_velocity}.

\section{Relationship between the discrete velocity divergence and the residual of the pressure Poisson equation} \label{appendix_div_u_residual}
\setcounter{equation}{0}
\setcounter{figure}{0}

When the linear system of the pressure Poisson equation \eqref{eq_pressure_2} is solved using an iterative method, the absolute residual vector $\boldsymbol{\varepsilon}$ of the linear system at the $m$-th iteration can be written as
\begin{equation}
\begin{split}
\boldsymbol{\varepsilon}^{m}
&= \frac{1}{4} \left(\frac{\rho}{\Delta t} {\rm L} {U}^{*} - {\rm L} \delta{P}^{n+1,\,m}\right) \\
&= \frac{1}{4} \left(\frac{\rho}{\Delta t} {\rm L} {U}^{*} - {\rm L}{\rm G} p^{n+1,\,m}\right).
\end{split}
\end{equation}
Here, ${\rm L}$ denotes the discrete Laplacian operator obtained by differentiating the MLS polynomial, and ${\rm G}$ denotes the discrete gradient-type operator that maps $p$ to $\delta P$ according to Eq.~\eqref{eq_deltaP_ij}. The same discrete Laplacian operator ${\rm L}$ can be applied to both the velocity and pressure terms because the edge-based quantities $U$ and $\delta P$ are reconstructed using the same MLS procedure, including the boundary treatment at prescribed-velocity boundaries, as defined in Eqs.~\eqref{eq_mls_U_with_bc} and \eqref{eq_mls_deltaP_with_bc}. When convergence is reached at the $m_{\rm crit}$-th iteration, the corresponding iterative solution $p^{n+1,\,m_{\rm crit}}$ is taken as $p^{n+1}$. From Eq. \eqref{eq_projection_stg_2}, ${U}^{n+1}$ can be expressed as
\begin{align}
{U}^{n+1} = {U}^{*} - \frac{\Delta t}{\rho} {\rm G} p^{n+1}.
\end{align}
From these equations and Eq. \eqref{eq_div_uh_at_xi}, the discrete velocity divergence at the $(n+1)$-th time step, ${\rm D}\mathbf{u}^{n+1} \coloneqq \nabla\cdot\mathbf{u}^{h} |^{n+1} = \nabla^{2}{U}^{h} |^{n+1}/4$, can be expressed as
\begin{equation}
\begin{split}
{\rm D}\mathbf{u}^{n+1}
&= \frac{1}{4}{\rm L}{U}^{n+1} \\
&= \frac{1}{4} \left({\rm L}{U}^{*} - \frac{\Delta t}{\rho} {\rm LG} p^{n+1}\right) \\
&= \frac{\Delta t}{\rho} \boldsymbol{\varepsilon}^{m_{\rm crit}}.
\end{split}
\end{equation}
This relation indicates that the discrete velocity divergence obtained using the radial velocity component corresponds directly to the residual of the pressure Poisson equation through the pressure-projection procedure.

The absolute $L_{2}$ residual norm is defined as $\varepsilon_{2}^{m}=\|\boldsymbol{\varepsilon}^{m}\|_{2}$. Accordingly, the norm of the discrete velocity divergence at the $(n+1)$-th time step satisfies
\begin{align}
\|{\rm D} \mathbf{u}^{n+1}\|_{2}
= \frac{\Delta t}{\rho} \varepsilon_{2}^{m_{\rm crit}}.
\label{eq_divergence_residual_relation}
\end{align}
In the simulations presented in this study, convergence of the iterative solver was assessed using the relative residual norm $\varepsilon_{2}^{m}/\varepsilon_{2}^{0}$. The iteration was terminated at the first iteration $m_{\rm crit}$ satisfying
\begin{equation}
 \frac{\varepsilon_{2}^{m_{\rm crit}}} {\varepsilon_{2}^{0}}
 \leq \varepsilon_{\rm crit}.
 \label{eq_relative_residual_criterion}
\end{equation}
Thus, $\varepsilon_{2}^{m_{\rm crit}}$ denotes the absolute residual norm at the stopping iteration.

\begin{figure}
  \begin{center}
    \includegraphics[width=0.6\linewidth]{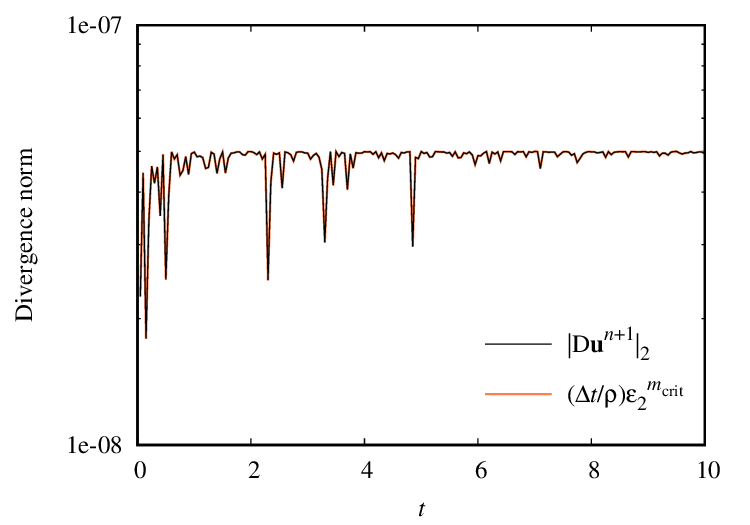}
    \caption{%
      Time evolution of $\|{\rm D}\mathbf{u}^{n+1}\|_{2}$ and $(\Delta t/\rho)\varepsilon_{2}^{m_{\rm crit}}$ for the stationary-cylinder flow problem in Section~\ref{subsec_stat_cyl}, using the staggered formulation with $\beta=0.99$, $\Delta t=\Delta t^{*}$ and $\varepsilon_{\rm crit}=10^{-5}$.
    }
    \label{fig_div_vel_residual}
  \end{center}
\end{figure}

Figure~\ref{fig_div_vel_residual} shows the time evolution of $\|{\rm D}\mathbf{u}^{n+1}\|_{2}$ and $(\Delta t/\rho)\varepsilon_{2}^{m_{\rm crit}}$ for the stationary-cylinder flow problem in Section~\ref{subsec_stat_cyl}, using the staggered formulation with $\beta=0.99$ and $\Delta t=\Delta t^{*}$. As predicted by Eq.~\eqref{eq_divergence_residual_relation}, the two quantities are observed to coincide to numerical precision. 

\end{appendix}

\bibliographystyle{elsarticle-num}  
\bibliography{ref}

\begin{thebibliography}{10}
\expandafter\ifx\csname url\endcsname\relax
  \def\url#1{\texttt{#1}}\fi
\expandafter\ifx\csname urlprefix\endcsname\relax\def\urlprefix{URL }\fi
\expandafter\ifx\csname href\endcsname\relax
  \def\href#1#2{#2} \def\path#1{#1}\fi

\bibitem{belytschko_ElementfreeGalerkinMethods_1994}
T.~Belytschko, Y.~Y. Lu, L.~Gu, Element-free {{Galerkin}} methods, Int. J. Numer. Methods Eng. 37~(2) (1994) 229--256.
\newblock \href {https://doi.org/10.1002/nme.1620370205} {\path{doi:10.1002/nme.1620370205}}.

\bibitem{wang_SuperiorityMixedElement_2012}
X.~Wang, J.~Ouyang, J.~Su, B.~Yang, On the superiority of the mixed element free {{Galerkin}} method for solving the steady incompressible flow problems, Eng. Anal. Bound. Elem. 36~(11) (2012) 1618--1630.
\newblock \href {https://doi.org/10.1016/j.enganabound.2012.05.006} {\path{doi:10.1016/j.enganabound.2012.05.006}}.

\bibitem{atluri_NewMeshlessLocal_1998}
S.~N. Atluri, T.~Zhu, A new {{Meshless Local Petrov-Galerkin}} ({{MLPG}}) approach in computational mechanics, Comput. Mech. 22~(2) (1998) 117--127.
\newblock \href {https://doi.org/10.1007/s004660050346} {\path{doi:10.1007/s004660050346}}.

\bibitem{abbaszadeh_DirectMeshlessLocal_2020}
M.~Abbaszadeh, M.~Dehghan, Direct meshless local {{Petrov}}--{{Galerkin}} ({{DMLPG}}) method for time-fractional fourth-order reaction--diffusion problem on complex domains, Comput. Math. Appl. 79~(3) (2020) 876--888.
\newblock \href {https://doi.org/10.1016/j.camwa.2019.08.001} {\path{doi:10.1016/j.camwa.2019.08.001}}.

\bibitem{liu_ReproducingKernelParticle_1995}
W.~K. Liu, S.~Jun, Y.~F. Zhang, Reproducing kernel particle methods, Int. J. Numer. Methods Fluids 20~(8-9) (1995) 1081--1106.
\newblock \href {https://doi.org/10.1002/fld.1650200824} {\path{doi:10.1002/fld.1650200824}}.

\bibitem{gingold_SmoothedParticleHydrodynamics_1977}
R.~A. Gingold, J.~J. Monaghan, Smoothed particle hydrodynamics: Theory and application to non-spherical stars, Mon. Not. R. Astron. Soc. 181~(3) (1977) 375--389.
\newblock \href {https://doi.org/10.1093/mnras/181.3.375} {\path{doi:10.1093/mnras/181.3.375}}.

\bibitem{koshizuka_MovingParticleSemiImplicitMethod_1996}
S.~Koshizuka, Y.~Oka, Moving-{{Particle Semi-Implicit Method}} for {{Fragmentation}} of {{Incompressible Fluid}}, Nucl. Sci. Eng. 123~(3) (1996) 421--434.
\newblock \href {https://doi.org/10.13182/NSE96-A24205} {\path{doi:10.13182/NSE96-A24205}}.

\bibitem{onate_StabilizedFinitePoint_1996}
E.~O{\~n}ate, S.~Idelsohn, O.~C. Zienkiewicz, R.~L. Taylor, C.~Sacco, A stabilized finite point method for analysis of fluid mechanics problems, Comput. Methods Appl. Mech. Eng. 139~(1) (1996) 315--346.
\newblock \href {https://doi.org/10.1016/S0045-7825(96)01088-2} {\path{doi:10.1016/S0045-7825(96)01088-2}}.

\bibitem{perrone_GeneralFiniteDifference_1975}
N.~Perrone, R.~Kao, A general finite difference method for arbitrary meshes, Comput. Struct. 5~(1) (1975) 45--57.
\newblock \href {https://doi.org/10.1016/0045-7949(75)90018-8} {\path{doi:10.1016/0045-7949(75)90018-8}}.

\bibitem{liszka_FiniteDifferenceMethod_1980}
T.~Liszka, J.~Orkisz, The finite difference method at arbitrary irregular grids and its application in applied mechanics, Comput. Struct. 11~(1) (1980) 83--95.
\newblock \href {https://doi.org/10.1016/0045-7949(80)90149-2} {\path{doi:10.1016/0045-7949(80)90149-2}}.

\bibitem{benito_SolvingParabolicHyperbolic_2007}
J.~J. Benito, F.~Ure{\~n}a, L.~Gavete, Solving parabolic and hyperbolic equations by the generalized finite difference method, J. Comput. Appl. Math. 209~(2) (2007) 208--233.
\newblock \href {https://doi.org/10.1016/j.cam.2006.10.090} {\path{doi:10.1016/j.cam.2006.10.090}}.

\bibitem{prieto_ApplicationGeneralizedFinite_2011}
F.~U. Prieto, J.~J. Benito~Mu{\~n}oz, L.~G. Corvinos, Application of the generalized finite difference method to solve the advection--diffusion equation, J. Comput. Appl. Math. 235~(7) (2011) 1849--1855.
\newblock \href {https://doi.org/10.1016/j.cam.2010.05.026} {\path{doi:10.1016/j.cam.2010.05.026}}.

\bibitem{tang_LeastsquaresGeneralizedFinite_2025}
Z.~Tang, H.~Pan, Z.~Fu, M.~Chen, L.~Ling, A least-squares generalized finite difference method for solving nonlinear reaction--diffusion systems, Eng. Anal. Bound. Elem. 179 (2025) 106351.
\newblock \href {https://doi.org/10.1016/j.enganabound.2025.106351} {\path{doi:10.1016/j.enganabound.2025.106351}}.

\bibitem{dilts_MovingleastsquaresparticleHydrodynamicsConsistency_1999}
G.~A. Dilts, Moving-least-squares-particle hydrodynamics---{{I}}. {{Consistency}} and stability, Int. J. Numer. Methods Eng. 44~(8) (1999) 1115--1155.
\newblock \href {https://doi.org/10.1002/(SICI)1097-0207(19990320)44:8<1115::AID-NME547>3.0.CO;2-L} {\path{doi:10.1002/(SICI)1097-0207(19990320)44:8<1115::AID-NME547>3.0.CO;2-L}}.

\bibitem{suchde_MeshfreeGFDMSolvers_2018}
P.~Suchde, J.~Kuhnert, S.~Tiwari, On meshfree {{GFDM}} solvers for the incompressible {{Navier}}--{{Stokes}} equations, Comput. Fluids 165 (2018) 1--12.
\newblock \href {https://doi.org/10.1016/j.compfluid.2018.01.008} {\path{doi:10.1016/j.compfluid.2018.01.008}}.

\bibitem{tamai_LeastSquaresMoving_2014}
T.~Tamai, S.~Koshizuka, Least squares moving particle semi-implicit method, Comput. Part. Mech. 1~(3) (2014) 277--305.
\newblock \href {https://doi.org/10.1007/s40571-014-0027-2} {\path{doi:10.1007/s40571-014-0027-2}}.

\bibitem{matsunaga_ImprovedTreatmentWall_2020}
T.~Matsunaga, A.~S{\"o}dersten, K.~Shibata, S.~Koshizuka, Improved treatment of wall boundary conditions for a particle method with consistent spatial discretization, Comput. Methods Appl. Mech. Eng. 358 (2020) 112624.
\newblock \href {https://doi.org/10.1016/j.cma.2019.112624} {\path{doi:10.1016/j.cma.2019.112624}}.

\bibitem{vasyliv_SimulatingIncompressibleFlow_2020}
Y.~Vasyliv, A.~Alexeev, Simulating incompressible flow on moving meshfree grids, Comput. Fluids 200 (2020) 104464.
\newblock \href {https://doi.org/10.1016/j.compfluid.2020.104464} {\path{doi:10.1016/j.compfluid.2020.104464}}.

\bibitem{matsuda_ParticlebasedMethodUsing_2022}
T.~Matsuda, K.~Tsukui, S.~Ii, A particle-based method using the mesh-constrained discrete point approach for two-dimensional {{Stokes}} flows, Mech. Eng. J. 9~(5) (2022) 22--00204.
\newblock \href {https://doi.org/10.1299/mej.22-00204} {\path{doi:10.1299/mej.22-00204}}.

\bibitem{matsuda_MeshconstrainedDiscretePoint_2025}
T.~Matsuda, S.~Ii, A mesh-constrained discrete point method for incompressible flows with moving boundaries, J. Comput. Phys. 532 (2025) 113945.
\newblock \href {https://doi.org/10.1016/j.jcp.2025.113945} {\path{doi:10.1016/j.jcp.2025.113945}}.

\bibitem{lancaster_SurfacesGeneratedMoving_1981}
P.~Lancaster, K.~Salkauskas, Surfaces generated by moving least squares methods, Math. Comput. 37~(155) (1981) 141--158.
\newblock \href {https://doi.org/10.1090/S0025-5718-1981-0616367-1} {\path{doi:10.1090/S0025-5718-1981-0616367-1}}.

\bibitem{matsunaga_StabilizedLSMPSMethod_2022}
T.~Matsunaga, S.~Koshizuka, Stabilized {{LSMPS}} method for complex free-surface flow simulation, Comput. Methods Appl. Mech. Eng. 389 (2022) 114416.
\newblock \href {https://doi.org/10.1016/j.cma.2021.114416} {\path{doi:10.1016/j.cma.2021.114416}}.

\bibitem{matsunaga_HighorderTimemarchingSchemes_2025}
T.~Matsunaga, High-order time-marching schemes for incompressible flow in particle methods, Comput. Methods Appl. Mech. Eng. 447 (2025) 118395.
\newblock \href {https://doi.org/10.1016/j.cma.2025.118395} {\path{doi:10.1016/j.cma.2025.118395}}.

\bibitem{hopman_QuantifyingCheckerboardProblem_2025}
J.~A. Hopman, D.~Santos, {\`A}.~{Alsalti-Baldellou}, J.~Rigola, F.~X. Trias, Quantifying the checkerboard problem to reduce numerical dissipation, J. Comput. Phys. 521 (2025) 113537.
\newblock \href {https://doi.org/10.1016/j.jcp.2024.113537} {\path{doi:10.1016/j.jcp.2024.113537}}.

\bibitem{swegle_AnalysisSmoothedParticle_1994}
J.~W. Swegle, S.~W. Attaway, M.~W. Heinstein, F.~J. Mello, D.~L. Hicks, An analysis of smoothed particle hydrodynamics, Tech. Rep. SAND--93-2513, Sandia National Labs., Albuquerque, NM (United States) (Feb. 1994).
\newblock \href {https://doi.org/10.2172/10159839} {\path{doi:10.2172/10159839}}.

\bibitem{liu_DualParticleApproachIncompressible_2024}
S.~Liu, X.~He, Y.~Guo, Y.~Chang, W.~Wang, A {{Dual-Particle Approach}} for {{Incompressible SPH Fluids}}, ACM Trans. Graph. 43~(3) (2024) 28:1--28:18.
\newblock \href {https://doi.org/10.1145/3649888} {\path{doi:10.1145/3649888}}.

\bibitem{SHAHANE2021110623}
S.~Shahane, A.~Radhakrishnan, S.~P. Vanka, A high-order accurate meshless method for solution of incompressible fluid flow problems, Journal of Computational Physics 445 (2021) 110623.
\newblock \href {https://doi.org/10.1016/j.jcp.2021.110623} {\path{doi:10.1016/j.jcp.2021.110623}}.

\bibitem{Unnikrishnan2022}
A.~Unnikrishnan, S.~Shahane, V.~Narayanan, S.~P. Vanka, Shear-driven flow in an elliptical enclosure generated by an inner rotating circular cylinder, Physics of Fluids 34~(1) (2022) 013607.
\newblock \href {https://doi.org/10.1063/5.0076537} {\path{doi:10.1063/5.0076537}}.

\bibitem{CHU2024112822}
T.~Chu, O.~T. Schmidt, Mesh-free hydrodynamic stability, Journal of Computational Physics 502 (2024) 112822.
\newblock \href {https://doi.org/10.1016/j.jcp.2024.112822} {\path{doi:10.1016/j.jcp.2024.112822}}.

\bibitem{Unnikrishnan2024}
A.~Unnikrishnan, V.~Narayanan, L.~P. Chamorro, S.~P. Vanka, Taylor-{C}ouette flow and heat transfer in an elliptical enclosure with a rotating inner cylinder, Physics of Fluids 36~(3) (2024) 033604.
\newblock \href {https://doi.org/10.1063/5.0190826} {\path{doi:10.1063/5.0190826}}.

\bibitem{dyka_ApproachTensionInstability_1995}
C.~T. Dyka, R.~P. Ingel, An approach for tension instability in smoothed particle hydrodynamics ({{SPH}}), Comput. Struct. 57~(4) (1995) 573--580.
\newblock \href {https://doi.org/10.1016/0045-7949(95)00059-P} {\path{doi:10.1016/0045-7949(95)00059-P}}.

\bibitem{randles_NormalizedSPHStress_2000}
P.~W. Randles, L.~D. Libersky, Normalized {{SPH}} with stress points, Int. J. Numer. Methods Eng. 48~(10) (2000) 1445--1462.
\newblock \href {https://doi.org/10.1002/1097-0207(20000810)48:10<1445::AID-NME831>3.0.CO;2-9} {\path{doi:10.1002/1097-0207(20000810)48:10<1445::AID-NME831>3.0.CO;2-9}}.

\bibitem{chalk_StressParticleSmoothedParticle_2020}
C.~M. Chalk, M.~Pastor, J.~Peakall, D.~J. Borman, P.~A. Sleigh, W.~Murphy, R.~Fuentes, Stress-{{Particle Smoothed Particle Hydrodynamics}}: {{An}} application to the failure and post-failure behaviour of slopes, Comput. Methods Appl. Mech. Eng. 366 (2020) 113034.
\newblock \href {https://doi.org/10.1016/j.cma.2020.113034} {\path{doi:10.1016/j.cma.2020.113034}}.

\bibitem{he_StaggeredMeshlessSolidfluid_2012}
X.~He, N.~Liu, G.~Wang, F.~Zhang, S.~Li, S.~Shao, H.~Wang, Staggered meshless solid-fluid coupling, ACM Trans. Graph. 31~(6) (2012) 149:1--149:12.
\newblock \href {https://doi.org/10.1145/2366145.2366168} {\path{doi:10.1145/2366145.2366168}}.

\bibitem{he_VariationalStaggeredParticle_2020}
X.~He, H.~Wang, G.~Wang, H.~Wang, E.~Wu, A {{Variational Staggered Particle Framework}} for {{Incompressible Free-Surface Flows}} (Jan. 2020).
\newblock \href {http://arxiv.org/abs/2001.09421} {\path{arXiv:2001.09421}}, \href {https://doi.org/10.48550/arXiv.2001.09421} {\path{doi:10.48550/arXiv.2001.09421}}.

\bibitem{park_ExistenceStabilityVirtual_2016}
S.-K. Park, G.~Jo, H.~J. Choe, Existence and stability in the virtual interpolation point method for the {{Stokes}} equations, J. Comput. Phys. 307 (2016) 535--549.
\newblock \href {https://doi.org/10.1016/j.jcp.2015.12.002} {\path{doi:10.1016/j.jcp.2015.12.002}}.

\bibitem{CHU2023111756}
T.~Chu, O.~T. Schmidt, {RBF}-{FD} discretization of the {N}avier-{S}tokes equations on scattered but staggered nodes, Journal of Computational Physics 474 (2023) 111756.
\newblock \href {https://doi.org/10.1016/j.jcp.2022.111756} {\path{doi:10.1016/j.jcp.2022.111756}}.

\bibitem{trask_HighorderStaggeredMeshless_2017}
N.~Trask, M.~Perego, P.~Bochev, A high-order staggered meshless method for elliptic problems, SIAM J. Sci. Comput. 39~(2) (2017) A479--A502.
\newblock \href {https://doi.org/10.1137/16M1055992} {\path{doi:10.1137/16M1055992}}.

\bibitem{trask_CompatibleHighorderMeshless_2018}
N.~Trask, M.~Maxey, X.~Hu, A compatible high-order meshless method for the {{Stokes}} equations with applications to suspension flows, J. Comput. Phys. 355 (2018) 310--326.
\newblock \href {https://doi.org/10.1016/j.jcp.2017.10.039} {\path{doi:10.1016/j.jcp.2017.10.039}}.

\bibitem{chorin_NumericalSolutionNavierStokes_1968}
A.~J. Chorin, Numerical solution of the {{Navier-Stokes}} equations, Math. Comput. 22~(104) (1968) 745--762.
\newblock \href {https://doi.org/10.1090/S0025-5718-1968-0242392-2} {\path{doi:10.1090/S0025-5718-1968-0242392-2}}.

\bibitem{xiao_SimpleCIPFinite_2004}
F.~Xiao, A {{Simple CIP Finite Volume Method}} for {{Incompressible Flows}}, JSME Int. J. Ser. B Fluids Thermal Eng. 47~(4) (2004) 664--671.
\newblock \href {https://doi.org/10.1299/jsmeb.47.664} {\path{doi:10.1299/jsmeb.47.664}}.

\bibitem{xiao_NumericalSimulationsFreeinterface_2005}
F.~Xiao, A.~Ikebata, T.~Hasegawa, Numerical simulations of free-interface fluids by a multi-integrated moment method, Computers \& Structures 83~(6) (2005) 409--423.
\newblock \href {https://doi.org/10.1016/j.compstruc.2004.06.005} {\path{doi:10.1016/j.compstruc.2004.06.005}}.

\bibitem{xiao_UnifiedFormulationCompressible_2006}
F.~Xiao, R.~Akoh, S.~Ii, Unified formulation for compressible and incompressible flows by using multi-integrated moments {{II}}: {{Multi-dimensional}} version for compressible and incompressible flows, Journal of Computational Physics 213~(1) (2006) 31--56.
\newblock \href {https://doi.org/10.1016/j.jcp.2005.08.002} {\path{doi:10.1016/j.jcp.2005.08.002}}.

\bibitem{rhie_NumericalStudyTurbulent_1983}
C.~M. Rhie, W.~L. Chow, Numerical study of the turbulent flow past an airfoil with trailing edge separation, AIAA J. 21~(11) (1983) 1525--1532.
\newblock \href {https://doi.org/10.2514/3.8284} {\path{doi:10.2514/3.8284}}.

\bibitem{zhang_GeneralizedFormulationsRhie_2014}
S.~Zhang, X.~Zhao, S.~Bayyuk, Generalized formulations for the {{Rhie}}--{{Chow}} interpolation, J. Comput. Phys. 258 (2014) 880--914.
\newblock \href {https://doi.org/10.1016/j.jcp.2013.11.006} {\path{doi:10.1016/j.jcp.2013.11.006}}.

\bibitem{harlow_ParticleincellMethodNumerical_1962}
F.~H. Harlow, The particle-in-cell method for numerical solution of problems in fluid dynamics, Tech. Rep. LADC-5288, Los Alamos Scientific Lab., N. Mex. (Feb. 1962).
\newblock \href {https://doi.org/10.2172/4769185} {\path{doi:10.2172/4769185}}.

\bibitem{brackbill_FLIPMethodAdaptively_1986}
J.~U. Brackbill, H.~M. Ruppel, {{FLIP}}: {{A}} method for adaptively zoned, particle-in-cell calculations of fluid flows in two dimensions, J. Comput. Phys. 65~(2) (1986) 314--343.
\newblock \href {https://doi.org/10.1016/0021-9991(86)90211-1} {\path{doi:10.1016/0021-9991(86)90211-1}}.

\bibitem{zhang_LeastsquaresMeshfreeMethod_2005}
X.~K. Zhang, K.-C. Kwon, S.-K. Youn, The least-squares meshfree method for the steady incompressible viscous flow, J. Comput. Phys. 206~(1) (2005) 182--207.
\newblock \href {https://doi.org/10.1016/j.jcp.2004.11.033} {\path{doi:10.1016/j.jcp.2004.11.033}}.

\bibitem{chandra_StabilizedMixedMaterial_2024}
B.~Chandra, R.~Hashimoto, S.~Matsumi, K.~Kamrin, K.~Soga, Stabilized mixed material point method for incompressible fluid flow analysis, Comput. Methods Appl. Mech. Eng. 419 (2024) 116644.
\newblock \href {https://doi.org/10.1016/j.cma.2023.116644} {\path{doi:10.1016/j.cma.2023.116644}}.

\bibitem{he_ImprovedMPMFormulation_2025}
K.-Y. He, Y.-F. Jin, X.-W. Zhou, Z.-Y. Yin, X.~Chen, An improved {{MPM}} formulation for free surface flow problems based on finite volume method, Comput. Methods Appl. Mech. Eng. 446 (2025) 118264.
\newblock \href {https://doi.org/10.1016/j.cma.2025.118264} {\path{doi:10.1016/j.cma.2025.118264}}.

\bibitem{xu_CoupledSPHFVM_2021}
Y.~Xu, G.~Yang, Y.~Zhu, D.~Hu, A coupled {{SPH}}--{{FVM}} method for simulating incompressible interfacial flows with large density difference, Eng. Anal. Bound. Elem. 128 (2021) 227--243.
\newblock \href {https://doi.org/10.1016/j.enganabound.2021.04.005} {\path{doi:10.1016/j.enganabound.2021.04.005}}.

\bibitem{nishiguchi_EulerianFiniteVolume_2024}
K.~Nishiguchi, T.~Shimada, C.~Peco, K.~Kondo, S.~Okazawa, M.~Tsubokura, Eulerian finite volume method using {{Lagrangian}} markers with reference map for incompressible fluid--structure interaction problems, Comput. Fluids 274 (2024) 106210.
\newblock \href {https://doi.org/10.1016/j.compfluid.2024.106210} {\path{doi:10.1016/j.compfluid.2024.106210}}.

\bibitem{katz_meshless_2009}
A.~Katz, A.~Jameson, \href{https://arc.aiaa.org/doi/10.2514/6.2009-3534}{A {Meshless} {Volume} {Scheme}}, in: 19th {AIAA} {Computational} {Fluid} {Dynamics}, Fluid {Dynamics} and {Co}-located {Conferences}, American Institute of Aeronautics and Astronautics, 2009.
\newblock \href {https://doi.org/10.2514/6.2009-3534} {\path{doi:10.2514/6.2009-3534}}.
\newline\urlprefix\url{https://arc.aiaa.org/doi/10.2514/6.2009-3534}

\bibitem{kwan-yu_chiu_conservative_2012}
E.~Kwan-yu Chiu, Q.~Wang, R.~Hu, A.~Jameson, \href{https://doi.org/10.1137/110842740}{A {Conservative} {Mesh}-{Free} {Scheme} and {Generalized} {Framework} for {Conservation} {Laws}}, SIAM J. Sci. Comput. 34~(6) (2012) A2896--A2916.
\newblock \href {https://doi.org/10.1137/110842740} {\path{doi:10.1137/110842740}}.
\newline\urlprefix\url{https://doi.org/10.1137/110842740}

\bibitem{suchde_flux_2017}
P.~Suchde, J.~Kuhnert, S.~Schröder, A.~Klar, \href{https://onlinelibrary.wiley.com/doi/abs/10.1002/nme.5511}{A flux conserving meshfree method for conservation laws}, Int. J. Numer. Methods Eng. 112~(3) (2017) 238--256.
\newblock \href {https://doi.org/10.1002/nme.5511} {\path{doi:10.1002/nme.5511}}.
\newline\urlprefix\url{https://onlinelibrary.wiley.com/doi/abs/10.1002/nme.5511}

\bibitem{TRASK2020109187}
N.~Trask, P.~Bochev, M.~Perego, \href{https://www.sciencedirect.com/science/article/pii/S0021999119308927}{A conservative, consistent, and scalable meshfree mimetic method}, Journal of Computational Physics 409 (2020) 109187.
\newblock \href {https://doi.org/10.1016/j.jcp.2019.109187} {\path{doi:10.1016/j.jcp.2019.109187}}.
\newline\urlprefix\url{https://www.sciencedirect.com/science/article/pii/S0021999119308927}

\bibitem{vandervorst_BiCGSTABFastSmoothly_1992}
H.~A. {van der Vorst}, Bi-{{CGSTAB}}: {{A Fast}} and {{Smoothly Converging Variant}} of {{Bi-CG}} for the {{Solution}} of {{Nonsymmetric Linear Systems}}, SIAM J. Sci. and Stat. Comput. 13~(2) (1992) 631--644.
\newblock \href {https://doi.org/10.1137/0913035} {\path{doi:10.1137/0913035}}.

\bibitem{liu_PreconditionedMultigridMethods_1998}
C.~Liu, X.~Zheng, C.~H. Sung, Preconditioned multigrid methods for unsteady incompressible flows, J. Comput. Phys. 139~(1) (1998) 35--57.
\newblock \href {https://doi.org/10.1006/jcph.1997.5859} {\path{doi:10.1006/jcph.1997.5859}}.

\bibitem{ding_NumericalSimulationFlows_2007}
H.~Ding, C.~Shu, K.~S. Yeo, D.~Xu, Numerical simulation of flows around two circular cylinders by mesh-free least square-based finite difference methods, Int. J. Numer. Methods Fluids 53~(2) (2007) 305--332.
\newblock \href {https://doi.org/10.1002/fld.1281} {\path{doi:10.1002/fld.1281}}.

\bibitem{prandtl_ErgebnisseAerodynamischenVersuchsanstalt_2009}
L.~Prandtl, C.~Wieselsberger, A.~Betz (Eds.), {Ergebnisse der aerodynamischen versuchsanstalt zu g\"ottingen - II. Lieferung}, Vol.~5 of {G\"ottinger Klassiker der Str\"omungsmechanik}, Universit\"atsverlag G\"ottingen, G\"ottingen, 2009.
\newblock \href {https://doi.org/10.17875/gup2009-101} {\path{doi:10.17875/gup2009-101}}.

\bibitem{williamson_ObliqueParallelModes_1989}
C.~H.~K. Williamson, Oblique and parallel modes of vortex shedding in the wake of a circular cylinder at low reynolds numbers, J. Fluid Mech. 206 (1989) 579--627.
\newblock \href {https://doi.org/10.1017/S0022112089002429} {\path{doi:10.1017/S0022112089002429}}.

\bibitem{ghia_HighReSolutionsIncompressible_1982}
U.~Ghia, K.~N. Ghia, C.~T. Shin, High-{{Re}} solutions for incompressible flow using the {{Navier-Stokes}} equations and a multigrid method, J. Comput. Phys. 48~(3) (1982) 387--411.
\newblock \href {https://doi.org/10.1016/0021-9991(82)90058-4} {\path{doi:10.1016/0021-9991(82)90058-4}}.

\bibitem{hou_SimulationCavityFlow_1995}
S.~Hou, Q.~Zou, S.~Chen, G.~Doolen, A.~C. Cogley, Simulation of {{Cavity Flow}} by the {{Lattice Boltzmann Method}}, J. Comput. Phys. 118~(2) (1995) 329--347.
\newblock \href {https://doi.org/10.1006/jcph.1995.1103} {\path{doi:10.1006/jcph.1995.1103}}.

\bibitem{peng_Transition2DLiddriven_2003}
Y.-F. Peng, Y.-H. Shiau, R.~R. Hwang, Transition in a 2-{{D}} lid-driven cavity flow, Comput. Fluids 32~(3) (2003) 337--352.
\newblock \href {https://doi.org/10.1016/S0045-7930(01)00053-6} {\path{doi:10.1016/S0045-7930(01)00053-6}}.

\bibitem{li_FlowStructuresTwodimensional_2025}
L.~Li, J.~F. Klausner, R.~Mei, Flow structures in two-dimensional lid-driven cavity flow: {{Benchmark}} numerical results for steady flows, Eur. J. Mech. B/Fluids 114 (2025) 204313.
\newblock \href {https://doi.org/10.1016/j.euromechflu.2025.204313} {\path{doi:10.1016/j.euromechflu.2025.204313}}.

\bibitem{dutsch_LowReynoldsnumberFlowOscillating_1998}
H.~D{\"u}tsch, F.~Durst, S.~Becker, H.~Lienhart, Low-{{Reynolds-number}} flow around an oscillating circular cylinder at low {{Keulegan}}--{{Carpenter}} numbers, J. Fluid Mech. 360 (1998) 249--271.
\newblock \href {https://doi.org/10.1017/S002211209800860X} {\path{doi:10.1017/S002211209800860X}}.

\bibitem{guilmineau_NUMERICALSIMULATIONVORTEX_2002}
E.~Guilmineau, P.~Queutey, A {{NUMERICAL SIMULATION OF VORTEX SHEDDING FROM AN OSCILLATING CIRCULAR CYLINDER}}, J. Fluid. Struct. 16~(6) (2002) 773--794.
\newblock \href {https://doi.org/10.1006/jfls.2002.0449} {\path{doi:10.1006/jfls.2002.0449}}.

\bibitem{chi_DirectionalGhostcellImmersed_2020}
C.~Chi, A.~Abdelsamie, D.~Th{\'e}venin, A directional ghost-cell immersed boundary method for incompressible flows, J. Comput. Phys. 404 (2020) 109122.
\newblock \href {https://doi.org/10.1016/j.jcp.2019.109122} {\path{doi:10.1016/j.jcp.2019.109122}}.

\bibitem{ghomizad_SharpInterfaceDirectforcing_2021}
M.~B. Ghomizad, H.~Kor, K.~Fukagata, A sharp interface direct-forcing immersed boundary method using the moving least square approximation, J. Fluid Sci. Technol. 16~(2) (2021) JFST0013--JFST0013.
\newblock \href {https://doi.org/10.1299/jfst.2021jfst0013} {\path{doi:10.1299/jfst.2021jfst0013}}.

\end{thebibliography}

\end{document}